\documentclass[EJP]{ejpecp} 

\usepackage[T1]{fontenc}
\usepackage[utf8]{inputenc}

\SHORTTITLE{Estimation of Local Anisotropy Based on Level Sets} 

\TITLE{Estimation of Local Anisotropy Based on Level Sets}

\AUTHORS{%
  Corinne~Berzin\footnote{Univ. Grenoble Alpes, CNRS, LJK, 38000 Grenoble, France. \BEMAIL{corinne.berzin@univ-grenoble-alpes.fr} \url{http://www-ljk.imag.fr/membres/Corinne.Berzin/}}}


\AMSSUBJ{62G10; 53C65; 62F12; 60G60} 
\AMSSUBJSECONDARY{60G10; 60G15} 

\SUBMITTED{\today} 
\ACCEPTED{} 

\VOLUME{0}
\YEAR{2017}
\PAPERNUM{0}
\DOI{10.1214/YY-TN}

\ABSTRACT{
    Consider an affine field $X: \reels^2 \to \reels$, that is a process equal in law to $Z(A.t)$, where $Z$ is isotropic and $A: \reels^2 \to \reels^2$ is a linear self-adjoint transformation. The field $X$ and transformation $A$ will be supposed to be respectively Gaussian and definite positive.
    Denote  $0 < \lambda :=\frac{\lambda_2}{\lambda_1} \le 1$ the ratio of the eigenvalues of $A$, let $\lambda_1$, $\lambda_2$ with $\lambda_2 \le \lambda_1$.
    This paper is aimed at testing the null hypothesis ``$X$ \textit{is isotropic}'' versus the alternative ``$X$ \textit{is affine}''.   
    Roughly speaking, this amounts to testing ``$\lambda=1$'' versus ``$\lambda <1$''.
    By setting level $u$ in $\reels$, this is implemented by the partial observations of process $X$ through some particular level functionals viewed over a square $T_n$, which grows to $\reels^2$.
    This leads us to provide estimators for the affinity parameters that are shown to be almost surely consistent.
    Their asymptotic normality results provide confidence intervals for parameters.
    
    This paper offered an important opportunity to study general level functionals near the level $u$ and for a fixed bounded rectangle $T$ of $\reels^2$, part of the difficulties arises from the fact that the topology of level set $\textrm{C}_{T, X}(u)=\left\{t \in T: X(t)=u \right\}$ can be irregular, even if the trajectories of $X$ are regular.
    A significant part of the paper is dedicated to show the $L^2$-continuity in the level $u$ of these general functionals.}

	\usepackage{wasysym}
	\usepackage{textcomp}
	\usepackage{xspace}
	\usepackage{dsfont}
	\renewcommand{\le}{\ensuremath{\leqslant}}
	\renewcommand{\ge}{\ensuremath{\geqslant}}
	\DeclareMathOperator{\sym}{sym}
	\DeclareMathOperator{\card}{card}
	\newcommand{\rv}{random variable\xspace}
	
	\reversemarginpar
	
	\newcommand{\introo}[3][,]{{]}{#2}#1{#3}{[}}
	\newcommand{\introf}[2]{{]}{#1},{#2}{]}}
	
	\newcommand{\intrff}[2]{{[}{#1},\,{#2}{]}}
	\def\ds{\displaystyle}
	\def\bk{\boldsymbol{k}}
	\def\bm{\boldsymbol{m}}

	\def\bsX{\boldsymbol{X}}
	
	\def\supp{{\rm supp}}
	\def\sign{{\rm sign}}
	\def\tr{{\rm tr}}
	\def\ud{{\rm\,d}}
	\def\reels{\mathbb R}
	\def\entiers{\mathbb Z}
	\def\naturels{\mathbb N}
	\def\1{\mathds1}
	\def\Esp{{\rm E}}
	\newcommand{\esp}[2][]
     		{\ensuremath{\Esp_{{#1}}\kern-2pt\left[{#2}\right]}}
	\def\Var{{\rm Var}}
	\newcommand{\var}[2][]
     		{\ensuremath{\Var_{{#1}}\kern-2pt\left[{#2}\right]}}

	\newcommand{\abs}[1]{\left| {#1} \right|}
	\newcommand{\normp}[2][p]{\left\| {#2} \right\|_{#1}}
	
	\newenvironment{proofarg}[1]
   {\par\noindent
    \textit{#1.}\ }
   {\hfill{$\Box$}\par\noindent}

\begin{document}



%
%
\section{Introduction}
\par The aim of the present paper is to test the null hypothesis that a given real Gaussian process $X$ indexed in $\reels^2$ and living in the class of affine processes is isotropic. Recall that a process $(X(t), t \in \reels^2)$ is said to be an affine process when it is equal in law to $(Z(A.t), t \in \reels^2)$, where $Z$ is isotropic (see definition in Section \ref{sec:notations}, page \pageref{isotropic}) and $A: \reels^2 \to \reels^2$ is a linear self-adjoint matrix.
We assume that $X$ is partially observed through some level functionals of its level curve $\textrm{C}_{T_n, X}(u)$ for a fixed level $u$, say $\textrm{C}_{T_n, X}(u)=\left\{t \in T_n: X(t)=u \right\}$.
The set $T_n$ is a bounded square of $\reels^2$, having the following shape $\introo{-n}{n}^2$, with $n$ a positive integer. 
We are interested in the asymptotic as the square $T_n$ tends to $\reels^2$. The eigenvector directions $v_1$, $v_2$ forming an orthonormal basis of $\reels^2$
and the associated eigenvalues $\lambda_1$, $\lambda_2$ ($0<\lambda_2 \le \lambda_1$) specify the affinity $A$. Also $0 < \lambda \le 1$ is defined as the quotient of the eigenvalues, $\lambda :=  {\lambda_2}/{\lambda_1}$. Saying that the process $X$ is isotropic means that $\lambda=1$ (see Adler and Taylor \cite[Section 5.7]{MR2319516}). A fixed vector of $\reels^2$ being fixed, say $v^{\star}$,
it can always be written in the previous basis as:
\begin{equation*}
v^{\star}= \cos(\theta) v_1+ \sin(\theta) v_2.
\end{equation*}
It is always possible to choose $-\frac{\ds \pi}{\ds 2} < \theta \le \frac{\ds \pi}{\ds 2}$ (see Figure \ref{fig:figure1}, page \pageref{fig:figure1}). We use Wschebor estimators \cite{MR0871689} of the affinity parameters, defined as $\ds \lambda$ and $\theta$, based on the shape of  the level curve $\textrm{C}_{T_n, X}(u)$ corresponding to a given level $u$.
Also, a lot of papers have been devoted to the study of two-dimensional deformed fields, these later naturally model spatial and physical phenomena neither isotropic nor stationary. They are obtained by deforming a fixed isotropic field, say $Z: \reels^2 \to \reels$, thanks to a deterministic function, say $\Upsilon: \reels^2 \to \reels^2$, that transforms bijectively the index set. This model of deformed fields, say $X_{\Upsilon}:=Z \circ \Upsilon$, was introduced by Sampson and Guttorp in  \cite{Sampson:1992} but with only a stationary assumption on $Z$.
Among the papers that focus on the study of these fields, let us quote Allard et al. \cite{MR3476704}, Anderes et al. \cite{MR2543694}, Anderes and Stein \cite{MR2396813}, Fouedjio et al. \cite{MR3368820}, Fournier \cite{MR3877250} and Perrin and Senoussi \cite{MR1767607}.
For example, in \cite{MR1767607} the authors studied the models through the covariance function.
In this vein a large number of authors was interested in estimate deformation $\Upsilon$ thanks to the observation of the deformed field $X_{\Upsilon}$. As  precursors, Sampson and Guttorp in  \cite{Sampson:1992} used several observations on a sparse grid. Another approach was then to use only one observation of $X_{\Upsilon}$ on a dense grid, we can cite for example the relevant research of \cite{MR2543694}, \cite{MR2396813} and \cite{MR3368820}.
Notice that in \cite{MR3476704}, the authors focus on the case where the deformation $\Upsilon$ is linear which will also be our aim. The matrix representation of this deformation being the product of a diagonal and a rotation matrix, the authors calling the produced effect, ``geometric anisotropy".
However, as far as we know, very little research is based on the use of level sets of the observed process $X_{\Upsilon}$ and except for \cite{MR0871689} those quoted previously are not an exception to the rule.
In some cases, we may be interested in the use of the information given by some functional of level sets of the process rather than given by a simple realization of the process itself.
This is done for example in \cite{MR3877250} where the author chose as particular functional the Euler characteristic of some of excursion sets of the deformed field.
The covariance function of the underlying field $Z$ and the deformation $\Upsilon$ are unknown. The problem consists in identifying $\Upsilon$ using sparse data, that is through observations of some excursion sets of $X_{\Upsilon}$ above fixed levels by focusing on the mean Euler characteristic of the excursion sets for multiple windows of observation. For this doing some assumptions on $Z$ are required, among others a Gaussian property is asked.
The ambition of Fournier was to identify the deformation $\Upsilon$. However, it had proved impossible to distinguish between $\Upsilon$ and another deformation $\widetilde{\Upsilon}$ if the random fields $X_{\Upsilon}$ and $X_{\widetilde{\Upsilon}}$ have the same law. Indeed and because of the property of isotropy for process $Z$, this case occurs when $\widetilde{\Upsilon}= \rho \circ \Upsilon +x$, where $\rho$ is a rotation of $\reels^2$ and $x$ a vector of $\reels^2$. Thus,  the author almost entirely identifies $\Upsilon$ by proposing an identification of the matrix parameter $\Upsilon$ up to composition with an unknown rotation and with a translation.
Also, in a Note \cite{MR0978469}, Berzin and Iribarren proposed estimators  for the spectral matrix $\Sigma$ of the second order moments for a stationary and centred Gaussian bi-parametric random field. The estimators are defined as functionals over the level curves of the field and are inspired from those given by Caba\~na \cite{MR0898839} to estimate the affinity parameters previously defined. Making so, the used techniques could allow to estimate the ratio $\sqrt{\frac{\lambda_{-}}{\lambda_{+}}}$, where $0<\lambda_{-} \le \lambda_{+}$ are the eigenvalues of $\Sigma$. This last ratio being nothing else that $\lambda =\frac{\lambda_2}{\lambda_1}$ when $X=Z \circ A$ is an affine process. We see here a posteriori that inference based on one level set by the pioneering work of Caba\~na could already have been useful at this time for estimate the affinity parameter $\lambda$.

Concerning the tests for the hypotheses of isotropy of general random fields many methods have been provided.
An overview of methods available to test the isotropy have been very well summarized in Weller et al. \cite{MR3552737} but as previously explained, to my knowledge except the one given by Caba\~na \cite{MR0898839}, none used the observation of the level curves. The author's basic idea was to induce a distortion in the domain that breaks isotropy. This deformation becoming apparent in the length of the level sets. This last quoted work has deeply inspired the present paper. 
In \cite{MR0898839}, according to his own words, the author defines the class of affine processes as a reasonable alternative to test the null hypothesis that a given almost surely differentiable stationary real process with parameter in $\reels^2$ is isotropic.
More precisely Caba\~na, considers an affine process $\{X(t), t \in \reels^2\}$ that is equal in law to $\{Z(A . t), t \in \reels^2\}$, where $Z$ is isotropic (no necessarily Gaussian) and $A: \reels^2 \to \reels^2$ is a linear self-adjoint transformation. The angle $\theta_0$ defining the eigenvector directions $(\cos(\theta_0), \sin(\theta_0))$ and $(-\sin(\theta_0), \cos(\theta_0))$ and the repetitive eigenvalues $\lambda_1$, $\lambda_2$ specify the affinity $A$. As explained by the author there is no loss of generality in assuming $\lambda_1 \ge \lambda_2 >0$ and $\lambda_1 \lambda_2=1$.
The statistics used are based on the size and shape of the random surface $X$, corresponding to a fixed level $u$.
The regions where $X$ is observed are the indefinitely increasing family of rectangles $R_{\varrho}:= \{\varrho t: t \in R\}$, $\varrho \in \reels^+$ obtained from a fixed rectangle $R \subset \reels^2$. In order that the level sets $\mathcal{C}_{\textrm{R}_{\varrho}, \textrm{X}}(u):= \left\{s \in R_{\varrho}: X(s)=u \right\}$ are curves, some restrictions had to be imposed to the process $X$, (see Proposition \ref{Mario} for the discussion about the topic). So,  the paths of $X$ are required to be of class $\textbf{C}^2$ on $\reels^2$. Also
the pair $(X(0), \nabla X(0))$ is supposed to have a joint density ${p}_{X(0), \nabla X(0)}(x, x^{\prime})$ in $\reels \times \reels^2$, bounded for $x^{\prime}$ varying in a compact subset of $\reels^2$ and $x$ in some neighbourhood of $u$, $\nabla X$ standing for the gradient of $X$. The stationary random fields $\Theta$ and $\normp[2]{\nabla X}$ are defining the polar representation of $\nabla X$, that is for $t \in \reels^2$, $\nabla X(t)= \normp[2]{\nabla X(t)} (\cos(\Theta(t)), \sin(\Theta(t)))$, where $\normp[2]{\cdot}$ stands for the Euclidean norm in $\reels^2$.\\
In the first part of his paper, Caba\~na proposes estimators of the affinity parameters defined as $k:=(1 - \lambda^2)^{\frac{1}{2}}$ with $\ds \lambda:=\frac{\lambda_2}{\lambda_1}$ and $\theta_0$, based on the $u$ level curve of $X$ in the following way. First, he defines for  $f:\introf{-\pi}{\pi}\to \reels$, a bounded and measurable function, the general line integral
\begin{equation*}
{\cal{F}}(f, R):= \abs{R}^{-1} \int_{\mathcal{C}_{\textrm{R}, \textrm{X}}(u)} f(\Theta(t)) \ud\sigma_1(t),
\end{equation*}
where $\abs{R}$ denote the area of $R$ that is $\sigma_2(R)$, and for $d=1, 2$, $\sigma_d$ denotes the Lebesgue measure on $\reels^d$.
(Note that the assumptions made on $X$ ensure that these line integrals are well defined, see the beginning of the paragraph \textit{General level functionals}  in Section \ref{sec:notations} page \pageref{general level functionals} for a commentary about this topic).\\
Then he defines three statistics obtained as particular line integrals, let
\begin{align*}
{\cal{L}}(R)&:= \abs{R}^{-1} \int_{\mathcal{C}_{\textrm{R}, \textrm{X}}(u)} 1 \ud\sigma_1(t), \\
{\cal{C}}(R)&:= \abs{R}^{-1} \int_{\mathcal{C}_{\textrm{R}, \textrm{X}}(u)} \cos(2\Theta(t)) \ud\sigma_1(t),\\
{\cal{S}}(R)&:= \abs{R}^{-1} \int_{\mathcal{C}_{\textrm{R}, \textrm{X}}(u)} \sin(2\Theta(t)) \ud\sigma_1(t).
\end{align*}
The author states that assumptions made on process $X$ imply that these functionals have first and second moments, that is, Rice formulas of orders one and two are ensured. Thus, he proves that the affinity parameters previously defined as $k$ and $\theta_0$ can be obtained via the straightforward computations of $\mathbb{E}[{\cal{L}}(R)]$, $\mathbb{E}[{\cal{C}}(R)]$ and of $\mathbb{E}[{\cal{S}}(R)]$.
In this way, the author highlights the fact that
\begin{eqnarray*}
\tan(2\theta_0)= \mathbb{E}[{\cal{S}}(R)]/ \mathbb{E}[{\cal{C}}(R)] \mbox{\, and \, } g(k)=\sqrt{(\mathbb{E}[{\cal{C}}(R)])^2+ (\mathbb{E}[{\cal{S}}(R)])^2}/\mathbb{E}[{\cal{L}}(R)] \, , \end{eqnarray*}
where $g$ is a known continuously increasing function on $\intrff{0}{1}$, such that $g(0)=0$.\\
In view of studying the asymptotic properties of future estimators of $k$ and $\theta_0$ that will naturally result of the two above equalities, the author adds a condition of uniform mixing for process $X$ that we recall here. The field $X$ is said to be {\it uniformly mixing}, when there exists a decreasing function $\delta: \reels^+ \to \reels^{+}$, $\lim_{y \to +\infty} \delta(y)=0$, such that if $S$, $S^{\prime}$ are two Borel sets in $\reels^2$, $d:=\mbox{distance}(S, S^{\prime})$, $\sigma(S)$ is the $\sigma$-fields generated by $\{X(t): t \in S \}$, then for any events $A \in \sigma(S)$, $B \in \sigma(S^{\prime})$, the inequality $\abs{\mathbb{P}(A/B) - \mathbb{P}(A)}\le \delta(d)$ holds. Under the uniform mixing condition for $X$, Iribarren \cite{MR990398} has shown that the asymptotic distribution of $\varrho \abs{R}^{\frac{1}{2}}({\cal{F}}(f, \varrho R)- \mathbb{E}[{\cal{F}}(f,R)])$ (as $\varrho$ goes to infinity) is Gaussian with mean zero and variance $\sigma^{2}_{f}$.\\
This leads Caba\~na to conclude that the estimators $\widehat{\theta}_{0, \varrho}$ of $\theta_{0}$ and $\widehat{k}_{\varrho}$ of $k$ defined as $\widehat{\theta}_{0, \varrho}:= (1/2) \arg({\cal{S}}(\varrho R)/{\cal{C}}(\varrho R))$ and $\widehat{k}_{\varrho}:=g^{-1}(\sqrt{{\cal{C}}^2(\varrho R)+ {\cal{S}}^2(\varrho R)}/{\cal{L}}(\varrho R))$ are consistent in probability.\\
In the second part of his work, by pointing out that an isotropic field $X$ has $k=0$ (see \cite[Section 5.7]{MR2319516}) and the farther is $X$ from isotropy, the greater is $k$, then he  motivates two isotropy tests for the null hypotheses ``$X$ is isotropic'' versus the alternative ``$X$ is affine'', as follows.
For the first test that will be truthfully our main focus, he first establishes that under the null hypothesis, the random variables ${\cal{C}}(\varrho R)$, ${\cal{S}}(\varrho R)$ and ${\cal{L}}(\varrho R)$ are asymptotically independent and also that $\mathbb{E}[{\cal{C}}(R)]=\mathbb{E}[{\cal{S}}(R)]=0$. He deduces from this result that $\varrho^2 \abs{R}({\cal{C}}^2(\varrho R)+ {\cal{S}}^2(\varrho R))$ is asymptotically distributed as $\sigma_0^2 \chi^2_2$, where $\sigma_0^2$ denotes the common value of $\sigma_f^2$ when $f(\theta)$ is $\cos(2\theta)$ or $\sin(2\theta)$. This last result suggests using the test variable
$$
{\cal{A}}_{\varrho}(R):= \varrho^2 \abs{R} \frac{{\cal{C}}^2(\varrho R)+ {\cal{S}}^2(\varrho R)}{{\cal{L}}^2(\varrho R)}
$$
to define the rejection region as ${\cal{A}}_{\varrho}(R)> \gamma$ for the null hypothesis ``$k=0$''.
Under $k=0$, the asymptotic distribution of ${\cal{A}}_{\varrho}(R)$ is $\sigma_0^2 \, (\mathbb{E}[{\cal{L}}(R)])^{-2}\, \chi^2_2$, since ${\cal{L}}(\varrho R)$ converges in probability to $\mathbb{E}[{\cal{L}}(R)]$. Furthermore when $k \neq 0$, $\varrho^{-2} \abs{R}^{-1} {\cal{A}}_{\varrho}(R)$ converges in probability to $g^{2}(k) >0$, and then ${\cal{A}}_{\varrho}(R)$ is stochastically unbounded for $\varrho$ tending to infinity. Thus, the proposed critical region provides a consistent test for any positive constant $\gamma$.\\
Nevertheless, as previously announced, the author develops a second isotropy test, an $F$-test, to remedy the imperfections of the first proposed test.
This second test is still based on the level sets of $X$ and ensues from the following argumentation.
In the case where the null hypothesis completely specifies the distribution of the process, then $\sigma_0^2$ and $\mathbb{E}[{\cal{L}}(R)]$ can be evaluated since they only depend on the distribution of the subjacent process $Z$ and on $\lambda_1$ that worth $1$ since $k=0$ and $\lambda_1\lambda_2=1$.
Thus, in this case, all the ingredients are gathered to perform the asymptotic test.
But this is rarely actually the case, since before deciding if the isotropic model is suitable, an estimation of joint distribution of the process is necessary beforehand.
Thus Caba\~na developed a one-way analysis of variance test in the following manner.
The rectangle $R$ is still fixed and is cut into $n$ rectangles of equal area, forming a partition ${\cal{F}}:=\{R_1, R_2, \dots, R_n\}$.
Under the same uniform mixing condition as supposed in the first part of his paper, the previous cited work \cite{MR990398} also allows to establish that since $R_1, R_2, \dots, R_n$ are disjoint then $\varrho \abs{R_i}^{\frac{1}{2}}({\cal{C}}(\varrho R_i)- \mathbb{E}[{\cal{C}}(R_i)])$, $\varrho \abs{R_i}^{\frac{1}{2}}({\cal{S}}(\varrho R_i)- \mathbb{E}[{\cal{S}}(R_i)])$, $i=1, 2, \dots, n$ are asymptotically independent.
Thus, the analysis of the model
\begin{eqnarray*}
{\cal{C}}(\varrho R_i)=\mathbb{E}[{\cal{C}}(R_i)] + \mbox{error, \,} {\cal{S}}(\varrho R_i)=\mathbb{E}[{\cal{S}}(R_i)] + \mbox{error \,} (i=1, 2, \dots, n)
\end{eqnarray*} 
suggested the test variable
\begin{align*}
{\cal{B}}_{\varrho}({\cal{F}})&:= \frac{n^{-1}(\sum_{i=1}^{n} {\cal{C}}(\varrho R_{i}))^2 + n^{-1}(\sum_{i=1}^{n} {\cal{S}}(\varrho R_{i}))^2}
{\sum_{i=1}^{n} {\cal{C}}^2(\varrho R_i)+ \sum_{i=1}^{n}{\cal{S}}^2(\varrho R_i)-n^{-1}(\sum_{i=1}^{n} {\cal{C}}(\varrho R_{i}))^2 - n^{-1}(\sum_{i=1}^{n} {\cal{S}}(\varrho R_{i}))^2}\\
&=\frac{n({\cal{C}}^{2}(\varrho R)+{\cal{S}}^2(\varrho R))}
{\sum_{i=1}^{n} ({\cal{C}}^2(\varrho R_i)+ {\cal{S}}^2(\varrho R_i))-n({\cal{C}}^{2}(\varrho R)+{\cal{S}}^2(\varrho R))}.
\end{align*}
Under the hypothesis ``$k=0$'', the law of ${\cal{B}}_{\varrho}({\cal{F}})$ tends to that of a $F_{2, 2n-2}$ distribution as $\varrho$ goes to infinity. 
When $k \neq 0$, that is when $X$ is an affine process but not isotropic, the limit law is a noncentral $F$-distribution.
The critical region is then ${\cal{B}}_{\varrho}({\cal{F}})>F_{2, 2n-2}(\alpha)$, where $F_{2, 2n-2}(\alpha)$ is the $(1-\alpha)$-quantile of the $F_{2, 2n-2}$ distribution.\\
As ultimate remark Caba\~na specifies that the size of the rectangles $R_i$ has to be large enough to ensure a fair approximation to the limit distribution. Doing so necessarily limits the value of the number $n$ of rectangles of the partition.\\

Bearing in mind to estimate the affinity parameters $\lambda$ and $\theta$ defined at the beginning of this Introduction, by using the level sets of the process $X_{\Upsilon}$ and in the case where the deformation $\Upsilon:=A$ is a linear self-adjoint and definite positive one, our starting point was an example from \cite[Chap 3.6, Paragraph F, page 79]{MR0871689}.

\indent  In this example Wschebor, proposes estimators of the affinity parameters of a process $X$  indexed in $\reels^2$ based on the shape of its level curve, corresponding to a given level $u$.
The author explains as preamble that his work was largely based on Caba\~na \cite{MR0898839}, although his estimators were different.
More precisely, Wschebor considers $\{X(t), t \in \reels^2\}$ a $C^2$-affine process (not necessarily Gaussian) and fixes $u$ a level in $\reels$ supposed to be a regular value of $X$, (see Proposition \ref{Mario} for much details on sufficient hypothesis ensuring such a property).
Moreover, the process $X$ is supposed to be sufficiently regular in order to verify the one and second order Rice formula. 
Chapter 3 of the book gives many details about the topic and on particular hypotheses $H_{1,2}$ and $H_{2,2}$.
Furthermore the process $X$ is supposed satisfying a condition of $\eta$-dependence. We recall here this notion. The process $(X, \nabla X)$ is said to be $\eta$-dependent if $\normp[2]{t_i-t^{\prime}_j} \ge \eta >0$ ($i=1, \dots, k$, $j=1, \dots, \ell$, $k, \ell=1, 2, \dots$), then the random vectors $((X(t_1), \nabla X(t_1)), \dots, (X(t_k), \nabla X(t_k)))$, $((X(t^{\prime}_1), \nabla X(t^{\prime}_1)), \dots, (X(t^{\prime}_{\ell}), \nabla X(t^{\prime}_{\ell})))$ are independent. Notice that the condition of $\eta$-dependence for $X$ implies that of uniformly mixing required by Caba\~na in \cite{MR0898839}.  The regions where $X$ is observed are the indefinitely increasing family of rectangles $R_{\varrho}:= \{(\varrho_1 t_1, \varrho_2 t_2): (t_1, t_2) \in R\}$, if $\varrho:=(\varrho_1, \varrho_2)$, $\varrho_1, \varrho_2 \in \reels^+$ obtained from a fixed rectangle $R \subset \reels^2$.
By defining the following general functional ${\cal{J}}(f, R)$ of the fixed level $u$ by:
\begin{equation*}
{\cal{J}}(f, R):=  \abs{R}^{-1} \int_{\mathcal{C}_{\textrm{R}, \textrm{X}}(u)} f({\frac{\nabla X(t)}{\normp[2]{\nabla X(t)}}}) \ud\sigma_1(t),
\end{equation*}
for $f: S^{1} \to \reels$ or $f: S^{1} \to \reels^2$, a bounded function where $S^1$ is the boundary of the unit ball of $\reels^2$, the author considers two particulars functionals of the level set $u$.
The first one where $f\equiv {\bf 1}$ corresponds to the particular functional being the measure of dimensional area of the level set by unit of volume, say ${\cal{J}}({\bf 1}, R)$. The second one, say  ${\cal{J}}(\overrightarrow{f^{\star}}, R)$, corresponds to the function $\overrightarrow{f^{\star}}$ taking values in $\reels^2$, that will be evaluated in the value of the normalized gradient of $X$ and being its value if this one lives in the same half-plane as that of $v^{\star}$, a fixed unitary vector in $\reels^2$, and minus this value if not, that is $\overrightarrow{f^{\star}}(\omega):= \omega \times (\1_{\{\langle\omega, v^{\star}\rangle \ge 0\}}-\1_{\{\langle\omega, v^{\star}\rangle < 0\}})$ if $\omega$ lives in $S^1$. Note that with the notations of Caba\~na, this functional can be expressed as 
\begin{equation*}
{\cal{J}}(\overrightarrow{f^{\star}}, R)= \abs{R}^{-1} \int_{\mathcal{C}_{\textrm{R}, \textrm{X}}(u)} (\cos(\Theta(t)), \sin(\Theta(t)))\, \sign(\cos(\Theta(t)-\theta_1)) \ud\sigma_1(t),
\end{equation*}
if $v^{\star}=(\cos(\theta_1), \sin(\theta_1))$.\\
Under the hypotheses that $\inf(\varrho_1, \varrho_2) \to +\infty$, Wschebor showed that the ratio $ {\cal{Q}}(R_{\varrho}):= \frac{{\cal{J}}(\overrightarrow{f^{\star}}, R_{\varrho})}{{\cal{J}}({\bf 1}, R_{\varrho})}$ tends in probability toward $ \frac{ \mathbb{E}[{\cal{J}}(\overrightarrow{f^{\star}}, R)]}{\mathbb{E}[{\cal{J}}({\bf 1}, R)]}$. 
This convergence enabled him to propose probability consistent estimators of direction $\theta$, say $\widehat{\theta}_{\varrho}$ and also of $\lambda$, say $\widehat{\lambda}_{\varrho}$, based on the observation of the quotient ${\cal{Q}}(R_{\varrho})$.
More precisely by using the one order Rice formula, he proves that the ratio $\frac{ \mathbb{E}[{\cal{J}}(\overrightarrow{f^{\star}}, R)]}{\mathbb{E}[{\cal{J}}({\bf 1}, R)]}$ can be written in the fixed direct orthonormal basis $(v^{\star}, v^{\star\star})$ (see Figure \ref{fig:figure1} page \pageref{fig:figure1}) as $F_1(\lambda, \theta) v^{\star} +F_2(\lambda, \theta) v^{\star\star}$, while functions $F_i$, $i=1, 2$ are explicitly known (see equation (\ref{function F}) page \pageref{F1 F2}). Thus, his idea was to observe the ratio $\ds {\cal{Q}}(R_{\varrho})$ in the fixed basis $(v^{\star}, v^{\star\star})$, say $(X_\varrho, Y_\varrho)$. Following this way Wschebor has shown that the system of equations 
\begin{align}
\label{systemBis}
\begin{cases}
X_\varrho &= F_1(\lambda, \theta)\\
Y_\varrho & =   F_2(\lambda, \theta)
\end{cases}
\end{align} admits one and only one solution, say $(\widehat{\lambda}_{\varrho}, \widehat{\theta}_{\varrho}) \in \, \introo{0}{1}\times \introf{-\frac{\pi}{2}}{\frac{\pi}{2}}$ as long as $(X_\varrho, Y_\varrho)$ is such that $X_\varrho \neq \frac{2}{\pi}$ or $Y_\varrho \neq 0$ (see Figure \ref{fig:figure2} page \pageref{fig:figure2} and Figure \ref{fig:figure3} page \pageref{fig:figure3}). It is important to note that such an approach furnished estimators $(\widehat{\lambda}_{\varrho}, \widehat{\theta}_{\varrho})$ of the affinity parameters $(\lambda, \theta)$ in an implicitly manner which is not the case for those provided by Caba\~na, the latter not being linked together and can therefore be calculated independently of each other. Then the author claims that the derived estimators are consistent in probability. Note that this consistence property is valid under the hypothesis that $0 < \lambda <1$ and that $-\frac{\pi}{2} < \theta < \frac{\pi}{2}$, that is for $(\lambda, \theta)$ belonging to an open set. Nothing is suggested by the author in case where the parameters $(\lambda, \theta)$ belong to the boundary of the set; neither for the convergence rate of these estimators, nor concerning those of Caba\~na's paper. Moreover, Wschebor does not suggest any isotropic test since in fact his example was part of those intended to illustrate the Rice formula and was not there intended to be dealt with in details. 

Thus, our goal  was to revisit Wschebor example in \cite{MR0871689} by showing among other things the almost sure convergence of the estimators $(\widehat{\lambda}_{\varrho}, \widehat{\theta}_{\varrho})$ toward the true parameters $(\lambda, \theta)$, by exhibiting their rate of convergence via a central limit theorem, and also by paying particular attention to the case where the parameters live on the border of their definition domain. That naturally lead to propose an isotropic test by considering the ideas developed in Caba\~na \cite{MR0898839}.

\paragraph{Main contribution of the paper} In the present work following the way opened by Caba\~na in \cite{MR0898839}, we still work with a $C^2$-affine process, say $\{X(t), t \in \reels^2\}$, that is equal in law to $\{Z(A . t), t \in \reels^2\}$, where $Z$ is isotropic and $A: \reels^2 \to \reels^2$ is a linear self-adjoint transformation, additionally supposed to be a positive definite one.
We get rid of the $\eta$-dependence hypothesis for the process $X$ requested by Wschebor in \cite{MR0871689} and we add the hypothesis that the subjacent process $Z$ is Gaussian and some technical assumptions on its covariance.
We consider his proposed estimators of the affinity parameters $\lambda$ and $\theta$ by considering the particular following situation.
We select the fixed rectangle $R:=\introo{-1}{1} \times ]-1, 1[$, the sequence $\varrho=(\varrho_1, \varrho_2):=(n, n)$ with $n$ positive integer, and we consider the observation windows $T_n:=R_{\varrho}=\introo{-n}{n}^2$ with $n$ tending toward infinity.
For these observation windows and as previously explained in the second part of the Introduction we define the following general functional of the fixed level $u \in \reels$
\begin{equation*}
J_{f}^{(n)}(u):= \frac{1}{(2n)^2} \int_{\mathcal{C}_{\textrm{T}_n, \textrm{X}}(u)}  f({\frac{\nabla X(t)}{\normp[2]{\nabla X(t)}}}) \ud\sigma_1(t),
\end{equation*} 
for $f: S^{1} \to \reels^d$, $d=1, 2$ a bounded continuous function. We observe the ratio ${\cal{Q}}(T_n):= \frac{J_{\overrightarrow{f^{\star}}}^{(n)}(u)}{J_{\bf 1}^{(n)}(u)}$ in the fixed basis $(v^{\star}, v^{\star\star})$, say $(X_n, Y_n)$. Considering the system of equations (\ref{systemBis}), in the same way as proposed by Wschebor we construct in an implicit manner the estimators $(\widehat{\lambda}_n, \widehat{\theta}_{n})$ (see Proposition \ref{def estimateurs}).

Our main contribution consists, in the one hand, to establish the almost sure consistence of these estimators (Theorem \ref{convergence presque-sure}) in case where $0< \lambda <1$ and $-\frac{\pi}{2}< \theta < \frac{\pi}{2}$.
We need to prove Theorem \ref{convergence sure} stated below using an ergodic theorem (Adler \cite[\S 6.5]{MR0611857}) and Rice formula (see the seminal work of Rice \cite{MR0010932}).
\paragraph{Theorem \ref{convergence sure}}
For $f: S^{1} \to \reels$ a continuous and bounded function,
\begin{equation*}
J_{f}^{(n)}(u) \xrightarrow[n\to+\infty]{a.s.}\mathbb{E}[J_{f}^{(1)}(u)].
\end{equation*}
This theorem applied to the particular functions $\overrightarrow{f^{\star}}$ and 
 to the function ${\bf 1}$ taking values in $\reels$ and identically equal to one implies that,
\begin{equation*}
\frac{J_{\overrightarrow{f^{\star}}}^{(n)}(u)}{J_{\bf 1}^{(n)}(u)}=X_n v^{\star}+Y_n v^{\star\star} \xrightarrow[n\to+\infty]{a.s.} \frac{\mathbb{E}[J_{\overrightarrow{f^{\star}}}^{(1)}(u)]}{\mathbb{E}[J_{\bf 1}^{(1)}(u)]}= F_1(\lambda, \theta) v^{\star} +F_2(\lambda, \theta) v^{\star\star},
\end{equation*}
where the functions $F_1$, $F_2$ are defined by equation (\ref{function F}) page \pageref{F1 F2} (Proposition \ref{decomposition base}).
\noindent By noting $\overrightarrow{F}:= (F_1, F_2)$ and by using the fact that $\overrightarrow{F}$ is a $C^2$-diffeomorphism (Proposition \ref{diffeomorphisme}), we deduce from this theorem the almost sure convergence of $(X_n, Y_n)= \overrightarrow{F}(\widehat{\lambda}_n, \widehat{\theta}_{n})$ toward $\overrightarrow{F}(\lambda, \theta)$ and then of that of $(\widehat{\lambda}_n, \widehat{\theta}_{n})$ toward $(\lambda, \theta)$ provided that these last parameters live in an open set.\\
In the second hand our contribution was to propose a Central Limit Theorem (CLT) (Theorem \ref{convergence lambda et theta}) for those estimators and some confidence intervals (Corollary \ref{convergence lambda et theta Bis}), still in case where $0< \lambda <1$ and $-\frac{\pi}{2}< \theta < \frac{\pi}{2}$.
For proving the asymptotic normality of the estimators, we use a CLT for general functionals $J_f^{(n)}(u)$, that was the subject of Theorem \ref{Peccati} stated here. 
\paragraph{Theorem \ref{Peccati}}
For $f: S^{1} \to \reels$ a continuous and bounded function,
we have the following convergence,
\begin{align*}
\xi_{f}^{(n)}(u):= 2n \left({J_{f}^{(n)}(u) - \mathbb{E}[J_{f}^{(n)}(u)}]\right) \xrightarrow[n \to +\infty]{Law} \mathcal{N}(0; \mathit{\Sigma}_{f, f}(u)).
\end{align*}
We deduce from Theorems \ref{convergence sure} and \ref{Peccati} the following proposition.
\paragraph{Proposition \ref{convergence en loi}}
\begin{align*}
2n\left({
 \tfrac{\ds J_{\overrightarrow{f^{\star}}}^{(n)}(u)}{\ds J_{\bf 1}^{(n)}(u)} - \ds \tfrac{\mathbb{E}[\ds J_{\overrightarrow{f^{\star}}}^{(1)}(u)]}{\mathbb{E}[\ds J_{\bf 1}^{(1)}(u)]}
}\right) 
\xrightarrow[n \to +\infty]{Law} \mathcal{N}(0; \mathit{\Sigma^{\star}}(u)).
\end{align*}
And writing this last convergence result in the basis $(v^{\star}, v^{\star\star})$, gives that 
\begin{multline*}
2n \left({
F_1(\widehat{\lambda}_{n} , \widehat{\theta}_{n}) - F_1(\lambda, \theta)
}\right) v^{\star} 
+2n \left({
F_2(\widehat{\lambda}_{n} , \widehat{\theta}_{n}) - F_2(\lambda, \theta)
}\right) v^{\star \star}
\xrightarrow[n \to +\infty]{Law} \mathcal{N}(0; \mathit{\Sigma^{\star}}(u)),
\end{multline*}
so that  
\begin{align*}
2n \left({
\overrightarrow{F}(\widehat{\lambda}_{n} , \widehat{\theta}_{n}) - \overrightarrow{F}(\lambda, \theta)
}\right)
\xrightarrow[n \to +\infty]{Law} \mathcal{N}(0; \mathit{\Sigma_Q^{\star}}(u),
\end{align*}
where $\mathit{\Sigma_Q^{\star}}(u):=Q \times \mathit{\Sigma^{\star}}(u) \times Q^{t}$ and $Q$ is the change of basis matrix from the canonical basis $(\vec{i}, \vec{j})$ to the basis $(v^\star, v^{\star \star})$.

Once again by using the fact that $\overrightarrow{F}$ is a $C^2$-diffeomorphism, we deduce a CLT for estimators $(\widehat{\lambda}_{n} , \widehat{\theta}_{n})$, as soon as the affinity parameters $\lambda$ and $\theta$ do not belong to the edges of the parameter space.
\label{parameters space}\\
\textbf{Remark:} this method for estimating the parameters $(\lambda, \theta)$ is available only when the space parameters for $X$ is $\reels^2$. Indeed, when this space is $\reels^d$ ($d \ge 2$), we deal with a function $\overrightarrow{F}$ taking values from $\reels^{2d-2}$ into $\reels^d$. The arguments of $\overrightarrow{F}$ would be $(\mu_2, \dots, \mu_{d}, \theta_1, \dots, \theta_{d-1})$, $\mu_i$ being the ratio of the decreasing eigenvalues $\lambda_i$ of $A$, let $\mu_i:={\lambda_i}/{\lambda_1}$, $i=2, d$, and $\theta_j$ the angles involved in the parametrization in spherical coordinates of vector $v^{\star}$, $j=1, d-1$. The arrival space of function $\overrightarrow{F}$ is $\reels^d$ since the gradient $\nabla X$ lives into $\reels^d$. Thus, if we want the function $\overrightarrow{F}$ to be an isomorphism we must require the following condition on $d$, that is $2d-2=d$, which is only satisfied when $d=2$.

The tools for proving the CLT in Theorem \ref{Peccati} is the use of the CLT technique for functionals belonging to the Wiener-It\^o chaos.
This method has been developed by Nourdin et al. \cite{MR2573557}, Nualart et al. \cite{MR2118863} and Peccati and Tudor \cite{MR2126978} among others.
The idea for proving such a CLT is inspired by the precursor work of Kratz and Le{\'o}n \cite{MR1860517}.
In the case where the process $X$ is a stationary Gaussian isotropic  process indexed by $\reels^2$ and the observation window is a fixed bounded rectangle $T$, the authors propose a way to approximate $\sigma_1(\mathcal{C}_{\textrm{T}, \textrm{X}}(u))/\sigma_2(T)$, say $J_{\bf 1}(u):={\cal{J}}({\bf 1}, T)$, by other functionals $J_{\bf 1}(u, \sigma)$ ($\sigma \to 0$) with the help of a kernel $K_{\sigma}$ tending to the delta-Dirac function in $u$.
This is done in such a way that the approximating functional can be expressed 	as stochastic integrals with respect to Hermite polynomials.
It consists in using the coarea formula (see Federer \cite[Theorem 3.2.12 p 249]{MR0257325} and also Berzin et al. \cite[Corollary 2.1.1]{Berzin:2017}), for $J_{\bf 1}(u, \sigma)$, transforming then this functional initially expressed on the level curve as a temporal functional on the rectangle $T$ and getting then its Hermite expansion in $H(X)$ the space of real square integral functionals of the field $X$.
Using this technique in the case where the observation window is still a fixed bounded rectangle $T$,  as $\sigma \to 0$ we were able to express in turn a general functional on the level curve $\mathcal{C}_{\textrm{T}, \textrm{X}}(u)$ as stochastic integral with respect to Hermite polynomials, let $J_f(u):={\cal{J}}(f, T)$, $f$ being a general continuous and bounded function (Theorem \ref{equality L2 Bis}).
Applying this result to the squares $T_n$,  as $n\to\infty$, we obtain the asymptotic variance $\mathit{\Sigma}_{f, f}(u)$ of the centred and suitably rescaled general functionals $J_f^{(n)}(u)$, say  $\xi_f^{(n)}(u)$, with (Proposition \ref{variance asymptotique xi}).
Finally, applying the Peccati-Tudor method  \cite{MR2126978} and expressing the functionals $\xi_f^{(n)}(u)$ into the Wiener-It\^o chaos, we obtain the CLT in Theorem \ref{Peccati}.
The way of proceeding is completely based on the methods developed into the paper Estrade and Le{\'o}n \cite{MR3572325} itself inspired by the article \cite{MR1860517}.
In this work the authors show a CLT for the Euler characteristic of the excursions above $u$ of the field $X$ on $T$ as $T$ grows to $\reels^d$, $X$ being a stationary Gaussian isotropic process indexed in $\reels^d$.

Our real contribution for proving the CLT, apart from showing the non-degeneration of the asymptotic limit matrix variance $\mathit{\Sigma^{\star}}(u)$ (Remark \ref{non degenerate}), was to rely on the two functionals $J_f(u)$ and $J_{f}(u, \sigma)$ its approximation via the kernel $K_{\sigma}$, that is to show that $J_{f}(u, \sigma)$ is an $L^2$-convergent approximation of $J_{f}(u)$ (Proposition \ref{convergence sigma}).
It was the opportunity to obtain as a first bonus the $L^2$-continuity in the level $u$ of the random variable $J_{f}(u)$ (Theorem \ref{continuite longueur de courbe}), which is a very nice interesting result in itself.
We did not find it in the literature and we believe that this result deserves consideration.
The proof is far from obvious and implements a number of ideas developed in Berzin et al. \cite{Berzin:2017}, from which a local parametrization of the level set $\textrm{C}_{T, X}(u+\delta)$ near the level $u$ (see \cite[Theorem 3.1.2]{Berzin:2017}). The second bonus was brought by obtaining the one order Rice formula for $J_{f}(u)$, that is the exact computation of $\mathbb{E}[J_{f}(u)]$ (Proposition \ref{rice}).\\
Having a closer look at how the estimators were obtained, we finally show apart in the final section the almost sure convergence of estimators $\widehat{\lambda}_n$ (Theorem \ref{convergence presque-sure Bis}) and also their rate of convergence (Theorem \ref{raccord}) in case where the true parameters $(\lambda, \theta)$ belong to the boundary of the set $\introo{0}{1} \times\introo{-\frac{\pi}{2}}{\frac{\pi}{2}}$.

By supposing that the covariance function of the subjacent process $Z$ is known, this leads us also to propose as in  Section 2.1 of Caba\~na's paper, statistical tests  for the null hypothesis ``$X$ {\it is isotropic}'' versus the alternative ``$X$ {\it is affine}'' (Theorem \ref{xi}).
Those ones are suggested by the convergence result previously stated in Proposition \ref{convergence en loi}. As explained in the first part of the Introduction testing the isotropy means to test the null hypothesis $H_0: ~\lambda=1$ against $H_1:~\lambda <1$.
Because
$$
\left({\frac{\mathbb{E}[J_{\overrightarrow{f^{\star}}}^{(1)}(u)]}{\mathbb{E}[J_{\bf 1}^{(1)}(u)]}=F_1(\lambda, \theta) v^{\star} +F_2(\lambda, \theta) v^{\star\star} \neq \frac{2}{\pi} v^{\star}}\right)\iff \left({ \lambda <1 }\right),
$$
Proposition \ref{convergence en loi} ensures that under the hypothesis $H_0$,
\begin{align*}
T_{\overrightarrow{f^{\star}}}^{(n)}(u):=2n\left({
 \tfrac{\ds J_{\overrightarrow{f^{\star}}}^{(n)}(u)}{\ds J_{\bf 1}^{(n)}(u)} - \ds \frac{2}{\pi} v^{\star}
}\right) 
\xrightarrow[n \to +\infty]{Law} \mathcal{N}(0; \mathit{\Gamma}(u, \tau)),
\end{align*}
(see Corollary \ref{base}), where the asymptotic matrix variance $\mathit{\Gamma}(u, \tau)$ depends only on $u$ and $\tau$ the common value of the eigenvalues of matrix $A$. Since $\mathbb{E}[J_{\bf 1}^{(1)}(u)]={\bf C} \tau$ under hypothesis $H_0$, with ${\bf C}$ a computable constant, we estimate the parameter $\tau$ by an almost sure convergence estimator $\widehat{\tau}_n$ obtained by the application of Theorem \ref{convergence sure} in the case where function $f$ is identically equal to ${\bf 1}$. Thus, under the hypothesis of isotropy
\begin{align*}
\Xi_{\overrightarrow{f^{\star}}}^{(n)}(u):= (S_{\overrightarrow{f^{\star}}}^{(n)}(u))^{t} \,S_{\overrightarrow{f^{\star}}}^{(n)}(u) \xrightarrow[n \to +\infty]{Law} \chi_2^2,
\end{align*}
where $S_{\overrightarrow{f^{\star}}}^{(n)}(u)):= \mathit{\Gamma^{-\frac{1}{2}}}(u, \widehat{\tau}_n) \cdot T_{\overrightarrow{f^{\star}}}^{(n)}(u)$.\\
 As in  \cite{MR0898839}, the last convergence result naturally suggests using the test variable $\Xi_{\overrightarrow{f^{\star}}}^{(n)}(u)$.
 We built a consistent test with rejection region $\Xi_{\overrightarrow{f^{\star}}}^{(n)}(u) > \gamma$, since the previous equivalence ensures that  $\frac{1}{(2n)^2} \Xi_{\overrightarrow{f^{\star}}}^{(n)}(u)$ converges in probability to $b>0$ in case of anisotropy (see Remark \ref{rejection}).
 
We end with a remark pinpointing the fact that such an isotropy test cannot be implemented in this way in the case where the $X$ process parameters are not in $\reels^2$.
In fact, Proposition \ref{convergence en loi} could be generalized without difficulty to the case where the parameter space is $\reels^d$ ($d \ge 3$).
Nevertheless, it is not enough for considering the idea of adapting the isotropy test proof to the new situation.
If we refer to the previous remark made on page \pageref{parameters space}, using the same notations, under the isotropy hypothesis the ratio ${\mathbb{E}[J_{\overrightarrow{f^{\star}}}^{(1)}(u)]}/{\mathbb{E}[J_{\bf 1}^{(1)}(u)]}$ is not characterized by equalities $\mu_2= \dots=\mu_d=1$.
Roughly speaking the isotropy hypothesis implies that ${\mathbb{E}[J_{\overrightarrow{f^{\star}}}^{(1)}(u)]}/{\mathbb{E}[J_{\bf 1}^{(1)}(u)]}= c(d) v^{\star}$ ($c(d)$ being a constant depending only of $d$ and computed under hypothesis of isotropy), but the converse is not always true.
Such a test would not be then consistent.

\paragraph{Outline of the paper}

Section \ref{sec:notations} contains some notations, among others definitions, affine process and isotropic process.
It gives explicitly the type of general functionals on the level set $u$ we are looking for, say $J^{(n)}_f(u)$, where the observation window $T_n$ grows to $\reels^2$ when $n$ goes to infinity. It is also an opportunity in this part to introduce some concepts, when the observation window is a fixed rectangle $T$ of $\reels^2$. They deal with regularity properties of level sets and Rice formulas, that are closed formulas for the first and second moments for functionals defined on the level set, say $J_f(u)$. Section \ref{sec:hypotheses} gives the assumptions on the interest model and some examples of processes satisfying these assumptions.

In Section \ref{fonctionnelles vues dans le chaos} the observation window is fixed by taking a rectangle $T$. We begin by showing in Section \ref{continuity in the level} the $L^2$-continuity in the level $u$ for $J_f(u)$, $f$ being a general continuous and bounded function.
In Section \ref{chaos section}, first by using the coarea formula for an approximation of the functional, say $J_f(u, \sigma) (\sigma \to 0)$, we express this last one as stochastic integral with respect to Hermite polynomials. Then by using the $L^2$-continuity in the level for $J_f(u)$ we establish the $L^2$-convergence of $J_f(u, \sigma)$ towards $J_f(u)$, from which we deduce an Hermite expansion for the initial functional $J_f(u)$.
As a derived product of the $L^2$-continuity in the level previously obtained we deduce in Section \ref{Rice} a one order Rice formula for such functionals, in other words we compute $\mathbb{E}[\smash{J_f(u)}]$.

In section \ref{convergence de fonctionnelles}, the observation windows $(T_n)_n$ is a sequence of open bounded squares of $\reels^2$, with the following form $T_n:= \introo{-n}{n}^2$ with $n \in \naturels^{\star}:=\{x \in \entiers, x >0\}$, and $n$ tends to infinity.
We focus on convergence results for $J^{(n)}_f(u)$.
Section \ref{almost sure convergence for the functionals} is devoted to establish the almost sure convergence of $J^{(n)}_f(u)$  for general function $f$. In Section \ref{convergence en loi pour les fonctionnelles}, by using the Hermite expansion of $J^{(n)}_f(u)$ derived from section \ref{chaos section}, we give the rate of this convergence. First the asymptotic variance as the squares $T_n$ grow to $\reels^2$ is expressed as a series and we give an explicit lower bound.
We then proved the asymptotic normality for the centred and suitably rescaled general functionals $J^{(n)}_f(u)$ through the Peccati-Tudor Theorem.
The results obtained in Section \ref{convergence de fonctionnelles} enable us, by considering in Section \ref{Les estimateurs} the particular function $\overrightarrow{f^{\star}}$, to give rise to a first result, the definition of the estimators $\widehat{\lambda}_{n}$ and $\widehat{\theta}_{n}$ of the affinity parameters $\lambda$ and $\theta$. Secondly, we prove their almost sure consistency and also their rate of convergence in law in case where the true parameters live in an open set. Also, the coefficients of the asymptotic matrix variance are computed in Appendix \ref{ann:proof} and a lower bound is given for its determinant.
This law convergence result gives rise to confidence intervals for parameters $\lambda$ and $\theta$ in the specific special case where the covariance $r_z$ of the underlying isotropic process is known.

In Section \ref{Complementary results} some complementary convergence results for $\widehat{\lambda}_{n}$ are proposed when the affinity parameters $\lambda$ and $\theta$ belong to the edges of the parameter space, including the particular case where $\lambda=1$.
Finally, supposing that the covariance $r_z$ is known, we conclude by proposing an isotropy test.

This paper is complemented with Appendix \ref{ann:proof} giving technical proofs of some lemmas.
%
%
\section{Notations and hypotheses}
\subsection{Notations}
\label{sec:notations}

Let us give some definitions, notations and some propositions and theorems related to properties of level sets and Rice formulas.

In the following $T$ is an open bounded rectangle of $\reels^2$.\\
Let $(\Omega, \cal{A}, P)$ be a probability space and $X:   \Omega \times T \subset \Omega \times \reels^2 \to \reels$ a process continuously differentiable on $T$, that is $X \in \textbf{C}^{1}(T)$.
\paragraph{Level sets}
\label{LevelSetsPar}
We denote $\textrm{D}^{\textrm{r}}_{\textrm{X}}$ the following set 
 \begin{align*}
\textrm{D}^{\textrm{r}}_{\textrm{X}}:=\left\{t \in
T : \nabla X(t)\,\,\mbox{is of  rank } 1
\right\}=
\left\{t \in
T : \normp[2]{\nabla X(t)} \neq 0
\right\},
\end{align*}
where $\normp[2]{\cdot}$ denotes the Euclidean norm in $\reels^{2}$, while $\nabla X$ stands for the Jacobian of $X$.\\
For $u \in \reels$ we define the level set at $u$ as:
$$
\mathcal{C}_{\textrm{T}, \textrm{X}}(u):=\left\{t \in T: X(t)=u \right\}.
$$ 
This set can be very irregular, but its intersection with $ \textrm{D}^{\textrm{r}}_{\textrm{X}}$, let
$\mathcal{C}_{\textrm{T}, \textrm{X}}^{\textrm{D}^{\textrm{r}}}(u):=\mathcal{C}_{\textrm{T}, \textrm{X}}(u)\cap \textrm{D}^{\textrm{r}}_{\textrm{X}}$ is a $\textbf{C}^1$-manifold of dimension one. Indeed, at each point $t \in\mathcal{C}_{\textrm{T}, \textrm{X}}^{\textrm{D}^{\textrm{r}}}(u)$ the Jacobian matrix of $X$ is of full rank one. Thus, an application of the implicit theorem in a neighbourhood of each point $t$ provides a local parametrization of the level set  $\mathcal{C}_{\textrm{T}, \textrm{X}}^{\textrm{D}^{r}}(u)$. For more details concerning the proof of the construction of an atlas of $\mathcal{C}_{\textrm{T}, \textrm{X}}^{\textrm{D}^{r}}(u)$, the reader can look at the beginning of the proof of \cite[Theorem 3.1.2]{Berzin:2017}. 

The following proposition (see Aza{\"\i}s and Wschebor \cite[Proposition 6.12]{MR2478201}) gives a sufficient condition to ensure that the level is a regular value of $X$.
\begin{proposition}
\label{Mario}
Let $u \in \reels^2$. We assume the following:
\begin{itemize}
\item The paths of $X$ are of class $\textbf{C}^2$ on $T$.
\item For each $t \in T$, the pair $(X(t), \nabla X(t))$ has a joint density ${p}_{X(t), \nabla X(t)}(x, x^{\prime})$ in $\reels \times \reels^2$, which is bounded for $(t, x^{\prime})$ varying in a compact subset of  \,$T \times \reels^2$ and $x$ in some neighbourhood of $u$.
\end{itemize}
Then 
\begin{align*}
\mathbb{P}(\omega \in \Omega, \exists\, t \in T, X(t)(\omega) = u, \normp[2]{\nabla X(t)(\omega)}=0)=0,
\end{align*}
that is almost surely
\begin{align}
\label{vide}
 \, \,\, \, \mathcal{C}_{\textrm{T}, \textrm{X}}^{\textrm{D}^{\textrm{r}}}(u)=\mathcal{C}_{\textrm{T}, \textrm{X}}(u).
\end{align}
\end{proposition}
\begin{remark}
\label{videbis}
Note that if $X$ is a Gaussian process twice continuously differentiable on $T$, that is $X \in \textbf{C}^2(T)$, such that for each $t \in T$, the vector $(X(t), \nabla X(t))$ has a joint density, then the last proposition hypotheses are verified and 
equality (\ref{vide}) holds.
In that case the reader can also be referred to \cite[Proposition 3.3.2]{Berzin:2017}.
\end{remark}
\paragraph{General level functionals}
$S^{1}$ is the boundary of the unit ball of $\reels^2$ and
for $d=1, 2$,  $\sigma_d$ denotes the Lebesgue measure on
$\reels^d$.
\label{general level functionals}
For $f: S^{1} \to \reels$ a continuous and bounded function, we define the following general functional $J_{f}(u)$ of the fixed level $u$ by:
\begin{equation*}
J_{f}(u):= \frac{1}{\sigma_{2}(T)} \int_{\mathcal{C}_{\textrm{T}, \textrm{X}}^{\textrm{D}^{\textrm{r}}}(u)} f(\nu_{X}(t)) \ud\sigma_1(t),
\end{equation*}
where $\nu_{X}(t):= {\frac{\nabla X(t)}{\normp[2]{\nabla X(t)}}}$.\\
\begin{remark}
\label{existence fonctionnelles}
Note that for general process $X$, the finiteness of this integral is not necessarily guaranteed. However, this integral makes sense, for example, under assumptions of Proposition \ref{Mario}. More precisely, we substitute hypotheses concerning the process $X$ on $T$ in the last proposition by those on $O$, where $O$ is any open set such that $T \subset \overline{T} \subset O$, where $\overline{T}$ stands for the closure of $T$. In that case we deduce that almost surely $\mathcal{C}_{{\scriptsize \overline{\textrm{T}}, \textrm{X}}}^{{\scriptsize \textrm{D}^{\textrm{r}}}}(u)=\mathcal{C}_{{\scriptsize \overline{\textrm{T}}, \textrm{X}}}(u)$. Then we construct a partition of unity for the compact manifold $\mathcal{C}_{{\scriptsize \overline{\textrm{T}}, \textrm{X}}}^{{\scriptsize \textrm{D}^{\textrm{r}}}}(u)$ of $\reels^2$, almost surely included in 
$
\textrm{D}^{{\scriptsize \textrm{r}}}_{{\scriptsize \textrm{O}, \textrm{X}}}:=
\left\{t \in
O : \normp[2]{\nabla X(t)} \neq 0
\right\}$.
In that way as mentioned above, by using a local parameterization of the level set $\mathcal{C}_{{\scriptsize \textrm{T}, \textrm{X}}}^{{\scriptsize \textrm{D}^{\textrm{r}}}}(u)$, we prove that almost surely,
$\sigma_1(\mathcal{C}_{{\scriptsize \textrm{T}, \textrm{X}}}^{{\scriptsize \textrm{D}^{\textrm{r}}}}(u))< +\infty$ (see also Remark \ref{finitude courbe de Niveau}). Since function $f$ is bounded, this construction guarantees the existence of functional $J_{f}(u)$. We also refer the reader to \cite[Theorem 3.1.2]{Berzin:2017} for such a partition of unity construction.
\end{remark}

\noindent Note that in this article, $T$ will increase to $\reels^2$.
In this regard, we will consider $T_n:=\introo{-n}{n}^2$, open squares of $\reels^2$, with $n \in \naturels^{\star}$, a positive integer, and let $n$ tends to infinity.
In this case, and with the comment made in Remark \ref{existence fonctionnelles}, the corresponding functional on the set $u$ will be denoted by $J_f^{(n)}(u)$, that is
\begin{equation*}
J_{f}^{(n)}(u):= \frac{1}{(2n)^2} \int_{\mathcal{C}_{\textrm{T}_n, \textrm{X}}^{\textrm{D}^{\textrm{r}}}(u)} f(\nu_{X}(t)) \ud\sigma_1(t).
\end{equation*}
The process $X:  \Omega \times \reels^2 \to \reels$ will at least be continuously differentiable on $\reels^2$, that is $X \in \textbf{C}^{1}(\reels^2)$ and $\textrm{D}^{\textrm{r}}_{\textrm{X}}$ will be the following set 
 \begin{align*}
\textrm{D}^{\textrm{r}}_{\textrm{X}}:=\left\{t \in
\reels^2 : \nabla X(t)\,\,\mbox{is of  rank } 1
\right\}=
\left\{t \in
\reels^2 : \normp[2]{\nabla X(t)} \neq 0
\right\}.
\end{align*}

\paragraph{Rice formulas} The following results  give Rice formulas for the first and second moments for general functionals defined on the level set of the process $X$ (see \cite{MR2478201} or \cite{Berzin:2017} for references). \\
By applying the coarea formula, for which a statement is expressed later in this text in Corollary \ref{coarea formula}, to general functionals defined on the level sets $\mathcal{C}_{\textrm{T}, \textrm{X}}^{\textrm{D}^{\textrm{r}}}(u)$ and taking expectation afterwards, one obtains the well-known Kac-Rice formula. This formula gives a computation of the expectation of these general functionals for almost any level $u$. A proof of this formula is given for example in \cite[Proposition 2.2.1]{Berzin:2017}. Let us recall it here by adapting its statement to our context. \\

Let $\widetilde{Y}: \Omega \times T \subset  \Omega \times \reels^2 \to \reels$ a process defined on $T$ and continuous on $ \textrm{D}^{\textrm{r}}_{\textrm{X}}$. Let us consider the following assumptions:
\begin{itemize}
\item ${\bf H_1}$: The function 
\begin{equation*}
\textbf{u} \longmapsto \mathbb{E}\left[\int_{\mathcal{C}_{\textrm{T}, \textrm{X}}^{\textrm{D}^{\textrm{r}}}(\textbf{u})}|{\widetilde{Y}(t)}| \ud \sigma_{1}(t)\right],
\end{equation*}
is a continuous function of the variable  $\textbf{u}$.\\
\item ${\bf H_2}$: The function 
\begin{equation*}
\textbf{u} \longmapsto \int_{T}\textrm{p}_{X(t)}(\textbf{u})\mathbb{E}[|{\widetilde{Y}(t)}|  \normp[2]{\nabla
X(t)}|X(t)=\textbf{u}]\, \ud t,
\end{equation*}
is a continuous function of the variable $\textbf{u}$. \\
\item  ${\bf H_3}$: For almost any
$t \in T$, the density of 
$X(t)\,,\textrm{p}_{X(t)}(\cdot)$ exists.\\
\end{itemize}
Proposition 2.2.1 implies the following proposition.
\begin{proposition}[\textit{Kac-Rice formula}]
If $X$ and $\widetilde{Y}$ satisfy the assumption (${\bf H_1}$ or ${\bf H_2}$) and ${\bf H_3}$, then for almost any $u \in \reels$ one has
\begin{equation}
\label{Kac-Rice formula}
\mathbb{E}\left[\int_{\mathcal{C}_{\textrm{T}, \textrm{X}}^{\textrm{D}^{\textrm{r}}}(u)}{\widetilde{Y}(t)}\ud \sigma_{1}(t)\right]
=
\int_{T}\textrm{p}_{X(t)}(u)\mathbb{E}[{\widetilde{Y}(t)}  \normp[2]{\nabla
X(t)}|X(t)=u]\, \ud t.
\end{equation}
\end{proposition}
However, in applications the interest is focused on a fixed level $u$. So, we add the two following assumptions to ensure the validity of the Kac-Rice formula for all level.
\begin{itemize}
\item ${\bf H_4}$: The function 
\begin{equation*}
\textbf{u} \longmapsto \mathbb{E}\left[\int_{\mathcal{C}_{\textrm{T}, \textrm{X}}^{\textrm{D}^{\textrm{r}}}(\textbf{u})}{\widetilde{Y}(t)}\ud \sigma_{1}(t)\right],
\end{equation*}
is a continuous function of the variable  $\textbf{u}$.\\
\item ${\bf H_5}$: The function 
\begin{equation*}
\textbf{u} \longmapsto \int_{T}\textrm{p}_{X(t)}(\textbf{u})\mathbb{E}[{\widetilde{Y}(t)} \normp[2]{\nabla
X(t)}|X(t)=\textbf{u}]\, \ud t,
\end{equation*}
is a continuous function of the variable  $\textbf{u}$.
\end{itemize}
Adding the two previous assumptions the two sides of  equality (\ref{Kac-Rice formula}) become now continuous as function of the level. Thus, we enunciate the following formula of Rice, for which a proof is given in \cite[Theorem 3.1.1]{Berzin:2017}.
\begin{theorem}[\textit{One order Rice formula}]
\label{one order Rice formula}
If $X$ and $\widetilde{Y}$ satisfy the assumption (${\bf H_1}$ or ${\bf H_2}$), ${\bf H_3}$, ${\bf H_4}$ and ${\bf H_5}$ then for any $u \in \reels$ one has
\begin{equation*}
\mathbb{E}\left[\int_{\mathcal{C}_{\textrm{T}, \textrm{X}}^{\textrm{D}^{\textrm{r}}}(u)}{\widetilde{Y}(t)}\ud \sigma_{1}(t)\right]
=
\int_{T}\textrm{p}_{X(t)}(u)\mathbb{E}[{\widetilde{Y}(t)}  \normp[2]{\nabla
X(t)}|X(t)=u]\, \ud t.
\end{equation*}
\end{theorem}
\begin{remark}
Note that Rice formula only requires for $\widetilde{Y}$ to be continuous on $ \textrm{D}^{\textrm{r}}_{\textrm{X}}$, that will be the case for the process $\widetilde{Y}:= f(\nu_X) \1_{\textrm{D}^{\textrm{r}}_{\textrm{X}}}$ defined in the paragraph \textit{General level functionals}, page  \pageref{general level functionals}.
\end{remark}
Now we give conditions under which a two-order Rice formula will be valid. The reader is referred to  \cite[Theorem 3.3.3]{Berzin:2017} and adapted to the future situation.

Let us consider the following assumptions:
\begin{itemize}
\item ${\bf H_a}$: $X: \Omega \times T \subset  \Omega \times \reels^2 \to \reels$ is a centred Gaussian random field twice continuously differentiable on $T$, that is $X \in \textbf{C}^{2}(T)$.
\item ${\bf H_b}$: For all $t \in T$, the vector $(X(t), \nabla X(t))$ has a density.
\item ${\bf H_c}$ For all $t_1, t_2 \in T$, such that $t_1 \neq t_2$, the vector $(X(t_1), X(t_2))$ has a density, say ${p}_{X(t_1), X(t_2)}(\cdot, \cdot)$. 
\item ${\bf H_d}$:
\begin{align}
\label{sup esperance}
\mathbb{E}[\sup_{t \in T} \normp[1, 2]{\nabla^2X(t)}^{(s)}]^{4} < +\infty,
\end{align}
where $\nabla^2X$ stands for the $2\times2$ Hessian matrix of $X$ and $\normp[1,2]{\cdot}^{(s)}$ stands for the norm in
$\mathfrak{L^2_{s}}(\reels^{2}, \reels)$, the vectorial space of symmetric  linear  continuous functions from $\reels^{2}$ to $\reels$.
\item ${\bf H_e}$: $Y: \Omega \times T \subset  \Omega \times \reels^2 \to \reels$ is a continuous process defined on $T$ as $Y(t):=G(t, \nabla X(t))$, where 
\begin{align*} 
G: T \times \reels^{2}& \longrightarrow \reels\\
(t, z) & \mapsto G(t, z),
\end{align*}
is a bounded continuous function on $T \times \reels^{2}$.
\end{itemize}
Let us now state the following assumption  ${\bf H_6}$.
\begin{itemize}
\item ${\bf H_6}$: For all $u \in \reels$,
\begin{multline*}
\int_{T \times T}
\mathbb{E} \left[\abs{Y(t_1)} \abs{Y(t_2)} \normp[2]{\nabla
X(t_1)} \normp[2]{\nabla
X(t_2)}|X(t_1)=X(t_2)=u\right]\\
 \times {p}_{X(t_1), X(t_2)}(u, u)\, \ud t_1\, \ud t_2 < +\infty.
\end{multline*}
\end{itemize}
Let us express the following two-order Rice formula.
\begin{theorem}[\textit{Two-order Rice formula}]
\label{two-order Rice formula}
If $X$ and $Y$ satisfy the assumptions ${\bf H_a}$, ${\bf H_b}$, ${\bf H_c}$, ${\bf H_d}$, ${\bf H_e}$ and ${\bf H_6}$, then for any $u \in \reels$ one has
\begin{align}
\label{formule de Rice d'ordre deux}
\lefteqn{\mathbb{E}[\int_{\mathcal{C}_{\textrm{T}, \textrm{X}}^{\textrm{D}^{\textrm{r}}}(u)} Y(t) \ud \sigma_1(t)]^2}\\ \nonumber
&= \int_{T \times T}
	\mathbb{E} \left[Y(t_1)\, Y(t_2) \normp[2]{\nabla
X(t_1)} \normp[2]{\nabla
X(t_2)}|X(t_1)=X(t_2)=u\right] \,\times \\ \nonumber
 &\qquad\qquad\qquad\hfill{p}_{X(t_1), X(t_2)}(u, u)\ud t_1\ud t_2.
\end{align}
\end{theorem}
From this theorem we deduce the following corollary.
\begin{corollary}
\label{gaussien stationnaire}
Let $X:  \Omega \times \reels^2 \to \reels$ be a centred stationary Gaussian random field twice continuously differentiable on $\reels^2$, that is $X \in \textbf{C}^{2}(\reels^2)$, with covariance $r_x$ (belonging necessarily to $\textbf{C}^4(\reels^2)$) satisfying the following hypotheses:
\begin{itemize}
\item $\nabla ^2 r_x(0)$ is a non-degenerate matrix, where $\nabla^2r_x$ stands for the $2\times2$ Hessian matrix of $r_x$ .
\item $r_x^2(0)-r_x^2(t) \neq 0$, for all $t \neq 0$.
\end{itemize}
Then assumptions ${\bf H_a}$, ${\bf H_b}$, ${\bf H_c}$ and ${\bf H_d}$ are checked so that for all $u \in \reels$, almost surely
\begin{align*}
\, \,\, \, \mathcal{C}_{\textrm{T}, \textrm{X}}^{\textrm{D}^{\textrm{r}}}(u)=\mathcal{C}_{\textrm{T}, \textrm{X}}(u).
\end{align*}
Furthermore if $X$ and $Y$ satisfy the assumption ${\bf H_e}$ then they satisfy also assumptions ${\bf H_6}$. Thus, for any $u \in \reels$ the two-order Rice formula (\ref{formule de Rice d'ordre deux}) is valid.
\end{corollary}
The proof of this corollary provides from the following argumentation.\\
Hypotheses made on $X$, ensure that $X$ is a Gaussian process twice continuously differentiable on $T$, such that for each $t \in T$ the vector $(X(t), \nabla X(t))$ has a joint density, since $X(t)$ and $\nabla X(t)$ are independent random variables. Thus 
Remark \ref{videbis} gives equality (\ref{vide}). Furthermore, these hypotheses ensure that assumption ${\bf H_c}$ is obviously checked. By using the Borel-TIS inequality given in Adler and Taylor \cite[Theorem 2.1.1]{MR2319516}, one can show that assumption ${\bf H_d}$ is also checked. \\
Now let $X$ and $Y$ satisfying the assumption ${\bf H_e}$. The fact that $G$ is a bounded function on $T \times \reels^{2}$, implies that $Y$ is a bounded process on $T$. Thus, to check assumption  ${\bf H_6}$, it is enough to establish that for all $u \in \reels$, one has \begin{align*}
\int_{T \times T}
\mathbb{E} \left[\normp[2]{\nabla
X(t_1)} \normp[2]{\nabla
X(t_2)}|X(t_1)=X(t_2)=u\right]
 \times {p}_{X(t_1), X(t_2)}(u, u)\, \ud t_1\, \ud t_2 < +\infty.
 \end{align*}
 The proof of this finiteness can be founded for example in Berzin and Wschebor \cite[Proposition 2]{MR1222362}. \\
 Thus $X$ and $Y$ satisfy assumptions given in Theorem \ref{two-order Rice formula} and the two-order Rice formula ensues from it.
 \begin{remark}
 \label{finitude courbe de Niveau}
 One can observe that under the hypotheses of Corollary \ref{gaussien stationnaire} for $X$, by taking $Y \equiv 1$, one has proved that  for any $u \in \reels$,
\begin{multline*}
 \mathbb{E}[ \sigma_1(\mathcal{C}_{\textrm{T}, \textrm{X}}^{\textrm{D}^{\textrm{r}}}(u))]^2 =\\
  \int_{T \times T}
	\mathbb{E} \left[\normp[2]{\nabla
X(t_1)} \normp[2]{\nabla
X(t_2)}|X(t_1)=X(t_2)=u\right] \,\times {p}_{X(t_1), X(t_2)}(u, u)\ud t_1\ud t_2 < +\infty.
\end{multline*}
Note that in that case we proved that almost surely $\mathcal{C}_{{\scriptsize \textrm{T}, \textrm{X}}}^{\scriptsize\textrm{D}^{\textrm{r}}}(u) < +\infty$ (resp. for all $n \in \naturels^{\star}$, $\mathcal{C}_{{\scriptsize \textrm{T}_n, \textrm{X}}}^{\scriptsize\textrm{D}^{\textrm{r}}}(u) < +\infty$). Thus, for $f: S^{1} \to \reels$ a continuous and bounded function, the functional $J_{f}(u)$ (resp. $J_{f}^{(n)}(u)$) is almost surely well defined.

Also note that under the hypotheses of Corollary \ref{gaussien stationnaire} for $X$, those slightly modified of Proposition \ref{Mario} are also checked (see Remark \ref{existence fonctionnelles}). If we proceed as indicated in the previous remark, we then obtain another way to establish the almost sure finiteness of $\mathcal{C}_{{\scriptsize \textrm{T}, \textrm{X}}}^{\scriptsize\textrm{D}^{\textrm{r}}}(u)$ and also that of $\mathcal{C}_{{\scriptsize \textrm{T}_n, \textrm{X}}}^{\scriptsize\textrm{D}^{\textrm{r}}}(u)$. This explicit way provides another manner for proving the almost sure existence of functional $J_{f}(u)$ (resp. $J_{f}^{(n)}(u)$).
 \end{remark}
 \begin{remark}
 \label{noncontinuite}
Note that process $\widetilde{Y}:= f(\nu_X) \1_{\textrm{D}^{\textrm{r}}_{\textrm{X}}}$ defined in the paragraph \textit{General level functionals} is such that processes $X$ and $\widetilde{Y}$ do not verify hypothesis  ${\bf H_e}$, since $\widetilde{Y}$ is continuous on $ \textrm{D}^{\textrm{r}}_{\textrm{X}}$ and a priori not continuous on $T$. So, to compute $\mathbb{E}[J_{f}(u)]^2$ we will have to play close to the vest and be creative by generalizing the two order Rice formula to our purpose (see forthcoming Lemma \ref{continuite Y}).
\end{remark}

\paragraph{Isotropic process} 
A process $(Z(t), t \in \reels^2)$ is said to be {\it isotropic} if it is a stationary process and if for any isometry $J$ in $\reels^2$, $k \in \naturels:=\{x \in \entiers, x \ge0\}$ and $t_1, \dots, t_k \in \reels^2$, the joint laws of $(Z(t_1), \dots, Z(t_k))$ and $(Z(J.t_1), \dots, Z(J.t_k))$ are the same.

\label{isotropic}

\paragraph{Affine process}
A process $(X(t), t \in \reels^2)$ is said to be an {\it affine} process when it is equal in law to $(Z(A.t), t \in \reels^2)$, where $Z$ is isotropic and $A: \reels^2 \to \reels^2$ is a linear self-adjoint transformation. \\
If the eigenvalues of $A$ are denoted by $\lambda_1, \lambda_2$, saying that the process $X$ is isotropic means that $\lambda_1=\lambda_2$ (see Adler and Taylor \cite[Section 5.7]{MR2319516}).

\paragraph{Hermite polynomials}
We use the Hermite polynomials $(H_q)_{q \in \naturels}$ defined by 
$$
H_q(x):= (-1)^q \phi^{-1}(x)  \ds \frac{\ud ^q}{\ud  x^q}(\phi(x) ),
$$
where $\phi$ denotes the standard Gaussian density on $\reels$.\\
They provide an orthogonal basis of $L^2(\reels, \phi(x) \ud x)$.
We also denote by $\phi_m$ the standard Gaussian density on $\reels^m$, for $m=2$ or $3$.

For ${\bk}:=(k_1, k_2, k_3) \in \naturels^3$ and $y:=(y_1, y_2, y_3) \in \reels^3$, we set $\abs{{\bk}}:=k_1+k_2+k_3$, ${\bk}!:=k_1!k_2!k_3! $ and  $\widetilde{H}_{{\bk}}(y):= \prod_{\substack{1 \le j \le 3}}H_{k_{j}}(y_j)$.

We denote by $\langle \cdot , \cdot \rangle$ the canonical scalar product in $\reels^2$.\\
${\bf C}$ is a generic constant that could change value while developing a proof.

\subsection{Hypotheses}
\label{sec:hypotheses}

Let $(\Omega, \cal{A}, P)$ be a probability space and $X:   \Omega \times \reels^2 \to \reels$ an affine process, equal in law to $(Z(A.t), t \in \reels^2)$, where $Z: \Omega \times \reels^2 \to \reels$ is isotropic and $A: \reels^2 \to \reels^2$ is a linear self-adjoint transformation. We additionally suppose that $A$ is a positive definite transformation.\\
We also assume  that $Z$ is a centred Gaussian process twice continuously differentiable in $\reels^2$, that is  $Z \in \textbf{C}^{2}(\reels^{2})$.\\
The eigenvalues of $A$ are denoted by $\lambda_1, \lambda_2$, $0 < \lambda_2 \le \lambda_1$.
Let $0 < \lambda \le 1$ be the quotient of the eigenvalues, $\lambda := {\lambda_2}/{\lambda_1}$.
The process $X$ is isotropic means that $\lambda=1$ (see Adler and Taylor \cite[Section 5.7]{MR2319516}).

Let $r_z$ and $r_x$ be the covariance function of $Z$  and $X$ respectively.
The regularity assumption on $Z$ implies that $r_{z} \in \textbf{C}^4(\reels^2)$.\\

Let $T$ be an open bounded rectangle of $\reels^2$. Let also $(T_n)_n$ be open bounded squares with the following form $T_n:= \introo{-n}{n}^2$ with $n \in \naturels^{\star}$, and let $n$ tends to infinity.

\paragraph{Assumptions on the covariance}
For any multidimensional index ${\bm}:= (i_1,$ $\dots, i_k)$ with $1 \le k \le 2$ and $1 \le i_j \le 2$, we write 
$\tfrac{\partial ^{{\bm}}r_{z}}{\partial t^{{\bm}}}(t):= \tfrac{\partial ^{k}r_{z}}{\partial t_{i_{1}}\cdots \partial t_{i_{k}}}(t).$ 
Let 
$$
\Psi(t) := \max\left\{\abs{r_z(t)}, \abs{\dfrac{\partial ^{{\bm}}r_{z}}{\partial t^{{\bm}}}(t)}, {\bm} \in \{1, 2\}^k, 1 \le k \le 2 \right\}.
$$

We make the following assumption.
\begin{itemize}
\item ${\bf H_{\Psi}}$:  $\Psi(t) \to 0$ when $\normp[2]{t} \to +\infty$, $\Psi  \in L^1(\reels^2)$ and $\int_{\reels^2} r_z(t) \ud t$ $>0$.
\end{itemize}

\begin{remark}
\label{spectral positivity}
Note that $r_z \in L^1(\reels^2)$ implies that $r_z \in L^q(\reels^2)$ for all $q \ge 1$ and hence that $Z$ (resp.
X) admits a spectral density $f_z$ (resp.
$f_x$) that is continuous and such that $f_z(0) >0$. We wish to emphasize that this fact implies that the following assumption ${\bf H_{r}}$ is always fulfilled
\begin{itemize}
\item ${\bf H_{r}}$: $r_z^2(0)-r_z^2(t) \neq 0$, for all $t \neq 0$. 
\end{itemize}
\end{remark}
\begin{remark}
\label{norme density}
Note also that $\int_{\reels^2} f_z(t)  \normp[2]{t}^2 \ud t < +\infty$.
\end{remark}

\paragraph{Comments on assumptions made}  Having the property that $f_z(0)>0$  ensures that $\nabla ^2 r_x(0)$ is a non-degenerate matrix. Indeed $-\nabla ^2 r_x(0)= \mu A^2$, where $\mu:= - \dfrac{\partial ^{2}r_{z}}{\partial t_{i}^2}(0)=\frac{1}{2}\int_{\reels^2} f_z(t)  \normp[2]{t}^2 \ud t $, for any $i=1, 2$. Thus $\mu >0$.\\
Hypothesis ${\bf H_{r}}$ being fulfilled ensures that the assumptions of Corollary \ref{gaussien stationnaire} are verified, so that its conclusions are also verified.
That is  for all $u \in \reels$, equality (\ref{vide}) is  almost surely valid (for $T$ or $T_n$). Furthermore if $f: S^{1} \to \reels$ is a continuous and bounded function, functionals $J_{f}(u)$ and $J_{f}^{(n)}(u)$ ($n \in \naturels^{\star}$) are almost surely well defined (see Remark \ref{finitude courbe de Niveau}).  Also if $X$ and $Y$ satisfies the assumption ${\bf H_e}$ with $T$ or $T_n$ then the two-order Rice formula (\ref{formule de Rice d'ordre deux}) is valid for $T$ or $T_n$.\\
The assumption ${\bf H_{\Psi}}$ will be useful to prove Proposition \ref{variance asymptotique xi} in which we exhibit the asymptotic variance of the centred and suitably rescaled general functionals $J_f^{(n)}(u)$ as $n$ tends to infinity. The fact that $f_z(0)>0$ ensures that this limit variance is not trivial.\\

These comments lead us to adopt the following new notations.

\paragraph{Notations} We shall omit in the future the symbol $\textrm{D}^{\textrm{r}}$. Furthermore for writing convenience if $u$ is a fixed level belonging to $\reels$, we will note $\mathcal{C}(u)$ in place of $\mathcal{C}_{\textrm{T}, \textrm{X}}^{\textrm{D}^{\textrm{r}}}(u)$  and $\mathcal{C}_n(u)$ in place of $\mathcal{C}_{\textrm{T}_n, \textrm{X}}^{\textrm{D}^{\textrm{r}}}(u)$. Also, for $B$ a rectangle set included in $T_n$, we will note $\mathcal{C}_B(u)$ in place of $\mathcal{C}_{\textrm{B}, \textrm{X}}^{\textrm{D}^{\textrm{r}}}(u)$.

 \paragraph{Particular isotropic processes} Here we give three examples of isotropic processes satisfying assumption ${\bf H_{\Psi}}$. \\
 Since $Z$ is a stationary Gaussian process, $Z$ will be isotropic if, and only if, its covariance $r_z$ will only depend on the norm.
 Furthermore $Z$ must have a spectral density $f_z$, so this last condition is equivalent to the one that $f_z$ will depend only on the norm. \\
{\bf Powered exponential covariance} The first example is given by a covariance function belonging to the class of powered exponential family (see Yaglom \cite[Chap. 4.22.2, example 3, p. 364]{MR0893393}), that is: $r_z(t):= C\,  \exp(-\alpha \normp[2]{t}^2)$, $t \in \reels^2$, $C>0$, $\alpha >0$ being a scale parameter.
One has, $\int_{\reels^2} r_z(t) \ud t = \frac{C \pi}{\alpha}>0$. Here the spectral density $f_z$ is $f_z(\lambda):= \frac{C}{4 \pi \alpha}  \exp(-\frac{1}{4 \alpha}\normp[2]{\lambda}^2)$, $\lambda \in \reels^2$.  The covariance $r_z$ satisfies assumption ${\bf H_{\Psi}}$.

\noindent {\bf Generalized Cauchy covariance}
The second example of covariance function is given by the one belonging to the Cauchy family (see \cite[Chap. 4.22.2, example 5, p. 365]{MR0893393}, see also Anderson \cite{MR1179637} and Lim and Teo \cite[equality (2.1) p. 1326]{MR2508576}), that is $r_z(t):= C\,(\alpha^2 +\normp[2]{t}^2)^{-\nu}$, $t \in \reels^2$, $C>0$, $\alpha >0$ being a scale parameter and $\nu >0$ a smoothness parameter controlling the long range dependence of the process $Z$. The spectral density is 
$$
f_z(\lambda):=\frac{C \alpha^{1-\nu}\normp[2]{\lambda}^{\nu -1}}{2^{\nu} \pi \Gamma(\nu)} K_{\nu -1}(\alpha \normp[2]{\lambda})
$$
(see Yaglom \cite[Chap. 4, equality (4.12'), p. 129]{MR09155577}, see also Lim and Teo \cite[p. 3]{MR2485858}). Here $K_{\nu}$ is the modified Bessel function of the third kind (see Gradshteyn and Ryzhik \cite[8.432 p. 917]{MR2360010} for definition) and $\Gamma$ is the Gamma function defined in Appendix (\ref{fonction Gamma}).\\
If we add the condition $\nu >1$, then $\int_{\reels^2} r_z(t) \ud t = \frac{C \pi \alpha^{2-2\nu}}{\nu -1}>0$, and covariance $r_z$ satisfies assumption ${\bf H_{\Psi}}$. 

\noindent {\bf Whittle-Mat\'ern covariance}
The third example of covariance function is given by the one belonging to the Mat\'ern class (see \cite[Chap. 4.22.2, example 2, p. 363]{MR0893393}, see also \cite[equality (1) p. 2]{MR2485858}).
The covariance function is given by $r_z(t):= C(\alpha \normp[2]{t})^{\nu} K_{\nu}(\alpha \normp[2]{t})$, $t \in \reels^2$, $C>0$, $\alpha >0$ being a scale parameter controlling the spatial range of the covariance and  $\nu >0$ is the smoothness parameter governing the level of smoothness of $Z(t)$.
The spectral density is $f_z(\lambda):= \frac{C 2^{\nu -1} \Gamma(\nu +1) \alpha^{2\nu}}{\pi (\alpha^2+ \normp[2]{\lambda}^2)^{\nu +1}}$ (see  \cite[equality (4.131) p. 364, Chap. 4.22.2, example 2]{MR0893393} and also \cite[Proposition 3.1 p. 1333]{MR2508576}).
Since $\nu > -1$, using the spectral density representation of the covariance in zero in the second example one can prove that $\int_{\reels^2} r_z(t) \ud t = C \pi \alpha^{-2}2^{\nu +1} \Gamma(\nu +1)>0$. In addition, if we suppose that $\nu >2$, one can easily show that $\int_{\reels^2} \normp[2]{\lambda}^4 f(\lambda) d\lambda <+\infty$, so that the covariance function $r_z  \in \textbf{C}^4(\reels^2)$ (one can also refer to \cite[(33) and (34), p. 9]{MR2485858}).
By using \cite[equalities 10 and 11 of 8.486(1) p. 929]{MR2360010} we obtain the expression of the derivatives of function $r_z$. In \cite{MR2360010}, equalities 8.485 p. 928 and 8.445 and 8. 446 p. 919 give the asymptotic behaviour of $K_{\nu}(z)$ when $z \to 0$  and  6. of 8.451 p. 920 the one when $z \to +\infty$.  With these properties one can establish in case where $\nu >2$ that $\Psi(t) \to 0$ when $\normp[2]{t} \to +\infty$ and that  $\Psi  \in L^1(\reels^2)$. In conclusion we deduce that when $\nu >2$ the covariance $r_z$ satisfies assumption ${\bf H_{\Psi}}$.

\paragraph{Commentaries about these examples} In the two first examples we remark that covariance $r_z$ is strictly decreasing as function of the distance from the origin. If the linear transformation $A$ is not an isometry, this property vanishes since points living in a circle centred at the origin with same covariance $r_z$ will be transformed by $A$ into points over an ellipse and thus with different covariances $r_x$. In other words and in a general way, meaning not just for the two previous examples, circles are deformed by stretching via the transformation $A$. A point $t \in \reels^2$ such that $\normp[2]{t}=1$, that is living in the unity circle, will be such that $A \cdot t$ will belong to the centred ellipse at the origin defined by the $A$-eigenvectors directions $(v_1, v_2)$ as axis and by the eigenvalues $\lambda_1$, $\lambda_2$ the magnitude of the axis. Measuring the range of this stretching will be part of the matter of the paper, through the construction of estimators of the ratio of the eigenvalues, say $\lambda={\lambda_2}/{\lambda_1}$, and of $\theta$, the angle rotation between the eigenvectors $(v_1, v_2)$ and a fixed direct orthonormal basis $(v^{\star}, v^{\star \star})$ (see Figure \ref{fig:figure1} page \pageref{fig:figure1}).

\section{Hermite expansion for level functionals}
\label{fonctionnelles vues dans le chaos}
In this section $T$ will be a fixed open bounded rectangle of $\reels^2$.\\
For a continuous and bounded function $f$, $f: S^{1} \to \reels$, recall that we have defined in section \ref{sec:notations} the following general functional $J_{f}(u)$ of the level $u$ as:
\begin{equation*}
J_{f}(u):= \frac{1}{\sigma_{2}(T)} \int_{\mathcal{C}(u)} f(\nu_{X}(t)) \ud\sigma_1(t),
\end{equation*}
where $\nu_{X}(t):= {\frac{\nabla X(t)}{\normp[2]{\nabla X(t)}}}$.\\
Our objective is for $u$ a fixed level to give an Hermite expansion of the random variable $J_{f}(u)$.\\
In this aim the idea consists in approaching functionals $J_{f}(u)$ by other functionals, say $J_{f}(u, \sigma)$ ($\sigma \to 0$), in such a way that the last ones can be expressed as stochastic integrals with respect to Hermite polynomials. This is possible through the use of a kernel $K_{\sigma}$ and via the coarea formula. Then the idea will consist in proving that $J_{f}(u, \sigma)$ tends in $L^2(\Omega)$ towards $J_{f}(u)$ as $\sigma \to 0$.\\
We emphasize that the convergence is not an obvious purpose since the difficulty is to prove the continuity of function $({\bf{y_1}}, {\bf{y_2}}) \mapsto \mathbb{E}[J_{f}({\bf{y_1}})\, J_{f}({\bf{y_2}})]$ that will be the aim of the next paragraph. This fact is far from being trivial and requires a number of results stated in section \ref{sec:notations}, as the second order Rice formula, and also studied in next section, as a in-depth study of the level curve $\mathcal{C}(u+\delta)$ in a neighbourhood of the level $u$.\\
By using the $L^2$-continuity in the level $u$ of functional $J_{f}(u)$, we will obtain as by product an expression for its first moment, that is a first order Rice formula for all level.
\subsection{$L^2$-continuity in the level}
\label{continuity in the level}
We prove the following theorem.

\begin{theorem}
\label{continuite longueur de courbe}
Let $f: S^{1} \to \reels$ be a continuous and bounded function.
Then the function
\begin{align*}
({\bf{y_1}}, {\bf{y_2}}) \mapsto \mathbb{E}[J_{f}({\bf{y_1}}) J_{f}({\bf{y_2}})],
\end{align*}
is continuous.
In particular the same holds for function,
\begin{align*}
\textbf{y} \mapsto  \mathbb{E}[J_{f}(\textbf{y})]^2.
\end{align*}
\end{theorem}
\begin{proofarg}{Proof of Theorem \ref{continuite longueur de courbe}}
As explained in Remark \ref{noncontinuite}, process $\widetilde{Y}:= f(\nu_X) \1_{\textrm{D}^{\textrm{r}}_{\textrm{X}}}$ is not continuous on $T$, so we can not apply Corollary \ref{gaussien stationnaire} to compute $\mathbb{E}[J_{f}(\textbf{y})]^2$. We need to be a little careful. In fact the following lemma will prove as by product that such computation holds for $\mathbb{E}[J_{f}(\textbf{y})]^2$ and that function $\textbf{y} \mapsto  \mathbb{E}[J_{f}(\textbf{y})]^2$ is continuous.\\
Let us state this lemma proved in Appendix and for which the proof is based on the application of  Corollary \ref{gaussien stationnaire}.
\begin{lemma}
\label{continuite Y}
Let  $F: T \times (\reels^{2})^{\star} \to \reels$ be a bounded continuous function of its arguments, then
\begin{align*}
\textbf{y} \mapsto  \mathbb{E}[\int_{\mathcal{C}(\textbf{y})} F(t, \nabla X(t)) \1_{\textrm{D}^{\textrm{r}}_{\textrm{X}}}(t) \ud \sigma_1(t)]^2
\end{align*}
is a continuous function, where we have noted $(\reels^2)^{\star}:=\{v \in \reels^2, \normp[2]{v} \neq 0\}$.\\
Furthermore, for any $u \in \reels$ one has
\begin{align*}
\lefteqn{\mathbb{E}[\int_{\mathcal{C}(u)} F(t, \nabla X(t)) \1_{\textrm{D}^{\textrm{r}}_{\textrm{X}}}(t) \ud \sigma_1(t)]^2}\\ \nonumber
&= \int_{T \times T}
	\mathbb{E} \left[F(t_1, \nabla X(t_1)) \1_{\textrm{D}^{\textrm{r}}_{\textrm{X}}}(t_1)\, F(t_2, \nabla X(t_2)) \1_{\textrm{D}^{\textrm{r}}_{\textrm{X}}}(t_2) \normp[2]{\nabla
X(t_1)} \normp[2]{\nabla
X(t_2)}|X(t_1)=X(t_2)=u\right]  \\ \nonumber
&\qquad\qquad\qquad \,\times \hfill{p}_{X(t_1), X(t_2)}(u, u)\ud t_1\ud t_2,
\end{align*}
where ${p}_{X(t_1), X(t_2)}(u, u)$ stands for the density of $(X(t_1), X(t_2))$ in point $(u, u)$.\end{lemma}

Let $y$ be fixed in $\reels$ and $(y_k)_{k \in \naturels}$ be a real sequence tending to $y$. \\
Let us note $
I_{f}(y):=  \int_{\mathcal{C}(y)} f(\nu_{X}(t)) \ud \sigma_1(t)
$.\\
Proving Theorem \ref{continuite longueur de courbe} is equivalent to show that $\normp[L^2(\Omega)]{I_{f}(y_k)- I_{f}(y)}$ tends to zero as $k$ goes to infinity, where  $\normp[L^2(\Omega)]{\cdot}$ stands for the $L^2(\Omega)$-norm.\\
Thus we try to apply the Scheff\'e's lemma but it is not so direct. \\
Indeed Lemma \ref{continuite Y} applied to function $F: \reels^2 \times (\reels^{2})^{\star} \to \reels$, such that $F(t, z):= f(z/\normp[2]{z})\1_{\{z \neq 0 \}}(z)$ implies that $\normp[L^2(\Omega)]{I_{f}(y_k)}$ tends to $\normp[L^2(\Omega)]{I_{f}(y)}$ as $k$ goes to infinity. But we need the missing following convergence, that is, almost surely
$
I_{f}(y_k)$ tends to $ I_{f}(y)$, with $k$.\\
Our goal is now to tempt to apply the following theorem stated and proved as \cite[Theorem 3.1.2]{Berzin:2017}. Note that in its statement the functions of interest are deterministic, that is they are not random functions.
\begin{theorem}
\label{Rice continuite}
Let $g:  T \subset \reels^2 \to \reels$ be a function belonging to $\textbf{C}^1(T)$ such that $\nabla g$ is Lipschitz on $T$. Let $h:T \subset \reels^2 \to \reels$ be a continuous function on $T$ such that $\supp(h) \subset \textrm{D}^{\textrm{r}}_{\textrm{g}}:=\left\{t \in
T :  \normp[2]{\nabla g(t)} \neq 0\right\}$, $\supp(h)$ being the support of function $h$. Then the function
\begin{align*}
\textbf{y} \to \int_{\mathcal{C}_{\textrm{T}, \textrm{g}}^{\textrm{D}^{\textrm{r}}}(\textbf{y})} h(t) \ud \sigma_1(t)
\end{align*}
is continous with respect to the variable $\textbf{y}$. \\
We have noted $\mathcal{C}_{\textrm{T}, \textrm{g}}^{\textrm{D}^{\textrm{r}}}(y):=\left\{t \in
T : g(t)=y, \normp[2]{\nabla g(t)} \neq 0 \right\}$.
\end{theorem}
Remark that process $X$ is such that $\nabla X$ is almost surely Lipschitz on $T$. Indeed we have already mentioned after the statement of Corollary \ref{gaussien stationnaire} that assumption ${\bf H_d}$ remains checked. Thus equality (\ref{sup esperance}) is valid, that is $
\mathbb{E}[\sup_{t \in T} \normp[1, 2]{\nabla^2X(t)}^{(s)}]^{4} < +\infty$. So we can conclude that almost surely $L_X:= \sup_{t \in T} \normp[1, 2]{\nabla^2X(t)}^{(s)} < +\infty$. Furthermore the Taylor formula applied to $X$ belonging to $\textbf{C}^2(T)$  ensures that almost surely $\nabla X$ is Lipschitz with Lipschitz constant $L_X$.\\
As already mentioned in section \ref{sec:notations} process $\widetilde{Y}:= f(\nu_X) \1_{\textrm{D}^{\textrm{r}}_{\textrm{X}}}$ is not continuous on $T$ and there is no reason for why $\supp(\widetilde{Y}) \subset \textrm{D}^{\textrm{r}}_{\textrm{X}}$.\\
In view of applying the Scheff\'e's lemma the idea consists in approximating the functional $I_{f}(y)$ by the following one. It will check hypotheses of Theorem \ref{Rice continuite} and also that of Lemma \ref{continuite Y} .
For $m \in \naturels^{\star}$, let 
\begin{align*}
I_{f, m}(y):=\int_{\mathcal{C}(y)} Y_m(t) \ud \sigma_1(t),
\end{align*}
where we define for $t \in T$
\begin{align*}
Y_m(t):= F_m(t, \nabla X(t)) \1_{\textrm{D}^{\textrm{r}}_{\textrm{X}}}(t),
\end{align*}
with function $F_m: T \times (\reels^{2})^{\star} \to \reels$, defined by
\begin{align*}
F_m(t, z):= g_m(t)\, \varphi\left({\frac{1}{m \normp[2]{z}}}\right) f(z/\normp[2]{z} ),
\end{align*}
where $\varphi$ is a continuous decreasing function on $\reels^{+}$, such that
\begin{equation*}
\varphi(t):=
\begin{cases}
1, & \mbox{if $0 \le t \le 1$} \\
0,  &\mbox{if $2 \le t $}
\end{cases}
\end{equation*}
and $(g_m)_{m \in \naturels^{*}}$ is a sequence of functions defined on $\reels^2$ to $\intrff{0}{1}$ in the following manner
\begin{equation*}
g_m( x ) := \frac{d( x , T^{2m})}{d( x , T^{2m})+ d( x , T^{(m)})},
\end{equation*}
where $d( x ,~B)$ stands for the distance between the point $ x $ and the set $B\subset\reels^2$.
The closed sets $T^{2m}$ and $T^{(m)}$ are defined by
\begin{equation*}
T^{2m}:=\left\{ x  \in \reels^2, d( x , T^{c}) \le \tfrac{1}{2m}\right\}  \mbox{ and } T^{(m)}:=\left\{ x  \in \reels^2, d( x , T^{c}) \ge \tfrac{1}{m}\right\},
\end{equation*}
$T^c$ denoting the complement of $T$ on $\reels^2$.\\
We have shown in \cite[Lemmas 3.2.1 and 3.2.2]{Berzin:2017} that the functions $(g_m)_{m \in \naturels^{\star}}$ are well defined, continuous. Thus using the crucial central term in $\varphi$ defining $Y_m$, $m \in \naturels^{\star}$, one can prove that the r.v. $Y_m$ is continuous on $T$. Furthermore we established that the support of function $g_m$ restricted to the set $T$ is contained in $T$ for each $m \in \naturels^{\star}$. This fact implies that
$\supp(Y_m) \subset \textrm{D}^{\textrm{r}}_{\textrm{X}}$.\\
At last this sequence $(g_m)_{m \in \naturels^{\star}}$ is bounded by one and tends to one when $m$ goes to infinity. This implies that the sequence $(Y_m)_{m \in \naturels^{\star}}$ is bounded on $T$ and that  almost surely for all $t \in T$, $Y_m(t) \xrightarrow[m \to +\infty]{} \widetilde{Y}(t)=f(\nu_{X}(t)) \1_{\textrm{D}^{\textrm{r}}_{\textrm{X}}}(t)$.\\
 By summarizing, we have
 \begin{enumerate}
 \item $g:=X$ and $h:=Y_m$ verifies assumptions of Theorem \ref{Rice continuite}
\item $F_m$ and $F -F_m$ verify assumptions of Lemma \ref{continuite Y} 
\item Almost surely for all $t \in T$, $Y_m(t) \xrightarrow[m \to +\infty]{} \widetilde{Y}(t)$ 
\end{enumerate}

Assertion 1. implies  that for all $m \in \naturels^{\star}$, almost surely
\begin{align*}
I_{f, m}(y_k) \xrightarrow[k \to +\infty]{} I_{f, m}(y).
\end{align*}
Note that this almost sure convergence is what we need to move forward.
Furthermore, assertion 2. implies that for all $m \in \naturels^{\star}$,
\begin{align*}
\normp[L^2(\Omega)]{I_{f, m}(y_k)} \xrightarrow[k \to +\infty]{} \normp[L^2(\Omega)]{I_{f, m}(y)},
\end{align*}
so that Scheff\'e's lemma allows to conclude that for all $m \in \naturels^{\star}$,
\begin{align}
\label{egalite (5)}
\lim_{k \to +\infty}\normp[L^2(\Omega)]{I_{f, m}(y_k)- I_{f, m}(y)} = 0.
\end{align}
Now, assertion 2. also implies that for all $m \in \naturels^{\star}$,
\begin{align}
\label{egalite (6)}
 \lim_{k \to +\infty}  \normp[L^2(\Omega)]{I_{f}(y_k)- I_{f, m}(y_k)}=
  \normp[L^2(\Omega)]{I_{f}(y)- I_{f, m}(y)}.
\end{align}
An upper bound is 
\begin{multline*}
\normp[L^2(\Omega)]{I_{f}(y_k)- I_{f}(y)}  \le \normp[L^2(\Omega)]{I_{f}(y_k)- I_{f, m}(y_k)} + \\
\normp[L^2(\Omega)]{I_{f, m}(y_k)- I_{f, m}(y)} +\normp[L^2(\Omega)]{I_{f, m}(y)- I_{f}(y)}.
\end{multline*}
Applying convergence obtained in (\ref{egalite (5)}) and (\ref{egalite (6)}), it yields, for all $m \in \naturels^{\star}$,
\begin{align}
\label{inegalite I}
\underset{k \to +\infty}{\overline{\lim}} \normp[L^2(\Omega)]{I_{f}(y_k)- I_{f}(y)}  \le 2  \normp[L^2(\Omega)]{I_{f}(y)- I_{f, m}(y)}.
\end{align}
Let us set,
\begin{align*}
I_{f}(y)- I_{f, m}(y)=\int_{\mathcal{C}(\textbf{y})} \widetilde{Y_m}(t) \ud \sigma_1(t),
\end{align*}
where $\widetilde{Y_m}(t):= \widetilde{Y}(t)-Y_m(t)$.\\
By using convergence given in point 3., we obtain that
almost surely for all $t \in \textrm{D}^{\textrm{r}}_{\textrm{X}}$, $\lim_{m \to +\infty} \widetilde{Y}_m(t)=0$.\\
Furthermore for all $t \in \textrm{D}^{\textrm{r}}_{\textrm{X}}$, one has $\abs{\widetilde{Y}_m(t)} \le {\bf C}$.
Moreover by applying  Remark \ref{finitude courbe de Niveau}, one has
\begin{align}
\label{finitude E}
 \mathbb{E}[\sigma_{1}(\textrm{C}(y))]^2 < +\infty,
\end{align} 
and then almost surely $\sigma_{1}(\textrm{C}(y)) < +\infty$.
The Lebesgue convergence theorem induces that 
\begin{itemize}
\item Almost surely, $\lim_{m \to +\infty} (I_{f}(y)- I_{f, m}(y))=0$.
Furthermore
\item Almost surely, $\abs{I_{f}(y)- I_{f, m}(y)} \le {\bf C}\, \sigma_{1}(\textrm{C}(y)) \in L^{2}(\ud P)$,
last belonging providing from (\ref{finitude E}).\
\end{itemize}
Applying once again the Lebesgue convergence theorem, one has finally proved that $\lim_{m \to +\infty} \normp[L^2(\Omega)]{I_{f, m}(y)- I_{f}(y)}=0$.
Thus by using inequality (\ref{inegalite I}), one concludes that
\begin{align*}
\underset{k \to +\infty}{\overline{\lim}} \normp[L^2(\Omega)]{I_{f}(y_k)- I_{f}(y)}  =0,
\end{align*}
yielding Theorem \ref{continuite longueur de courbe}.
\end{proofarg}

We are now ready to define a $L^2(\Omega)$-expansion for $J_f(u)$.

\subsection{Hermite expansion for general functionals}
\label{chaos section}
Let $f: S^{1} \to \reels$ be a continuous and bounded function and $u$ a fixed level in $\reels$. We define an approximation of $J_{f}(u)$.
For $\sigma >0$, we define 
\begin{equation*}
J_{f}(u, \sigma):= \frac{1}{\sigma} \int_{-\infty}^{+\infty} K(\tfrac{u-v}{\sigma})\, J_{f}(v) \ud v,
\end{equation*}
where $K$ is a continuous density function with a compact support in $[-1, 1]$. In view of applying the coarea formula to 
$J_{f}(u, \sigma)$ we
recall a statement of this formula, adapted to our situation and that can be found in Federer \cite[Theorem 3.2.12]{MR0257325} and also in Berzin et al. \cite[Corollary 2.1.1]{Berzin:2017}. 

\begin{corollary} (\it{Coarea formula})\\
\label{coarea formula}
Let $h$ a measurable function, $h: \reels^2 \times \reels \to \reels$ and $g: T \subset \reels^2 \to \reels$ be a function belonging to $\textbf{C}^1(T)$.
We have
\begin{align*}
\int_{T} h(t, g(t)) \normp[2]{\nabla g(t)} \ud t = \int_{-\infty}^{+\infty} \left[{\int_{\mathcal{C}_{\textrm{T}, \textrm{g}}^{\textrm{D}^{\textrm{r}}}(y)}h(t, y) \ud \sigma_{1}(t)}\right]  \ud y,
\end{align*}
provided that one of the two integrals is finite.
\end{corollary}
By applying last corollary to the measurable function $h: \reels^2 \times \reels \to \reels$ defined as :
\begin{align*}
h(t, y):= \frac{1}{\sigma_2(T)}f(\nu_{X}(t))\, \1_{\{ \nabla X(t) \neq 0\}}\,\frac{1}{\sigma}\, K(\tfrac{u-y}{\sigma}),
\end{align*}
to the $C^1$-function $g:=X: T \subset \reels^2 \to \reels$, we get:
\begin{align*}
J_{f}(u, \sigma)=  \frac{1}{\sigma_{2}(T)} \int_{T} f(\nu_{X}(t)) \, \1_{\{ \nabla X(t) \neq 0\}}\, \frac{1}{\sigma}\, K(\tfrac{u-X(t)}{\sigma})\, \normp[2]{\nabla X(t)}\, \ud t < +\infty.
\end{align*}
Let
$\Lambda:= \begin{pmatrix}
\lambda_1 & 0\\
0 & \lambda_2
\end{pmatrix}$, 
a matrix that factorizes $A$ such that $A:=P \Lambda P^t$, where $P$ is a unitary matrix and $^t$ stands for the transpose symbol. 

The matrix $\mathit{\Gamma}^{\bsX}$ stands for the covariance matrix of the $3$-dimensional Gaussian vector 
\begin{equation*}
\bsX(t):= \left({\nabla X(t), X(t) }\right)^t.
\end{equation*}
We recall that since $Z$ is isotropic, $\mathbb{E}[\nabla Z(0) \nabla Z(0)^t]= \mu\, I_2$, where $\mu:= - \dfrac{\partial ^{2}r_{z}}{\partial t_{i}^2}(0)$ for any $i=1, 2$ and $I_2$ stands for the identity matrix in $\reels^2$.

We denote $\Delta$ for the $3\times3$ matrix,
$
\Delta:=
\begin{pmatrix}
\sqrt{\mu} P \Lambda&0\\
0 & \sqrt{r_z(0)},
\end{pmatrix}$.
It is such that $\Delta \Delta^t= \mathit{\Gamma}^{{\bsX}}$.
We have seen in Section \ref{sec:hypotheses} (see Remark \ref{spectral positivity}) that $Z$ admits a spectral density $f_z$ and then the same fact occurs for $X$.
We have noted $f_x$ the spectral density for $X$.\\
We obtain the following spectral representations: for $t \in \reels^2$,
\begin{align*}
X(t)= \int_{\reels^2} e^{i \,\langle t ,\, \lambda \rangle} \sqrt{f_x(\lambda)}\ud W(\lambda),
\end{align*}
where $W$ stands for the standard Brownian motion.\\
Thus, for any $\lambda:=(\lambda_i)_{1 \le i \le 2}$ in $\reels^2$, we let
\begin{align*}
\nu(\lambda):=(i\lambda_1, i \lambda_2, 1), 
\end{align*}
so that for any $t \in \reels^2$,
\begin{align*}
U(t)= (U_{i}(t))_{1 \le i \le 3}:=\Delta^{-1}.{\bsX}(t)=\int_{\reels^2} e^{i \,\langle t ,\, \lambda \rangle} \sqrt{f_x(\lambda)}\, \Delta^{-1}. \nu(\lambda)\ud W(\lambda).
\end{align*}

We can write for any $t \in \reels^2$, ${\bsX}(t)= \Delta.U(t)$, where $U(t)$ is a $3$-dimensional standard Gaussian vector.\\
With these notations, one obtains:
\begin{multline*}
J_{f}(u, \sigma) =  \frac{\sqrt{\mu}}{\sigma_{2}(T)} \frac{1}{\sigma}
	\int\limits_{T}\,f\left({\frac{P\Lambda.
		\begin{pmatrix} U_1(t)\\
		U_2(t)
		\end{pmatrix}}
		{\sqrt{\lambda_1^2U_1^2(t) + \lambda_2^2 U_2^2(t)}}}\right) \1_{\{U_1(t)\neq 0 \, or \,U_2(t)\neq 0\}}\\
		\times K\left({\frac{u-\sqrt{r_z(0)}\,U_{3}(t)}{\sigma}}\right)
\sqrt{\lambda_1^2U_1^2(t)+\lambda_2^2U_2^2(t)}\, \ud t.
\end{multline*}
We define for $\sigma >0$ and $y=(y_i)_{1\le i \le 3} \in \reels^3$, the map $g_{\sigma}$ as:
\begin{align*}
g_{\sigma}(y):= \frac{\sqrt{\mu}}{\sigma} K\left({\frac{u-\sqrt{r_z(0)}\,y_{3}}{\sigma}}\right)\, f\left({\frac{P\Lambda.\begin{pmatrix} y_1\\y_2 \end{pmatrix}}{\sqrt{\lambda_1^2y_1^2+\lambda_2^2y_2^2}}}\right) 
\sqrt{\lambda_1^2y_1^2+\lambda_2^2y_2^2}.
\end{align*}
Since the map belongs to $L^2(\reels^3, \phi_3(y)\ud y)$, the following expansion converges in this space:
\begin{align*}
g_{\sigma}(y)= \sum_{q=0}^{+\infty} \sum_{
\substack{
{\bk} \in \naturels^3 \\ \abs{{\bk}}=q
}
} a_{f, \,\sigma}({\bk}, u)\, \widetilde{H}_{{\bk}}(y),
\end{align*}
while taking ${\bk}:=(k_i)_{1 \le i \le 3} \in \naturels^3$,
\begin{align}
\label{afsigma}
a_{f, \,\sigma}({\bk}, u):= a_{f}(k_1, k_2)\, a_{\sigma}(k_3, u),
\end{align}
where for $(k_1, k_2) \in \naturels^2$,
\begin{multline}
\label{a f k}
a_f(k_1, k_2)
	:= \frac{\sqrt{\mu}}{k_1!k_2!} 
		\int\limits_{\reels^2} f\left({\frac{P\Lambda. \begin{pmatrix} y_1\\y_2 \end{pmatrix}}{\sqrt{\lambda_1^2y_1^2+\lambda_2^2y_2^2}}}\right) 
\sqrt{\lambda_1^2y_1^2+\lambda_2^2y_2^2} \, H_{k_1}(y_1) \,\phi(y_1)\\
H_{k_2}(y_2) \,\phi(y_2) \ud y_1 \ud y_2\
\end{multline}
and
\begin{align*}
a_{\sigma}(k_3, u):=\frac{1}{k_3!} \int\limits_{-\infty}^{+\infty} \frac{1}{\sigma} \,K\left({\tfrac{\ds u-\sqrt{r_z(0)}\,y}{\ds \sigma}}\right)\,H_{k_3}(y)\, \phi(y)\ud y.
\end{align*}
In this way, we obtain the expansion of the functional $J_{f}(u, \sigma)$ in $L^2(\Omega)$, that is:
\begin{align}
\label{Jfnusigma}
J_{f}(u, \sigma) \stackrel{L^2(\Omega)}{=} \frac{1}{\sigma_{2}(T)} \sum_{q=0}^{+\infty} \sum_{
\substack{
{\bk} \in \naturels^3 \\ \abs{{\bk}}=q
}
} a_{f, \,\sigma}({\bk}, u)\, \int_{T} \widetilde{H}_{{\bk}}(U(t))\, \ud t.
\end{align}
Now observe that as $\sigma$ tends to zero, $a_{\sigma}(k_3, u) \to a(k_3, u)$,  where coefficient $a(k_3, u)$ is defined as
\begin{align}
\label{ak3}
a\left(k_3, u\right):=\frac{1}{k_3!} H_{k_3}\left(\tfrac{\ds u}{\ds \sqrt{r_z\left(0\right)}}\right)\, \phi\left(\tfrac{\ds u}{\ds \sqrt{r_z\left(0\right)}}\right)\, \tfrac{\ds 1}{\ds \sqrt{r_z\left(0\right)}}.
\end{align}
 This remark will yield via the forthcoming theorem the following expansion in $L^2(\Omega)$ of $J_{f}(u)$.
 \begin{theorem}
\label{equality L2 Bis}
For $f: S^{1} \to \reels$ a continuous and bounded function and $u$ a fixed level in $\reels$, one has the following expansion in $L^2(\Omega)$,
\begin{align*}
J_{f}(u) = \frac{1}{\sigma_{2}(T)} \sum_{q=0}^{+\infty} \sum_{
\substack{
{\bk} \in \naturels^3 \\ \abs{{\bk}}=q
}
} a_{f}({\bk}, u)\, \int_{T} \widetilde{H}_{{\bk}}(U(t))\, \ud t,
\end{align*}
where coefficients $a_{f}({\bk}, u)$ are defined by
\begin{align}
\label{af}
a_f({\bk}, u):= a_f(k_1, k_2) \,a(k_3, u),
\end{align}
with $a_f(k_1, k_2)$ and $a(k_3, u)$ previously respectively defined  by (\ref{a f k}) and (\ref{ak3}).
\end{theorem}
 
\begin{remark}
\label{serie dirac}
This equality has nothing obvious because we do not know a priori that the right member in the last equality really belongs to $L^2(\Omega)$.
This fact comes from the way we obtained
this expansion, the series $\sum_{k_{3}=0}^{+\infty} a^2(k_{3}, u) k_{3}!$ being equal to $+\infty$ as the Hermite development in $L^2(\reels, \phi(x) \ud x)$ of delta's Dirac function in point $u$.
\end{remark}
\begin{proofarg}{Proof of Theorem \ref{equality L2 Bis}}
The proof of this theorem strongly leans on that given  in Estrade and Le\'on \cite{MR3572325} and extensively uses the orthogonality of the various chaos.\\
We demonstrate Proposition \ref{convergence sigma}.
\begin{proposition}
For $f: S^{1} \to \reels$ a continuous and bounded function and $u$ a fixed level in $\reels$,
one has the following convergence,
\label{convergence sigma}
\begin{align*}
J_{f}(u, \sigma) \xrightarrow[\delta \to 0]{L^2(\Omega)} J_{f}(u).
\end{align*}
\end{proposition}
\begin{proofarg}{Proof of Proposition \ref{convergence sigma}}
One can easily see that the proof of this proposition follows immediately from that of Theorem \ref{continuite longueur de courbe}.
\end{proofarg}
As in the proof of \cite[Proposition 1.3]{MR3572325}, let us define formally
\begin{align*}
\eta(T):=\frac{1}{\sigma_{2}(T)} \sum_{q=0}^{+\infty} \sum_{
\substack{
{\bk} \in \naturels^3 \\ \abs{{\bk}}=q
}
} a_{f}({\bk}, u)\, \int_{T} \widetilde{H}_{{\bk}}(U(t))\, \ud t,
\end{align*}
and prove the following lemma.
\begin{lemma}
\label{eta T}
$\eta(T) \in L^{2}(\Omega)$.
\end{lemma}
\begin{proofarg}{Proof of Lemma \ref{eta T}}
First, remark that if ${\bk}, {\bm} \in \naturels^3$ are such that $\abs{{\bk}} \neq \abs{{\bm}}$, then for all $s, t \in T \times T$,
$\mathbb{E}[\widetilde{H}_{{\bk}}(U(t)) \widetilde{H}_{{\bm}}(U(s))]=0$ and the above expression, $\eta(T)$ turns out to be a sum of orthogonal terms in $L^2(\Omega)$.

Indeed, to prove this we need a generalization of Mehler's formula given in Aza{\"\i}s and Wschebor \cite[Lemma 10.7, part (b), page 269]{MR2478201} via the following lemma proved in Appendix \ref{ann:proof}.
\begin{lemma}
\label{mehler}
Let $X=(X_i)_{i=1,\,2,\,3}$ and $Y=(Y_j)_{j=1,\,2,\,3}$ be two centred standard Gaussian vectors in $\reels^3$ such that for $1 \le i, j \le 3$, $\mathbb{E}[X_iY_j]=\rho_{ij}$, then for ${\bk}, {\bm} \in \naturels^3$, one has
\begin{align*}
\mathbb{E}[\widetilde{H}_{{\bk}}(X) \widetilde{H}_{{\bm}}(Y)]=
\bigg(\sum_{
\substack{
d_{ij} \ge 0 \\ \sum_{j}d_{ij}=k_i \\ \sum_{i}d_{ij}=m_j
}
} {\bk}!  {\bm}! \prod_{1 \le i, j \le 3} \frac{\rho_{ij}^{d_{ij}}}{d_{ij}!}\,\,\bigg)\, \1_{\abs{{\bk}}=\abs{{\bm}}}.
\end{align*}
\end{lemma}
As in \cite{MR3572325}, let us fix $Q \in \naturels$ and let us denote by $\pi^{Q}$ the projection onto the first $Q$ chaos in $L^2(\Omega)$ and by $\pi_Q$ the projection onto the remaining one.
With these notations one has
\begin{align*}
\pi^Q(\eta(T))= \frac{1}{\sigma_{2}(T)} \sum_{q=0}^{Q} \sum_{
\substack{
{\bk} \in \naturels^3 \\ \abs{{\bk}}=q
}
} a_{f}({\bk}, u)\, \int_{T} \widetilde{H}_{{\bk}}(U(t))\, \ud t,
\end{align*}
and $\pi_Q(\eta(T))=\eta(T)- \pi^Q(\eta(T))$.\\
By the precedent remark, we have
\begin{align*}
\mathbb{E}[\pi^Q(\eta(T))]^2= \sum_{q=0}^Q \mathbb{E}[\frac{1}{\sigma_{2}(T)} \sum_{
\substack{
{\bk} \in \naturels^3 \\ \abs{{\bk}}=q
}
} a_{f}({\bk}, u)\, \int_{T} \widetilde{H}_{{\bk}}(U(t))\, \ud t]^2.
\end{align*}
Remember that coefficients $a_{f, \sigma}({\bk}, u)$ have been defined by (\ref{afsigma}).

Applying the Fatou's lemma, since $\lim_{\sigma \to 0}a_{f, \sigma}({\bk}, u)=a_{f}({\bk}, u)$, we obtain 
\begin{align*}
\mathbb{E}[\pi^Q(\eta(T))]^2 &\le \mathop {\underline{\lim}}\limits_{\sigma \to 0} \sum_{q=0}^Q \mathbb{E}[\frac{1}{\sigma_{2}(T)} \sum_{
\substack{
{\bk} \in \naturels^3 \\ \abs{{\bk}}=q
}
} a_{f, \sigma}({\bk}, u)\, \int_{T} \widetilde{H}_{{\bk}}(U(t))\, \ud t]^2\\
&\le  \mathop {\underline{\lim}}\limits_{\sigma \to 0} \sum_{q=0}^{+\infty} \mathbb{E}[\frac{1}{\sigma_{2}(T)} \sum_{
\substack{
{\bk} \in \naturels^3 \\ \abs{{\bk}}=q
}
} a_{f, \sigma}({\bk}, u)\, \int_{T} \widetilde{H}_{{\bk}}(U(t))\, \ud t]^2\\
&= \mathop {\underline{\lim}}\limits_{\sigma \to 0} \mathbb{E}[
\frac{1}{\sigma_{2}(T)} \sum_{q=0}^{+\infty}
\sum_{
\substack{
{\bk} \in \naturels^3 \\ \abs{{\bk}}=q
}
} a_{f, \sigma}({\bk}, u)\, \int_{T} \widetilde{H}_{{\bk}}(U(t))\, \ud t
]^2\\
&=\mathop {\underline{\lim}}\limits_{\sigma \to 0} \mathbb{E}[J_{f}(u, \sigma)]^2 = \mathbb{E}[J_{f}(u)]^2 < +\infty,
\end{align*}
the two last equalities providing from equality (\ref{Jfnusigma}) and from Proposition \ref{convergence sigma} and the third last equality from the expansion in $L^2(\Omega)$ of $J_{f}(u, \sigma)$ (see (\ref{Jfnusigma})) and precedent remark. The finiteness of $\mathbb{E}[J_{f}(u)]^2$ comes from Lemma \ref{continuite Y}.\\
Thus the \rv
$\eta(T) \in L^2(\Omega)$ and 
\begin{align}
\label{esperance eta}
\mathbb{E}[\eta(T)]^2= \sum_{q=0}^{+\infty} \mathbb{E}\bigg[{\frac{1}{\sigma_{2}(T)} \sum_{
\substack{
{\bk} \in \naturels^3 \\ \abs{{\bk}}=q
}
} a_{f}({\bk}, u)\, \int_{T} \widetilde{H}_{{\bk}}(U(t))\, \ud t}\bigg]^2.
\end{align}
This achieves proof of Lemma \ref{eta T}.
\end{proofarg}
It remains to prove that $J_{f}(u)= \eta(T)$ in $L^2(\Omega)$.\\
As in the proof of Theorem \ref{continuite longueur de courbe}, we write $\normp[L^2(\Omega)]{\cdot}$ for the $L^2(\Omega)$-norm.
\\
For fixed $Q \in \naturels$ and $\sigma >0$, one has the following inequalities
\begin{multline*}
\normp[L^2(\Omega)]{J_{f}(u)-\eta(T)} \le
 \normp[L^2(\Omega)]{\pi_Q(J_{f}(u)-\eta(T))} \\+ \normp[L^2(\Omega)]{\pi^Q(J_{f}(u)-J_{f}(u, \sigma))}
+ \normp[L^2(\Omega)]{\pi^Q(J_{f}(u, \sigma)-\eta(T))}\\
\le  \normp[L^2(\Omega)]{\pi_Q(J_{f}(u))} + \normp[L^2(\Omega)]{\pi_Q(\eta(T))} +\normp[L^2(\Omega)]{J_{f}(u)-J_{f}(u, \sigma)} \\
+\normp[L^2(\Omega)]{\pi^Q(J_{f}(u, \sigma)-\eta(T))}.
\end{multline*}
Now by Lemma \ref{continuite Y} and Lemma \ref{eta T}, $J_{f}(u)$ and $\eta(T)$ belong to $L^2(\Omega)$, thus \\
$\lim_{Q \to +\infty}  \normp[L^2(\Omega)]{\pi_Q(J_{f}(u))} = \lim_{Q \to +\infty} \normp[L^2(\Omega)]{\pi_Q(\eta(T))}=0$.\\
Furthermore due to Proposition \ref{convergence sigma}, $\mathop{\lim}\limits_{\sigma \to 0} \normp[L^2(\Omega)]{J_{f}(u)-J_{f}(u, \sigma)} =0$ and for fixed $Q \in \naturels$ and since $\lim_{\sigma \to 0} a_{f, \sigma}({\bk}, u)=a_{f}({\bk}, u)$, $\mathop{\lim}\limits_{\sigma \to 0} \normp[L^2(\Omega)]{\pi^Q(J_{f}(u, \sigma)-\eta(T))}=0$.\\
Hence, for fixed $Q \in \naturels$ and taking limit as $\sigma$ tends to zero one obtains,
\begin{align*}
\normp[L^2(\Omega)]{J_{f}(u)-\eta(T)} \le
  \normp[L^2(\Omega)]{\pi_Q(J_{f}(u))} + \normp[L^2(\Omega)]{\pi_Q(\eta(T))}.
\end{align*}
Then, taking limit as $Q$ tends to infinity one finally gets
\begin{align*}
\normp[L^2(\Omega)]{J_{f}(u)-\eta(T)}=0.
\end{align*}
Theorem \ref{equality L2 Bis} ensues.
\end{proofarg}
\subsection{One order Rice formula}
\label{Rice}
As by product of Theorem \ref{continuite longueur de courbe}, we obtain the one order Rice formula for the general functional $J_f(u)$, $u$ being a fixed level in $\reels$. More precisely we have the
 \begin{proposition}
 \label{rice}
 For $f: S^{1} \to \reels$ a continuous and bounded function and $u$ a fixed level in $\reels$, one has
$$\mathbb{E}[J_{f}(u)] = p_{X(0)}(u)\, \mathbb{E}[f\left(\frac{\nabla X(0)}{\normp[2]{\nabla X(0)}}\right) \normp[2]{\nabla X(0)}],$$
where ${p}_{X(0)}(\cdot)$ stands for the density of $X(0)$.
\end{proposition}
\begin{proofarg}{Proof of Proposition \ref{rice}}
We refer to Theorem \ref{one order Rice formula} stated in section \ref{sec:notations} and we check that processes $X$ and $\widetilde{Y}:=f(\nu_{X}) \1_{\textrm{D}^{\textrm{r}}_{\textrm{X}}}$ verify assumptions given in this theorem, that is assumptions ${\bf H_1}$, ${\bf H_3}$, ${\bf H_4}$ and ${\bf H_5}$.\\
Theorem \ref{continuite longueur de courbe} implies that assumptions ${\bf H_1}$ and ${\bf H_4}$ are satisfied, assumption ${\bf H_3}$ being fulfilled.\\
Thus we just check if the hypothesis ${\bf H_5}$ is verified.
In this context, let us compute the following integral:
\begin{align*}
\int_T{p}_{X(t)}(u)\,
\mathbb{E} \left[\widetilde{Y}(t) \normp[2]{\nabla
X(t)}|X(t)=u\right] \ud t=\\
\sigma_2(T)\, p_{X(0)}(u)\, \mathbb{E}[f\left(\tfrac{\nabla X(0)}{\normp[2]{\nabla X(0)}}\right) \normp[2]{\nabla X(0)}] < +\infty.
\end{align*}
Since function $\textbf{u} \longmapsto p_{X(0)}(\textbf{u} )$ is continuous, hypothesis ${\bf H_5}$ is satisfied.
All conditions are met to apply the Rice formula, and one gets,
\begin{align*}
\mathbb{E}[J_{f}(u)] &= \frac{1}{\sigma_2(T)} \int_T{p}_{X(t)}(u)\,
\mathbb{E} \left[f\left(\tfrac{\nabla X(t)}{\normp[2]{\nabla X(t)}}\right) \1_{\textrm{D}^{\textrm{r}}_{\textrm{X}}}(t) \normp[2]{\nabla
X(t)}|X(t)=u\right] \ud t\\
&=p_{X(0)}(u)\, \mathbb{E}[f\left(\tfrac{\nabla X(0)}{\normp[2]{\nabla X(0)}}\right) \normp[2]{\nabla X(0)}].
\end{align*}
The proof is completed.
\end{proofarg}

\section{Convergence of general level functionals}
\label{convergence de fonctionnelles}
In this section  $(T_n)_n$ will be open bounded squares of $\reels^2$, with the following form $T_n:= \introo{-n}{n}^2$ with $n \in \naturels^{\star}$, and $n$ tends to infinity. Also let $u$ be a fixed level in $\reels$.\\
Remember that for a continuous and bounded function $f$, $f: S^{1} \to \reels$, we have defined in section \ref{sec:notations} the following general functional $J_{f}^{(n)}(u)$ of the level $u$ as:
\begin{equation*}
J_{f}^{(n)}(u):= \frac{1}{(2n)^2} \int_{\mathcal{C}_{n}(u)} f(\nu_{X}(t)) \ud\sigma_1(t),
\end{equation*}
where $\nu_{X}(t)= {\frac{\nabla X(t)}{\normp[2]{\nabla X(t)}}}$.\\
Our objective consists now to establish the almost sure convergence of such functionals and also their asymptotic normality, when the observation window $T_n$ tends to $\reels^2$ as $n$ tends to infinity.
To this end we will use the results established in previous section \ref{fonctionnelles vues dans le chaos}, that is the Hermite expansion of the random variables $J_{f}^{(n)}(u)$ and also the one order Rice formula in order to compute their expectation.

\subsection{ Almost sure convergence for $J_{f}^{(n)}(u)$}
\label{almost sure convergence for the functionals}
By applying an ergodic theorem for stationary processes (Adler \cite[\S 6.5)]{MR0611857}), we shall show the following general  almost sure convergence theorem.
\begin{theorem}
\label{convergence sure} 
For $f: S^{1} \to \reels$ a continuous and bounded function,
\begin{equation*}
J_{f}^{(n)}(u) \xrightarrow[n\to+\infty]{a.s.}\mathbb{E}[J_{f}^{(1)}(u)].
\end{equation*}
\end{theorem}
\begin{proofarg}{Proof of Theorem \ref{convergence sure}}
Let $f: S^1 \to \reels$ a bounded and continuous fonction.
As a first step we suppose that function $f$ is positive and that the square $T_n$ has the following shape: $T_n:=\introo{0}{n}\times \introo{0}{n}$.
Lemma \ref{lemme cramer} is proved in Appendix \ref{ann:proof}.
\begin{lemma}
Let $f: S^{1} \to \reels$ be a positive continuous and bounded function. One has
\label{lemme cramer}
\begin{multline*}
\int_{0}^{n-1} \int_{0}^{n-1} \int\limits_{\mathcal{C}_{[t\,, \,t+1[\times [s\, ,\,s +1[}(u)} f(\nu_{X}(x)) \ud \sigma_{1}(x)\, \ud t \ud s \\
\hfill- \int_{0}^{1}\int_{0}^{1} \int_{\mathcal{C}_{\introo{0}{t} \times \introo{0}{s}}(u)} f(\nu_{X}(x)) \ud \sigma_{1}(x)\, \ud t \ud s \\ \le\, 
 \int_{\mathcal{C}_{T_n}(u)} f(\nu_{X}(x)) \ud \sigma_1(x) \le \\
 \hfill \int_{0}^{n+1} \int_{0}^{n+1} \int_{\textrm{C}_{[t-1\,,\,t[\times[s-1\,,\,s[}(u)} f(\nu_{X}(x)) \ud \sigma_{1}(x)\, \ud t \ud s
 \hfill
\end{multline*}
\end{lemma}
Let us note $H(t,s):=  \int_{\mathcal{C}_{\introo{0}{t} \times \introo{0}{s}}(u)} f(\nu_{X}(x)) \ud \sigma_{1}(x)$.

On the one hand, since function $f$ is bounded, we have the convergence that follows
\begin{align*}
\frac{1}{(2n)^2} \int_{0}^{1} \int_{0}^{1}  H(t, s) \ud t \ud s \le \frac{{\bf C}}{(2n)^2} \,\sigma_{1}\left({\mathcal{C}_{\introo{0}{1}^2}(u)}\right)\xrightarrow[n\to+\infty]{a.s.} 0,
\end{align*}
last convergence providing from Remark \ref{finitude courbe de Niveau}.
Indeed, since \\$\mathbb{E}[\sigma_{1}\left({\mathcal{C}_{\introo{0}{1}^2}(u)}\right)]^2 <+\infty$, we deduce that $\sigma_{1}\left({\mathcal{C}_{\introo{0}{1}^2}(u)}\right)$ is almost surely finite.
On the other hand, noting by 
$$
\xi(t, s):= \int\limits_{\mathcal{C}_{[t\,, \,t+1[\times [s\, ,\,s +1[}(u)} f(\nu_{X}(x)) \ud \sigma_{1}(x),
$$ one has
\begin{multline*}
\frac{1}{(2n)^2} \int_{0}^{n-1} \int_{0}^{n-1} \int\limits_{\mathcal{C}_{[t\,, \,t+1[\times [s\, ,\,s +1[}(u)} f(\nu_{X}(x)) \ud \sigma_{1}(x)\, \ud t \ud s=\\
\left({\frac{n-1}{2n} }\right)^2{\frac{1}{(n-1)^2} \int_{0}^{n-1} \int_{0}^{n-1} \xi(t, s) \ud t \ud s}.
\end{multline*}
Since process $X$ is a centred stationary Gaussian process with continuous trajectories such that $r_x(t)$ tends to zero as $\normp[2]{t}$ tends to $+\infty$, we deduce from \cite[Theorem 6.5.4]{MR0611857} that process $X$ is ergodic.
Now since process $X$ is strictly stationary and ergodic, we deduce from \cite[Theorem 6.5.2]{MR0611857}, that process $\xi(t,s)$ is a strictly stationary ergodic process. Note that the set $[t, t+1[ \times [s, s+1[$ is not an open set of $\reels^2$. However Proposition \ref{rice} still remains valid. That is the closure of this rectangle does not play role in the expression of the one order Rice formula.
Thus from Proposition \ref{rice} we know that $\mathbb{E}[\abs{\xi(t,s)}] < +\infty$.
By \cite[Theorem 6.5.1]{MR0611857}, we deduce that
\begin{multline*}
\left({\frac{n-1}{2n} }\right)^2 {\frac{1}{(n-1)^2} \int_{0}^{n-1} \int_{0}^{n-1} \xi(t, s) \ud t \ud s}\xrightarrow[n\to+\infty]{a.s.} \frac{1}{4} \,\mathbb{E}[\xi(0,0)]
=\frac{1}{4}\, \mathbb{E}[J_{f}^{(1)}(u)].
\end{multline*}
In the same way, one obtains
\begin{multline*}
\left({\frac{n+1}{2n} }\right)^2 {\frac{1}{(n+1)^2} \int_{0}^{n+1} \int_{0}^{n+1}  \int_{\textrm{C}_{[t-1\,,\,t[\times[s-1\,,\,s[}(u)} f(\nu_{X}(x)) \ud \sigma_{1}(x)\, \ud t \ud s}\\
\xrightarrow[n\to+\infty]{a.s.} \frac{1}{4}\, \mathbb{E}[J_{f}^{(1)}(u)].
\end{multline*}
Finally, by using Lemma \ref{lemme cramer}, one proved that
\begin{align*}
\frac{1}{(2n)^2} \int\limits_{\mathcal{C}_{\introo{0}{n}^2}(u)} f(\nu_{X}(x)) \ud \sigma_{1}(x) \xrightarrow[n\to+\infty]{a.s.}  \frac{1}{4}\, \mathbb{E}[J_{f}^{(1)}(u)].
\end{align*}
Now, working in a similar way succesively with $T_n:=]-n, 0] \times \introo{0}{n}$ or $T_n:=]-n, 0] \times ]-n,0]$ or still with $T_n=\introo{0}{n} \times ]-n,0]$, one should prove the same convergence result for each square $T_n$, even if $T_n$ is not an open set of $\reels^2$.
Finally, if $f$ is a positive function, by using the linearity of the interest functional one have proved that
\begin{align*}
J_{f}^{(n)}(u) \xrightarrow[n\to+\infty]{a.s.} 4 \times \frac{1}{4}\, \mathbb{E}[J_{f}^{(1)}(u)]=\mathbb{E}[J_{f}^{(1)}(u)].
\end{align*}
To conclude the proof of Theorem  \ref{convergence sure}, we decompose function $f$ into its negative part and into its positive part, that is under the shape, $f=f^+-f^-$, and we apply the previous result to each of the interest functionals
$J_{f^{+}}^{(n)}(u)$ and $J_{f^{-}}^{(n)}(u)$.
\end{proofarg}
\subsection{Convergence in law for $\xi_{f}^{(n)}(u)$}
\label{convergence en loi pour les fonctionnelles}
We establish a CLT for a centred and suitably rescalled general functional $J_{f}^{(n)}(u)$, function $f: S^{1} \to \reels$ being any continuous and bounded function and $u$ a fixed level in $\reels$.
Roughly speaking we decided to give the rate of convergence in Theorem \ref{convergence sure}.
In this aim we define the \rv
$\xi_{f}^{(n)}(u)$ by 
\begin{equation}
\label{xi(f)}
\xi_{f}^{(n)}(u):= 2n \left({J_{f}^{(n)}(u) - \mathbb{E}[J_{f}^{(n)}(u)}]\right).
\end{equation}

First we compute the asymptotic variance of $\xi_{f}^{(n)}(u)$ as $n$ goes to infinity, which is the object of the following proposition.
The proof as well as the following remark very closely follow the ones given in \cite[Proposition 2.1]{MR3572325}.
\subsubsection{Asymptotic variance for $\xi_{f}^{(n)}(u)$}
The functionals $\xi_{f}^{(n)}(u)$ are also orthogonal in $L^2(\Omega)$.
This is a crucial fact for computing its variance.
Using the Arcones inequality (see \cite[Lemma 1, p{.} 2245]{MR1331224}), we deduce the asymptotic variance of $\xi_{f}^{(n)}(u)$ as $T_n$ grows to $\reels^2$, this variance depending on the level $u$ as follows.
\begin{proposition}
\label{variance asymptotique xi}
For $f: S^{1} \to \reels$ a continuous and bounded function,
we have the following convergence,
\begin{align*}
\lim_{n \to +\infty} \Var[\xi_{f}^{(n)}(u)] =\mathit{\Sigma}_{f, f}(u),
\end{align*}
$\mathit{\Sigma}_{f, f}(u)$ being defined by 
\begin{align}
\label{Sigma ff}
\mathit{\Sigma}_{f, f}(u):= \sum_{q=1}^{+\infty} \sum_{
\substack{
{\bf k, m} \in \naturels^3 \\ \abs{{\bk}}= \abs{{\bm}}=q
}
} a_f({\bk}, u)\, a_f({\bm}, u)\, R({\bk}, {\bm}),
\end{align}
where coefficients $a_f({\bk}, u)$ have been defined by equality (\ref{af}),
while $R({\bk}, {\bm})$ is defined as
\begin{align}
\label{R(k, m)}
R({\bk}, {\bm}):= \int_{\reels^2} \mathbb{E}[\widetilde{H}_{{\bk}}(U(0)) \widetilde{H}_{{\bm}}(U(v))] \ud v.
\end{align}
\end{proposition}

\begin{remark}
\label{non degeneree}
If $f$ is a function with constant sign, then $\mathit{\Sigma}_{f, f}(u) >0$.
\end{remark}

\begin{proofarg}{Proof of Proposition \ref{variance asymptotique xi}}
One has the following expansion in $L^2(\Omega)$,
\begin{align*}
\xi_{f}^{(n)}(u) = \frac{1}{\sqrt{\sigma_{2}(T_n)}} \sum_{q=1}^{+\infty} \sum_{
\substack{
{\bk} \in \naturels^3 \\ \abs{{\bk}}=q
}
} a_{f}({\bk}, u)\, \int_{T_n} \widetilde{H}_{{\bk}}(U(t))\, \ud t.
\end{align*}
Indeed by using Theorem \ref{equality L2 Bis} it remains to establish that $\mathbb{E}[J_{f}^{(n)}(u)]=a_f({\bf 0}, u)$,
where ${\bf 0}:=(0, 0, 0) \in \naturels^3$.
Thus it is enough to remark that $a(0, u)= \phi\bigg(\dfrac{u}{\sqrt{r_z(0)}}\bigg)\dfrac{1}{\sqrt{r_z(0)}}=p_{X(0)}(u)$ and that
$a_{f}(0, 0)=
\mathbb{E}[f\left(\dfrac{\nabla X(0)}{\normp[2]{\nabla X(0)}}\right) \normp[2]{\nabla X(0)}]$, since
$$
\normp[2]{\nabla X(0)}=\sqrt{\mu} \normp[2]{P \Lambda . \begin{pmatrix} U_1(0)\\ U_2(0) \end{pmatrix}}= \sqrt{\mu} \normp[2]{ \Lambda . \begin{pmatrix} U_1(0)\\ U_2(0) \end{pmatrix}}.
$$
Proposition \ref{rice} gives the result.\\

Since the \rv
$\xi_{f}^{(n)}(u)$ is a centred one, using equality given in (\ref{esperance eta}) and Mehler's formula (see Lemma \ref{mehler} of Section \ref{chaos section}), we obtain
\begin{align}
\label{variance xi}
\var[]{\xi_{f}^{(n)}(u)}=\mathbb{E}[\xi_{f}^{(n)}(u)]^2
= \sum_{q=1}^{+\infty}\sum_{
\substack{
{\bk}, {\bm} \in \naturels^3 \\ \abs{{\bk}}= \abs{{\bm}}=q
}
} a_{f}({\bk}, u)\, a_{f}({\bm}, u)\, R_n({\bk}, {\bm})
\end{align}
with
\begin{align*}
R_n({\bk}, {\bm}):= \frac{1}{(2n)^2} \int_{\introo{-n}{n}^2}  \int_{\introo{-n}{n}^2} \mathbb{E}[\widetilde{H}_{{\bk}}(U(s)) \widetilde{H}_{{\bm}}(U(t))]\ud s\, \ud t.
\end{align*}
Thus and since $U$ is a stationary process, we have
\begin{align*}
R_n({\bk}, {\bm})=  \int_{\introo{-2n}{2n}^2} \mathbb{E}[\widetilde{H}_{{\bk}}(U(0)) \widetilde{H}_{{\bm}}(U(v))]\, \left(1-\frac{\abs{v_1}}{2n}\right)\left(1-\frac{\abs{v_2}}{2n}\right) \ud v.
\end{align*}
Now by applying Lemma \ref{mehler} to $X:=U(0)$ and $Y:=U(v)$, one has for $\abs{{\bk}}=
\abs{{\bm}}$,
\begin{align*}
\mathbb{E}[\widetilde{H}_{{\bk}}(U(0)) \widetilde{H}_{{\bm}}(U(v))]=\sum_{
\substack{
d_{ij} \ge 0 \\ \sum_{j}d_{ij}=k_i \\ \sum_{i}d_{ij}=m_j
}
} {\bk}!  {\bm}! \prod_{1 \le i, j \le 3} \frac{(\mathit{\Gamma}^U_{ij}(v))^{d_{ij}}}{d_{ij}!},
\end{align*}
where 
\begin{align}
\label{Gamma Y v}
\mathit{\Gamma}^U_{ij}(v):= \mathbb{E}[U_i(0)U_j(v)].
\end{align}
Since 
\begin{align*}
\mathit{\Gamma}^U(v)&=(\mathit{\Gamma}^U_{ij}(v))_{1\le i, j \le 3}\\
&=
\begin{pmatrix}
-\ds \frac{1}{\mu} \,P^t (\dfrac{\partial ^{2}r_{z}}{\partial v_{i}\partial v_j}(A.v))_{1\le i, j \le 2}\, P & 
-\ds \frac{1}{\sqrt{\mu r_z(0)}} \,P^t (\dfrac{\partial r_{z}}{\partial v_{i}}(A.v))_{1\le i \le 2}\\
-\ds \frac{1}{\sqrt{\mu r_z(0)}} (\dfrac{\partial r_{z}}{\partial v_{i}}(A.v))_{1\le i \le 2}^t\, P & \ds \frac{r_z(A.v)}{r_z(0)}
\end{pmatrix}
\mbox{,}
\end{align*}
we have for any $v \in \reels^2$,
\begin{align}
\label{majoration psi}
\sup_{1\le i, j \le 3} \abs{\mathit{\Gamma}_{ij}^U(v)} \le {\bf L}\, \Psi(A.v),
\end{align}
where the function $\Psi$ has been introduced in Section \ref{sec:hypotheses} and ${\bf L}$ is some positive constant.
Hence, for $\abs{{\bk}}=\abs{{\bm}}=q$, with $q \in \naturels^{\star}$,
\begin{align*}
\abs{\mathbb{E}[\widetilde{H}_{{\bk}}(U(0)) \widetilde{H}_{{\bm}}(U(v))]} \le {\bf L'}\, \Psi^q(A.v),
\end{align*}
where ${\bf L'}$ is some constant depending on $q$.  \\
By the covariance assumption ${\bf H_{\Psi}}$ previously stated in Section \ref{sec:hypotheses}, $\Psi \in L^1(\reels^2)$ and then $\Psi^q(A\cdot) \in L^1(\reels^2)$.
We can apply the Lebesgue convergence theorem and obtain, for ${\bk}, {\bm} \in (\naturels^{3})^{\star}$,
\begin{align*}
\lim_{n \to +\infty}R_n({\bk}, {\bm}) = R({\bk}, {\bm}):= \int_{\reels^2} \mathbb{E}[\widetilde{H}_{{\bk}}(U(0)) \widetilde{H}_{{\bm}}(U(v))] \ud v.
\end{align*}
Now turning back to (\ref{variance xi}), we write
\begin{align*}
\var[]{\xi_f^{(n)}(u)}= \sum_{q=1}^{+\infty} V_q^{(n)}(u),
\end{align*}
and according to what we have just seen, for all $q \in \naturels^{\star}$,
\begin{align}
\label{def Vq}
V_q(u):=\lim_{n \to +\infty} V_q^{(n)}(u)=\sum_{
\substack{
{\bk}, {\bm} \in \naturels^3 \\ \abs{{\bk}}= \abs{{\bm}}=q
}
} a_{f}({\bk}, u)\, a_{f}({\bm}, u)\, R({\bk}, {\bm}).
\end{align}
Note that for any $q$, $V_q^{(n)}(u) \ge 0$ and so $V_q(u)$.

Thus, if we prove that $\lim_{Q \to +\infty} \sup_{n} \sum_{q=Q+1}^{+\infty} V_q^{(n)}(u)=0$, Fatou's lemma implies that $\lim_{Q \to +\infty} \sum_{q=Q+1}^{+\infty} V_q(u)=0$.
Thus $(\sum_{q=1}^{Q}V_q(u))_Q$ is an upper bounded increasing sequence and consequently a converging sequence, that is the series $\mathit{\Sigma}_{f, f}(u):=\sum_{q=1}^{+\infty} V_q(u)$ will be convergent.
Also the first convergence will imply that $\var[]{\xi_f^{(n)}(u)}$ tends to $\mathit{\Sigma}_{f, f}(u)$.
The proof of Proposition  \ref{variance asymptotique xi} will be completed.
Thus let us prove this convergence via a lemma.
\begin{lemma}
\label{sup variance}
For $q, n \in \naturels^{\star}$, let 
\begin{align*}
V_q^{(n)}(u):= \sum_{
\substack{
{\bk}, {\bm} \in \naturels^3 \\ \abs{{\bk}}= \abs{{\bm}}=q
}
} a_{f}({\bk}, u)\, a_{f}({\bm}, u)\, R_n({\bk}, {\bm}),
\end{align*} 
with $f: S^{1} \to \reels$ a continuous and bounded function. One has the following convergence
\begin{align}
\label{convergence uniforme}
\lim_{Q \to +\infty} \sup_{n \ge 1} \sum_{q=Q+1}^{+\infty} V_q^{(n)}(u)=0.
\end{align}
\end{lemma}
\begin{proofarg}{Proof of Lemma \ref{sup variance}}
First, let us remark that the convergence in (\ref{convergence uniforme}) is equivalent to the following one:
\begin{align*}
\lim_{Q \to +\infty} \var[]{\pi_Q(\xi_{f}^{(n)}(u))}=0,
\end{align*}
uniformly with respect to $n$, where $\pi_Q$ stands for the projection onto the chaos of strictly upper order in $Q$.

For the sake of simplicity of writing, let us note $V_{n, Q}:= \var[]{\pi_Q(\xi_{f}^{(n)}(u))}$.\\
Let $s \in \reels^2$ and $\theta_s$ be the shift operator associated with the field $X$, that is, $\theta_sX:=X_{s+\cdot}$.
Let us also introduce the set of indices $I_n:=[-n, n[^2 \,\cap \,\entiers^2$, clearly we have
\begin{align*}
\pi_Q(\xi_{f}^{(n)}(u))= \frac{1}{2n}\, \sum_{s \in I_n} \theta_s(\pi_Q(\xi_{f,1}(u))),
\end{align*}
where the \rv $\xi_{f,1}(u)$ is 
\begin{align*}
\xi_{f,1}(u):= \sum_{q=1}^{+\infty} \sum_{
\substack{
{\bk} \in \naturels^3 \\ \abs{{\bk}}=q
}
} a_{f}({\bk}, u)\, \int_{\introo{0}{1}^2} \widetilde{H}_{{\bk}}(U(t))\, \ud t.
\end{align*}
The stationarity of $X$ leads to
\begin{align*}
V_{n, Q}= \left({\frac{1}{2n}}\right)^2 \sum_{s \in I_{2n}} \alpha_s(n) \mathbb{E}[\pi_Q(\xi_{f,1}(u))\theta_s(\pi_Q(\xi_{f,1}(u)))],
\end{align*}
where $\alpha_s(n):=\mbox{card}\{t \in I_n,~t-s \in I_n\}$.\\
Obviously, one has $\alpha_s(n) \le (2n)^2$.\\
Now, on the one hand, by the covariance assumption ${\bf H_{\Psi}}$ made in Section \ref{sec:hypotheses},\\ $\lim_{\normp[2]{x} \to +\infty}\Psi(x)=0$, and since the eigenvalues of $A$ are strictly positive one also has\\
 $\lim_{\normp[2]{x} \to +\infty}\Psi(A.x)=0$.\\
On the other hand, let $0< \rho <1$ such that 
\begin{align}
\label{*K}
\rho\, {\bf L}(u )<1 \mbox{\, where \,} {\bf L}(u):= 2\frac{u^2}{r_z(0)}+1,
\end{align}
and $a >0$ such that $\normp[2]{x} \ge a$ implies
\begin{align}
\label{rho}
3\, {\bf L} \Psi(Ax) \le \rho < 1,
\end{align}
where ${\bf L}$ is defined by (\ref{majoration psi}).

We split $V_{n, Q}$ into $V_{n, Q}= V_{n, Q}^{(1)}+V_{n, Q}^{(2)}$, where in $V_{n, Q}^{(1)}$ the sum runs over the indices $s$ in $\{s \in I_{2n}, \normp[\infty]{s} < a+3 \}$ and in $V_{n, Q}^{(2)}$ over $s$ in $\{s \in I_{2n}, \normp[\infty]{s} \ge a+3 \}$, $\normp[\infty]{\cdot}$ standing for the supremum norm.
By Schwarz inequality and since $\alpha_s(n) \le (2n)^2$, using the stationarity of $X$ one has the following upper bound,
\begin{align*}
\abs{V_{n, Q}^{(1)}} \le (2(a+3))^2\, \mathbb{E}[\pi_Q(\xi_{f,1}(u))]^2,
\end{align*}
which goes to zero as $Q$ goes to infinity uniformly with respect to $n$, since adapting equality (\ref{esperance eta}) to the present situation one can proved that $\lim_{Q \to +\infty} \mathbb{E}[\pi_Q(\xi_{f,1}(u))]^2=0$.\\
We proved that $\lim_{Q \to +\infty} \sup_{n} V_{n, Q}^{(1)}=0$.

Now, let us consider  the term $V_{n, Q}^{(2)}$.
\begin{align*}
V_{n, Q}^{(2)}:= \left({\frac{1}{2n}}\right)^2 \sum_{\substack{s \in I_{2n} \\ \normp[\infty]{s} \ge a+3}} \alpha_s(n) \mathbb{E}[\pi_Q(\xi_{f,1}(u))\theta_s(\pi_Q(\xi_{f,1}(u)))].
\end{align*}
For $q \in \naturels^{\star}$, let us define function $F_q$ by
\begin{align*}
F_q(x):=\sum_{
\substack{
{\bk} \in \naturels^3 \\ \abs{{\bk}}=q
}
} a_{f}({\bk}, u)\, \widetilde{H}_{{\bk}}(x),\, x \in \reels^3.
\end{align*}
For $s \in I_{2n}$ such that $\normp[\infty]{s} \ge a+3$,
\begin{multline*}
\mathbb{E}[\pi_Q(\xi_{f,1}(u))\theta_s(\pi_Q(\xi_{f,1}(u)))]=\\
 \sum_{q=Q+1}^{+\infty} \int_{\introo{0}{1}^2}  \int_{\introo{0}{1}^2} \mathbb{E}[F_q(U(t))F_q(U(s+v))]\, \ud t \ud v.
\end{multline*}
At this step of the proof we want to propose a bound for $\mathbb{E}[F_q(U(t))F_q(U(s+v))]$, $t, v \in \introo{0}{1}^2$.\\
Note that if $\sum_{
\substack{
{\bk} \in \naturels^3 \\ \abs{{\bk}}=q
}
} a^2_{f}({\bk}, u) {\bk}! = 0$, then $F_q(x)=0$ for all $x \in \reels^3$ and a trivial bound is zero.
So let us suppose that $\sum_{
\substack{
{\bk} \in \naturels^3 \\ \abs{{\bk}}=q
}
} a^2_{f}({\bk}, u) {\bk}! \neq 0$.\\
We are going to apply Arcones inequality (see \cite[Lemma 1 p.
2245]{MR1331224}).
By using notations of this lemma, we apply it to $f:=F_q$ and to $X=(X^{(j)})_{1\le j \le 3}:= U(t)$ and $Y=(Y^{(k)})_{1\le k \le 3}:= U(s+v)$, with $d=3$, such that $r^{(j, k)}= \mathbb{E}[X^{(j)}Y^{(k)}]:= \mathit{\Gamma}_{jk}^U(s-t+v)$, $\mathit{\Gamma}^U$ being defined in (\ref{Gamma Y v}).\\
Now by using inequalities given in (\ref{majoration psi}) and (\ref{rho}),
\begin{align}
\label{arcones inequality}
\psi:= \left({\sup_{1\le j \le 3} \sum_{k=1}^3 \abs{r^{(j, k)}} }\right) \vee \left({\sup_{1\le k \le 3} \sum_{j=1}^3 \abs{r^{(j, k)}} }\right) \nonumber \\
\le 3\, {\bf L}\, \Psi(A.(s-t+v)) \le \rho <1.
\end{align}
It remains to verify that $F_q$ function on $\reels^3$ has finite second moment and rank $q$.\\
In the first place by Lemma \ref{mehler} given in Section \ref{chaos section} one has
\begin{align*}
\mathbb{E}[F_q(X)]^2=\mathbb{E}[F_q(U(t))]^2= \sum_{
\substack{
{\bk} \in \naturels^3 \\ \abs{{\bk}}=q
}
} a^2_{f}({\bk}, u) {\bk}! < +\infty.
\end{align*}
In the second place and since $\sum_{
\substack{
{\bk} \in \naturels^3 \\ \abs{{\bk}}=q
}
} a^2_{f}({\bk}, u) {\bk}! \neq 0$, this last equality ensures that rank $F_q \le q$. Furthermore let ${\bm} \in \naturels^3$ such that $\mathbb{E}[F_q(X)\widetilde{H}_{{\bm}}(X)] \neq 0$.
By Lemma \ref{mehler}, $\mathbb{E}[F_q(X)\widetilde{H}_{{\bm}}(X)]= \sum_{
\substack{
{\bk} \in \naturels^3 \\ \abs{{\bk}}=q
}
} a_{f}({\bk}, u) {\bk}! \1_{\abs{{\bk}}=\abs{{\bm}}}$, which implies ${\abs{\bm}}=q$ and rank $F_q= q$.\\\\
Thus we have all the ingredients to apply Arcones inequality.
For $q \ge 1$, using inequality given in (\ref{arcones inequality}) we get the bound
\begin{multline*}
\mathbb{E}[F_q(U(t))F_q(U(s+v))] \le \psi^q\, \mathbb{E}[F_q(U(t))]^2\\
\le \rho^{q-1}\, (3{\bf L})\, \Psi(A.(s-t+v))\, (\sum_{
\substack{
{\bk} \in \naturels^3 \\ \abs{{\bk}}=q
}
} a^2_{f}({\bk}, u) {\bk}!).
\end{multline*}
As already pointed out in Remark \ref{serie dirac}, the series $\sum_{{\bk} \in \naturels^3} a^2_{f}({\bk}, u) {\bk}! =+\infty$, so that we have to tread carefully in what follows.\\
Finally and since $\alpha_s(n) \le (2n)^2$, one has
\begin{align*}
\abs{V_{n, Q}^{(2)}} \le {\bf C} \sum_{q=Q+1}^{+\infty} \sum_{
\substack{
{\bk} \in \naturels^3 \\ \abs{{\bk}}=q
}
} \rho^{q-1}\, a^2_{f}({\bk}, u) {\bk}!\, \sum_{s \in I_{2n}} \int_{\introo{0}{1}^2} \int_{\introo{0}{1}^2}\, \Psi(A.(s-t+v))\, \ud t\, \ud v.
\end{align*}
Using that 
\begin{multline*}
\sum_{s \in I_{2n}} \int_{\introo{0}{1}^2} \int_{\introo{0}{1}^2}\, \Psi(A.(s-t+v))\, \ud t\, \ud v \le \int_{\reels^2} \Psi(A.v)\, \ud v\\ \le {\bf C}  \int_{\reels^2} \Psi(v)\, \ud v < +\infty,
\end{multline*}
last finiteness providing from assumption ${\bf H_{\Psi}}$ made in Section \ref{sec:hypotheses}, one has
\begin{align*}
\abs{V_{n, Q}^{(2)}} \le {\bf C} \sum_{q=Q+1}^{+\infty} \sum_{
\substack{
{\bk} \in \naturels^3 \\ \abs{{\bk}}=q
}
} \rho^{q-1}\, a^2_{f}({\bk}, u) {\bk}!.
\end{align*}
To conclude the proof of this lemma we just have to check that
\begin{align*}
\sum_{q=1}^{+\infty} \sum_{
\substack{
{\bk} \in \naturels^3 \\ \abs{{\bk}}=q
}
} \rho^{q-1}\, a^2_{f}({\bk}, u) {\bk}! < +\infty.
\end{align*}
Now remember that for ${\bk}=(k_i)_{1 \le i \le 3} \in \naturels^3$, 
\begin{align*}
a_f({\bk}, u)= a_f(k_1, k_2) \,a(k_3, u),
\end{align*}
with $a_f(k_1, k_2)$ and $a(k_3, u)$ respectively defined by equalities (\ref{a f k}) and (\ref{ak3}).

On the one hand, since the function
$$
h: (y_1, y_2) \mapsto f\Big(\frac
	{P\Lambda .{\footnotesize \begin{pmatrix} y_1\\y_2 \end{pmatrix}}}
	{\sqrt{\lambda_1^2y_1^2+\lambda_2^2y_2^2}}
	\Big)
	\sqrt{\lambda_1^2y_1^2+\lambda_2^2y_2^2}
$$
is such that $h \in L(\reels^2, \phi_2(y) \ud y)$, we deduce that,
\begin{align}
\label{finitude k1 k2}
\sum_{k_1, k_2 \in \naturels} a_f^2(k_1, k_2) k_1!\, k_2! < +\infty.
\end{align}
On the other hand by using the expression of Hermite's polynomials given in Appendix \ref{ann:proof} by (\ref{hermite form}) and (\ref{H2l+1}), for all $k \in \naturels$ and for all $x \in \reels$ we get the bound
\begin{align*}
H_k^2(x) \le (k+1)!\, (2x^2+1)^k,
\end{align*}
so that for all ${\bk} \in \naturels^3$ such that $\abs{{\bk}}=q$ and remembering that ${\bf L}(u)$ has been defined in (\ref{*K}), one has
\begin{align*}
k_3!\, a^2(k_3, u) \le {\bf C}\, (k_3+1)\, {\bf L}^{k_3}(u)
\le {\bf C}\, (q+1)\, {\bf L}^q(u).
\end{align*}
Thus by inequality (\ref{finitude k1 k2}), finally one obtains
\begin{align*}
\sum_{q=1}^{+\infty} \sum_{
\substack{
{\bk} \in \naturels^3 \\ \abs{{\bk}}=q
}
} \rho^{q-1}\, a^2_{f}({\bk}, u) {\bk}! \le {\bf C}\, \sum_{q=1}^{+\infty} \rho^{q-1}\, (q+1)^2\, {\bf L}^q(u)< +\infty,
\end{align*}
last finiteness providing from inequality (\ref{*K}).\\
This yields Lemma \ref{sup variance}.
\end{proofarg}
Proposition \ref{variance asymptotique xi} ensues.
\end{proofarg}
\begin{proofarg}{Proof of Remark \ref{non degeneree}}
Remark ensues from the following argumentation.\\
We have seen in the proof of Proposition \ref{variance asymptotique xi} that 
$\mathit{\Sigma}_{f, f}(u)= \sum_{q=1}^{+\infty} V_q(u)$, with 
$$V_q(u) = \sum_{
\substack{
{\bf k, m} \in \naturels^3 \\ \abs{{\bk}}= \abs{{\bm}}=q
}
} a_f({\bk}, u)\, a_f({\bm}, u)\, R({\bk}, {\bm})\ge 0,$$ for all $q \ge 1$.
Thus 
\begin{align*}
\mathit{\Sigma}_{f, f}(u) \ge V_1(u)+ V_2(u).
\end{align*}
By using Lemma \ref{mehler} and the inversion formula, a computation gives that for $\abs{{\bk}}= \abs{{\bm}}=1$, $R({\bk}, {\bm})=0$ except when $\bk=\bm= (0, 0, 1)$ and in this case one has \begin{equation*}
R((0, 0, 1), (0, 0, 1))= \frac{1}{r_x(0)} \int_{\reels^2} r_x(v)\, \ud v= (2\pi)^2\, \frac{f_z(0)}{\lambda_1 \lambda_2 r_z(0)}.
\end{equation*}
Thus 
\begin{align*}
V_1(u)=a^2_f(0, 0)\, \frac{u^2}{r_z^2(0)}\, \phi^2\Big(\frac{u}{\sqrt{r_z(0)}}\Big) (2\pi)^2\, \frac{f_z(0)}{\lambda_1 \lambda_2 r_z(0)} >0,
\end{align*}
if $u \neq 0$, since $f$ is supposed to have constant sign and $f_z(0) >0$ (see Remark \ref{spectral positivity} given in Section \ref{sec:hypotheses}).\\
Using arguments similar to the previous ones, the fact that $\int_{\reels^2} f_z(t)  \normp[2]{t}^2 \ud t < +\infty$ (see Remark \ref{norme density}) and Parseval equality, straightforward calculations provide that
\begin{multline*}
\lefteqn{V_2(u) = 2 \times (2\pi)^2 \int_{\reels^2} f_x^2(t) \times \bigg[{ a_f((1, 1, 0), u) (d_{11}d_{21}t_1^2 }}\\ 
\left.{+ (d_{11}d_{22}+d_{12}d_{21})t_1t_2
+d_{12}d_{22}t_2^2) + a_f((1, 0, 1), u) \frac{1}{\sqrt{r_z(0)}} (d_{11}t_1+d_{12}t_2) }\right.\\
\left.{+ a_f((0, 1, 1), u) \frac{1}{\sqrt{r_z(0)}} (d_{21}t_1+d_{22}t_2)+
a_f((2, 0, 0), u) (d_{11}t_1+d_{12}t_2)^2 }\right.\\
{ +a_f((0, 2, 0), u) (d_{21}t_1+d_{22}t_2)^2+ a_f((0, 0, 2), u) \frac{1}{r_z(0)}
  }\bigg]^2 \ud t \ge 0,
\end{multline*}
where $(d_{ij})_{1\le i, j \le 2}=D:=\ds \frac{1}{\sqrt{\mu}} \,\Lambda^{-1}P^t$.
\begin{remark}
Note that in the case where the process $ X $ is isotropic our result contains that of Kratz and Le\'on \cite[Theorem 3]{MR1860517}.
\end{remark}
Since $\det(D) \neq 0$, one gets the following equivalence:
$$
(V_2(u)=0) \iff (a_f({\bk}, u)=0, \mbox{\, for all\, } {\bk} \in \naturels^{3} \mbox{\, such that  } \abs{\bf{k}}=2)
$$
In particular, since $f$ has a constant sign, $a_f(0,0) \neq 0$ so that $a_f((0, 0, 2), 0) \neq 0$ and $V_2(0) >0$.\\
Finally we proved that for $u \neq 0$, $\mathit{\Sigma}_{f, f}(u) \ge V_1(u)+ V_2(u) \ge V_1(u) >0$ and for $u=0$, $\mathit{\Sigma}_{f, f}(0) \ge V_1(0)+ V_2(0) \ge V_2(0) >0$.\\
The proof of Remark \ref{non degeneree} is completed.
\end{proofarg}

Now, we have got all the tools to prove that the \rv $\xi_f^{(n)}(u)$ converges in law as $n$ tends to infinity to a centred Gaussian variable with finite variance $\mathit{\Sigma}_{f, f}(u)$ given by (\ref{Sigma ff}), see Theorem \ref{Peccati}.
The proof we give for this theorem, is inspired by the one presented in \cite[Proposition 2.4]{MR3572325}.
\subsubsection{General level functionals viewed into the Wiener-It\^o chaos}
\label{section xi}
Using the Peccati and Tudor theorem (see \cite{MR2126978}), we obtain the following theorem.
\begin{theorem}
\label{Peccati}
For $f: S^{1} \to \reels$ a continuous and bounded function,
we have the following convergence,
\begin{align*}
\xi_{f}^{(n)}(u) \xrightarrow[n \to +\infty]{Law} \mathcal{N}(0; \mathit{\Sigma}_{f, f}(u)).
\end{align*}
\end{theorem}
\begin{remark}
For example if $f \equiv 1$, we find that the normalized centred curve length converges in law to a non degenerate Gaussian \rv.\\
\end{remark}
Note that for all real numbers $a$ and $b$ and for all continuous and bounded functions, $f_1$ and $f_2: S^{1} \to \reels$, one has
$\xi_{af_1+bf_2}^{(n)}(u)= a \xi_{f_1}^{(n)}(u)+b \xi_{f_2}^{(n)}(u)$.\\
By generalizing the definition given in (\ref{xi(f)}), we define the following functional. For $k \in \naturels^{\star}$ and $\overrightarrow{f}:=(f_1, f_2, \cdots, f_k): S^{1} \to \reels^k$ continuous and bounded function let
\begin{align*}
\xi_{\overrightarrow{f}}^{(n)}(u):= 2n \left({J_{\overrightarrow{f}}^{(n)}(u) -  \mathbb{E}[J_{\overrightarrow{f}}^{(n)}(u)]}\right)
\end{align*}
where
\begin{align}
\label{J(f) Bis}
J_{\overrightarrow{f}}^{(n)}(u) -  \mathbb{E}[J_{\overrightarrow{f}}^{(n)}(u)]
:= \left({ J_{f_i}^{(n)}(u) - \mathbb{E}[J_{f_i}^{(n)}(u)]}\right)_{1 \le i \le k}.
\end{align}
We define
\begin{align*}
\mathit{\Sigma}_{f_i, f_j}(u):= \sum_{q=1}^{+\infty} \sum_{
\substack{
{\bk,\,\bm} \in \naturels^3 \\ \abs{{\bk}}= \abs{{\bm}}=q
}
} a_{f_i}({\bk}, u)\, a_{f_j}({\bm}, u)\, R({\bk}, {\bm}),
\end{align*}
where $R({\bk}, {\bm})$ is defined by (\ref{R(k, m)}), 
and 
\begin{align}
\label{Sigma vectoriel}
\mathit{\Sigma}_{\overrightarrow{f}}(u):= \left({\mathit{\Sigma}_{f_i, f_j}(u)}\right)_{1 \le i, j \le k}.
\end{align}
By Cramer-Wold device we readily get Corollary \ref{convergence triple}.\\
\begin{corollary}
\label{convergence triple}
For all $k \in \naturels^{\star}$ and all $\overrightarrow{f}:=(f_1, \cdots, f_k): S^{1} \to \reels^k$ continuous and bounded function, one has
\begin{align*}
\xi_{\overrightarrow{f}}^{(n)}(u) \xrightarrow[n \to +\infty]{Law} \mathcal{N}(0; \mathit{\Sigma_{\overrightarrow{f}}}(u)).
\end{align*}
\end{corollary}
\begin{proofarg}{Proof of Theorem \ref{Peccati}}
First, let $Q$ a fixed integer in $\naturels^{\star}$ and let us consider the projection of the \rv $\xi_f^{(n)}(u)$ onto the first $Q$ chaos in $L^2(\Omega)$ that is
\begin{align*}
\pi^Q(\xi_f^{(n)}(u)):= \frac{1}{2n} \sum_{q=1}^{Q} \sum_{
\substack{
{\bk} \in \naturels^3 \\ \abs{{\bk}}=q
}
} a_{f}({\bk}, u)\, \int_{T_n} \widetilde{H}_{{\bk}}(U(t))\, \ud t.
\end{align*}
We will show the asymptotic normality of this sequence as $n$ tends to infinity.
For this purpose and in order to apply the Peccati and Tudor theorem (see \cite[Theorem 1]{MR2126978}), we will give an expansion of this \rv
into the Wiener-It\^o chaos of order less or equal to $Q$.\\
To this end, remember that in Section \ref{chaos section} for any $t \in \reels^2$ one has defined  the $3$-dimensional standard Gaussian vector $U(t)$ as
\begin{align*}
U(t)=\int_{\reels^2} e^{i \,\langle t ,\, \lambda \rangle} \sqrt{f_x(\lambda)}\, \Delta^{-1}.\nu(\lambda)\ud W(\lambda)= (U_{i}(t))_{1 \le i \le 3},
\end{align*}
where $W$ stands for the standard Brownian motion and for any $\lambda=(\lambda_i)_{1 \le i \le 2}$ in $\reels^2$, 
\begin{align*}
\nu(\lambda)=(i\lambda_1, i \lambda_2, 1). 
\end{align*}
In what follows, for any $t \in\reels^2$ and $j=1, 2, 3$, we denote by $\varphi_{t, j}$ the square integrable map on $\reels^2$ such that,
\begin{align*}
U_j(t)=\int_{\reels^2} \varphi_{t, j}(\lambda)\ud W(\lambda).
\end{align*}
Since $(\varphi_{t, j})_{1 \le j \le 3}$ is an orthonormal system in $L^2(\reels^2)$, using It\^o's formula (see Major \cite[Theorem 4.2 p.
30]{MR0611334}), for fixed ${\bk}=(k_i)_{1 \le i \le 3} \in \naturels^3$ such that $\abs{{\bk}}=q$,
\begin{align*}
\widetilde{H}_{{\bk}}(U(t))
	&= \prod_{j=1}^3 H_{k_j}(U_j(t))\\
	&= \int_{\reels^{2q}} \varphi_{t, 1}^{\otimes k_1} \otimes \varphi_{t, 2}^{\otimes k_2} \otimes \varphi_{t, 3}^{\otimes k_3}	(\lambda_1, \dots, \lambda_q)\ud W(\lambda_1) \cdots \ud W(\lambda_q)\\
	&=I_q(\varphi_{t, 1}^{\otimes k_1} \otimes \varphi_{t, 2}^{\otimes k_2} \otimes \varphi_{t, 3}^{\otimes k_3}),
\end{align*}
where $I_q$ stands for the Wiener-It\^o integral of order $q$.\\
We shall use notations introduced in Slud \cite{MR1303648}.\\
For each $q \ge 1$, let consider $L^2_{\sym}((\reels^2)^q)$ the complex Hilbert-space
\begin{align*}
L^2_{\sym}((\reels^2)^q):=\{f_q \in L^2((\reels^2)^q), \mbox{\, such that for all \,} x \in (\reels^2)^q, \\
f_q(x)=\overline{f_q(-x)}, f_q(x_1, \dots, x_q)=f_q(x_{\pi(1)}, \dots, x_{\pi(q)}), \forall \pi \in S_q \}, 
\end{align*}
where $S_q$ denotes the symmetric group of permutations of $\{1, \dots, q\}$.\\
For $f_q \in L^2((\reels^2)^q)$, we denote by $\sym(f_q)$ the symmetrization of $f_q$, that is for $x_1, \dots, x_q \in \reels^2$,
\begin{align*}
\sym(f_q)(x_1, \dots, x_q):=\frac{1}{q!}\sum_{\pi \in S_q} f_q(x_{\pi(1)}, \dots, x_{\pi(q)}).
\end{align*}
Observe that for $f_q \in L^2((\reels^2)^q)$ such that for all $x \in (\reels^2)^q,
f_q(x)=\overline{f_q(-x)}$, one has 
\begin{align}
\label{Iq}
I_q(f_q)=I_q(\sym(f_q)).
\end{align}
For $q \ge 1$, $f_q \in L^2_{\sym}((\reels^2)^q)$ and $p=1, \dots, q$, we will write $f_q \otimes_p f_q$ for the $p$-th contraction of $f_q$ defined as
\begin{align*}
f_q \otimes_p f_q(\lambda_1, \dots, \lambda_{2q-2p}):=\int_{(\reels^2)^p} f_q(\lambda_1, \dots,\lambda_{q-p}, x_1, \dots, x_p)\\
f_q(\lambda_{q-p+1}, \dots,\lambda_{2q-2p}, -x_1, \dots, -x_p)\, \ud x_1 \cdots \ud x_p.
\end{align*}
Using the property of $I_q$ given in (\ref{Iq}), the \rv
of interest can be written as
\begin{align*}
\pi^Q(\xi_f^{(n)}(u))=\sum_{q=1}^Q I_q(f_q^{(n)}),
\end{align*}
where for $n, q \in \naturels^{\star}$
\begin{align*}
f_q^{(n)}:= \frac{1}{2n} \sum_{
\substack{
{\bk} \in \naturels^3 \\ \abs{{\bk}}=q
}
} a_{f}({\bk}, u)\, \int_{T_n} \sym(\varphi_{t, 1}^{\otimes k_1} \otimes \varphi_{t, 2}^{\otimes k_2} \otimes \varphi_{t, 3}^{\otimes k_3})\, \ud t.
\end{align*}
We  symmetrized the function $f_q^{(n)}$ with the aim of applying \cite[Theorem 1]{MR2126978}.\\
Symmetrizing the function complicates a lot the calculations in the study of the contractions.
So we are going to write function $f_q^{(n)}$ in another way.\\
For ${\bk} =(k_i)_{1 \le i  \le 3}$ such that $\abs{{\bk}}=q$, we define
\begin{align*}
\mathcal{A}_{{\bk}}:= \{m=(m_1, \dots, m_q) \in \{1, 2, 3\}^q, 
\forall i=1, 2, 3,\, \sum_{j=1}^q \1_{\{ i\}}(m_j)=k_i \},
\end{align*} 
one\ has $\card\{\mathcal{A}_{{\bk}}\}={{\bk}!}/{q!}$.
Let us remark that the sets $(\mathcal{A}_{{\bk}})_{{\bk} \in \naturels^3, \abs{{\bk}}=q}$ provide a partition of $\{1, 2, 3\}^q$.\\
With these notations one has
\begin{multline*}
\sum_{
\substack{
{\bk} \in \naturels^3 \\ \abs{{\bk}}=q
}
} a_{f}({\bk}, u)\,  \sym(\varphi_{t, 1}^{\otimes k_1} \otimes \varphi_{t, 2}^{\otimes k_2} \otimes \varphi_{t, 3}^{\otimes k_3})=\\
 \sym\Big(\sum_{
\substack{
{\bk} \in \naturels^3 \\ \abs{{\bk}}=q
}
} \sum_{m \in \mathcal{A}_{{\bk}}} \,\frac{a_{f}({\bk}, u)}{\card(\mathcal{A}_{{\bk}})}\, \varphi_{t, m_1} \dots \varphi_{t, m_q}\Big).
\end{multline*}
Since $(\mathcal{A}_{{\bk}})_{{\bk} \in \naturels^3, \abs{{\bk}}=q}$ provides a partition of $\{1, 2, 3\}^q$, then for all $m \in \{1, 2, 3\}^q, \exists ! {\bk} \in \naturels^3$ such that $\abs{{\bk}}=q$ and $m \in \mathcal{A}_{{\bk}}$.
So for fixed $m \in \{1, 2, 3\}^q$, let 
$
b_{f}(m, u):= \ds \frac{a_{f}({\bk}, u)}{\card(\mathcal{A}_{{\bk}})}
$.
Thus
\begin{align*}
\lefteqn{\sum_{
\substack{
{\bk} \in \naturels^3 \\ \abs{{\bk}}=q
}
} a_{f}({\bk}, u)\,  \sym(\varphi_{t, 1}^{\otimes k_1} \otimes \varphi_{t, 2}^{\otimes k_2} \otimes \varphi_{t, 3}^{\otimes k_3})}\\
&\qquad=
 \sym\Big(\sum_{m \in \{1, 2, 3\}^q} b_{f}(m, u) \varphi_{t, m_1} \cdots \varphi_{t, m_q}\Big)\\
&\qquad=\sum_{m \in \{1, 2, 3\}^q} b_{f}(m, u) \varphi_{t, m_1} \cdots \varphi_{t, m_q}
\end{align*}
since $m \mapsto b_{f}(m, u)$ is symmetric on $\{1, 2, 3\}^q$.\\
Finally the \rv
$\pi^Q(\xi_f^{(n)}(u))$ can be written as 
\begin{align*}
\pi^Q(\xi_f^{(n)}(u))= \sum_{q=1}^Q I_q(f_q^{(n)}),
\end{align*}
where for $n, q \in \naturels^{\star}$
\begin{align*}
f_q^{(n)}= \frac{1}{2n}\, \int_{]-n, n[^2}\,
 \sum_{m \in \{1, 2, 3\}^q} b_{f}(m, u) \varphi_{t, m_1} \cdots \varphi_{t, m_q}\, \ud t,
\end{align*}
that ends our first objective.\\
To obtain convergence of $\pi^Q(\xi_f^{(n)}(u))$, we use \cite[Theorem 1]{MR2126978}.
Convergence in Proposition \ref{variance asymptotique xi} gives the required conditions appearing at the beginning of this latter theorem.
So we just verify condition (i) in proving the following lemma.
\begin{lemma}
\label{contractions}
For fixed integers $q$ and $p$ such that $q \ge 2$ and $p=1, \dots, q-1$,
\begin{align*}
\lim_{n \to +\infty} \int_{(\reels^2)^{2(q-p)}} \abs{f_q^{(n)} \otimes_p f_q^{(n)}(\lambda_1, \dots, \lambda_{q-p}, \mu_1, \dots, \mu_{q-p})}^2\, \\
\ud\lambda_1 \dots \ud\lambda_{q-p} \ud\mu_{1}\dots\ud\mu_{q-p}=0.
\end{align*}
\end{lemma}
\begin{proofarg}{Proof of Lemma \ref{contractions}}
Let
\begin{align*}
C_n:=\int_{(\reels^2)^{2(q-p)}} \abs{f_q^{(n)} \otimes_p f_q^{(n)}(\lambda_1, \dots, \lambda_{q-p}, \mu_1, \dots, \mu_{q-p})}^2\, \\
\ud\lambda_1 \dots \ud\lambda_{q-p} \ud\mu_{1}\dots\ud\mu_{q-p}.
\end{align*}
Straightforwards calculations show that
\begin{multline*}
C_n= \left(\frac{1}{2n}\right)^4\int_{\left(]-n, n[^2\right)^4} \sum_{k, m \in \{1, 2, 3\}^q} \sum_{K, M \in \{1, 2, 3\}^q} b_{f}(k, u)\, b_{f}(m, u)\, \\
b_{f}(K, u)\, b_{f}(M, u)\, \mathit{\Gamma}^U_{k_1, K_1}(t_1-t_2)\dots \mathit{\Gamma}^U_{k_{q-p}, K_{q-p}}(t_1-t_2)\, \times\\
 \mathit{\Gamma}^U_{m_1, M_1}(s_1-s_2)\cdots \mathit{\Gamma}^U_{m_{q-p}, M_{q-p}}(s_1-s_2)\, 
 \mathit{\Gamma}^U_{k_{q-p+1}, m_{q-p+1}}(t_1-s_1) \cdots   \times\\
 \mathit{\Gamma}^U_{k_{q}, m_{q}}(t_1-s_1)\,
 \mathit{\Gamma}^U_{K_{q-p+1}, M_{q-p+1}}(t_2-s_2) \cdots \mathit{\Gamma}^U_{K_{q}, M_{q}}(t_2-s_2)\,
 \ud t_1\ud s_1\, \ud t_2\ud s_2\,,
 \end{multline*}
 where $\mathit{\Gamma}^U$ has been defined in (\ref{Gamma Y v}).\\
Using inequality (\ref{majoration psi}), we get the bound
 \begin{multline*}
C_n \le  \left(\frac{1}{2n}\right)^4 {\bf L}^{2q} \Big(\sum_{m \in \{1, 2, 3\}^q} \abs{b_{f}(m, u)}\Big)^4
\int_{\left(]-n, n[^2\right)^4} \Psi^{q-p}(A.(t_1-t_2)) \times \\
\Psi^{q-p}(A.(s_1-s_2)) \Psi^p(A.(t_1-s_1)) \Psi^p(A.(t_2-s_2)) \ud t_1\ud s_1\ud t_2\ud s_2.
\end{multline*}
Moreover, we have
\begin{align*}
\Psi^{q-p}(A.(t_1-t_2)) \Psi^{p}(A.(t_2-s_2)) \le \Psi^{q}(A.(t_1-t_2)) +\Psi^{q}(A.(t_2-s_2)).
\end{align*}
Furthermore for $r \in \naturels^{\star}$, one has
\begin{align*}
\int_{]-n,n[^2} \Psi^r(A.(u-v)) \ud u \ud v \le {\bf C} \int_{\reels^2} \Psi^r(u) \ud u < +\infty.
\end{align*}
Applying the penultimate and last inequalities to $p \ge 1$, $q\ge 1$ and $q-p \ge 1$, one obtains
\begin{align*}
C_n \le  {\bf C_q}\, \left(\frac{1}{2n}\right)^2\, \Big(\int_{\reels^2} \Psi^q(u)\,  du\Big)\Big(\int_{\reels^2} \Psi^{q-p}(u)\ud u\Big)
\Big(\int_{\reels^2} \Psi^p(u)\ud u\Big),
\end{align*}
thus we proved that $\lim_{n \to +\infty}C_n=0$, this achieves proof of Lemma \ref{contractions}.
\end{proofarg}

Hence we proved that,
\begin{itemize}
\item
for all $Q \ge 1$,
$
\pi^Q(\xi_f^{(n)}(u)) \xrightarrow[n \to +\infty]{Law} N(0; \sum_{q=1}^Q V_q(u)),
$
where $V_q(u)$ has been defined by (\ref{def Vq}).\\
On the other hand we proved in Lemma \ref{sup variance} that
\item for all $n \ge 1$,
$
\pi^Q(\xi_f^{(n)}(u)) \xrightarrow[Q \to +\infty]{L^2(\Omega)} \xi_f^{(n)}(u),
$\\
and that
\item 
$
\lim\limits_{Q \to +\infty} \sum\limits_{q=1}^Q V_q(u)=\mathit{\Sigma}_{f, f}(u).
$\\
Finally and by Proposition \ref{variance asymptotique xi} we also have
\item
$
\lim\limits_{Q\to +\infty} \lim\limits_{n \to +\infty} \normp[L^2(\Omega)]{\pi^Q(\xi_f^{(n)}(u)) - \xi_f^{(n)}(u)}=0.
$
\end{itemize}
Applying Dynkin \cite[Lemma 1.1]{MR0920254}, we can conclude that 
$
\xi_f^{(n)}(u) \xrightarrow[n \to +\infty]{Law} N(0; \mathit{\Sigma}_{f, f}(u)),
$ that achieves proof of Theorem \ref{Peccati}.
\end{proofarg}

\section{The affinity estimators}
\label{Les estimateurs}
In this section  $(T_n)_n$ will be still open bounded squares of $\reels^2$, with the following form $T_n:= \introo{-n}{n}^2$ with $n \in \naturels^{\star}$, and $n$ tends to infinity. Also $u$ will be a fixed level in $\reels$.\\
We propose estimators of the affinity parameters defined as $\ds \lambda$ and $\theta$, by considering the particular functional of the level set $u$ proposed by Wschebor \cite[chap. 3.6 F, pages 82-85]{MR0871689}. In this aim we reproduce in sections \ref{fonctionnelle d'interet}, \ref{le diffeomorphisme} and \ref{les parametres d'affinite} all the computations made by Wschebor including the three forthcoming attached pictures.  By using the general convergence results established in previous section \ref{convergence de fonctionnelles} we point out convergence properties of these estimators.
\subsection{The functional of interest}
\label{fonctionnelle d'interet}
Let $J_{\overrightarrow{f^{\star}}}^{(n)}(u)$ this particular functional (see notations given in (\ref{J(f) Bis}) for the definition of a general vector-functional of the level $u$), where the vector-function $\overrightarrow{f^{\star}}$ is defined as follows.
Let $v^{\star} \in S^{1}$ {\bf a fixed vector} and consider 
\begin{equation}
\label{fstar}
\begin{aligned}S^{1}& \longrightarrow S^{1}\\
\omega & \mapsto \overrightarrow{f^{\star}}(\omega):= \omega \times (\1_{\{\langle\omega, v^{\star}\rangle \ge 0\}}-\1_{\{\langle\omega, v^{\star}\rangle < 0\}}).
\end{aligned}
\end{equation}
\begin{remark}
Note that function $\overrightarrow{f^{\star}}$ is continuous except on the line given by equation $\langle\omega, v^{\star}\rangle = 0$.  Nevertheless Rice formulas will still remain valid for the associated functionals. To be convinced of that, one only needs to apply the same techniques as those used for showing Lemma \ref{continuite Y}.

\end{remark}
Remember that in section \ref{sec:hypotheses} we have denoted the eigenvalues of $A$ by $\lambda_1, \lambda_2$, $0 < \lambda_2 \le \lambda_1$. Also $0 < \lambda \le 1$ has been defined as the quotient of the eigenvalues, $\lambda := \ds{ \frac{\lambda_2}{\lambda_1}}$.

Let $P:=(v_1, v_2)$ be an orthonormal basis of eigenvectors of matrix $A$, such that $\lambda_1$ and $\lambda_2$ are their respectives eigenvalues.
The vector $v^{\star}$ can always be written in this basis:
\begin{equation*}
v^{\star}= \cos(\theta) v_1+ \sin(\theta) v_2.
\end{equation*}
It is always possible to choose $-\frac{\ds \pi}{\ds 2} < \theta \le \frac{\ds \pi}{\ds 2}$.

\begin{figure}[htbp]
\begin{center}
\includegraphics[width=.45\linewidth]{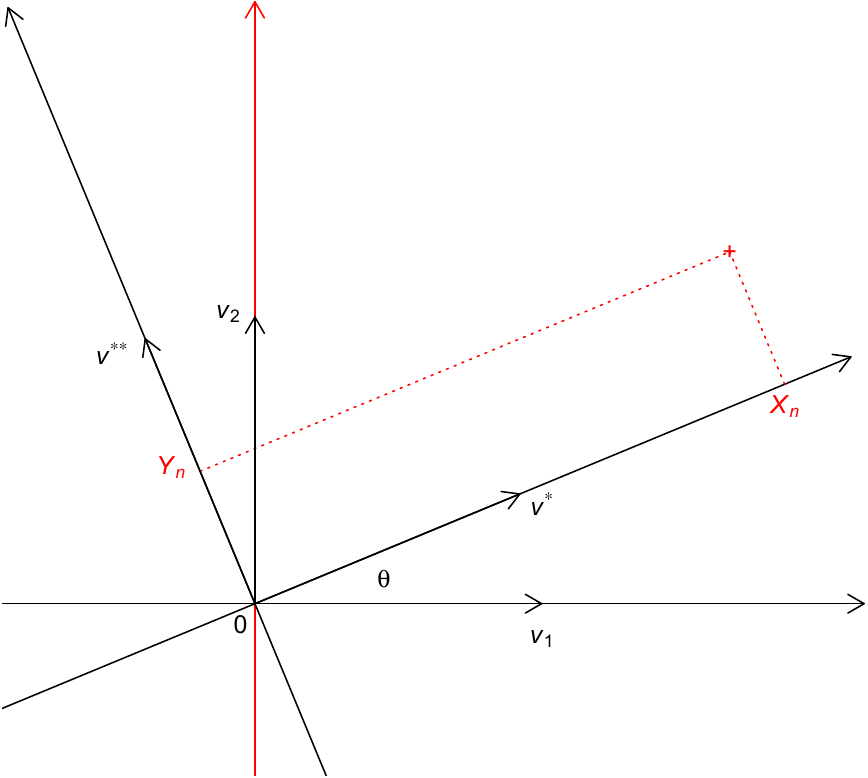}
\caption{Eigenvector directions}
\label{fig:figure1}
\end{center}
\end{figure}

Indeed, $\theta$ could be the angle between $v^\star$ and the eigenvector corresponding to the highest eigenvalue, because of the symmetry with respect to the point $(0, 0)$ and the fact that the mapping transforms  $\lambda$ into $\frac{1}{\lambda}$ and $\theta$ into $\frac{\pi}{2}-\theta$ has this effect.\\
By defining $v^{\star \star}:=\cos(\theta+\frac{\pi}{2}) \,v_1 +\sin(\theta+\frac{\pi}{2}) \,v_2= -\sin(\theta) \,v_1+\cos(\theta) \,v_2$, we thus defined a direct orthonormal basis $(v^{\star}, v^{\star \star})$, see Figure \ref{fig:figure1}.

Theorem \ref{convergence sure} applied to the particular functions $\overrightarrow{f^{\star}}$ and 
 to the function ${\bf 1}$ taking values in $\reels$ and identically equal to one implies that,
\begin{equation}
\label{convergence quotient esperance}
\frac{J_{\overrightarrow{f^{\star}}}^{(n)}(u)}{J_{\bf 1}^{(n)}(u)}\xrightarrow[n\to+\infty]{a.s.} \frac{\mathbb{E}[J_{\overrightarrow{f^{\star}}}^{(1)}(u)]}{\mathbb{E}[J_{\bf 1}^{(1)}(u)]}.
\end{equation}
Applying the Rice formula for this particular function $\overrightarrow{f^{\star}}$ and for the function ${\bf 1}$, we show the following proposition.
\begin{proposition}
\label{decomposition base}
\begin{align*}
\frac{J_{\overrightarrow{f^{\star}}}^{(n)}(u)}{J_{\bf 1}^{(n)}(u)}\xrightarrow[n\to+\infty]{a.s.} 
&\frac{1}{I(\lambda)} \left[{(\cos^{2}(\theta)+ \lambda^2 \sin^2(\theta))^{\frac{1}{2}}\, v^{\star} + \frac{\sin(\theta)\cos(\theta)(\lambda^2-1)}{(\cos^{2}(\theta)+ \lambda^2 \sin^2(\theta))^{\frac{1}{2}}} \,v^{\star \star}
}\right],  \nonumber 
\end{align*}
where $I(\lambda)$ is the elliptic integral
\begin{align}
\label{integrale elliptique}
I(\lambda):= \int_{0}^{\frac{\pi}{2}} (\cos^{2}(x)+ \lambda^2 \sin^2(x))^{\frac{1}{2}} \ud x.
\end{align}
\end{proposition}
\begin{remark}
Note that $J_{\bf 1}^{(n)}(u)=\ds  {\frac{ \sigma_{1}(\mathcal{C}_{n}(u))}{\sigma_{2}(T_n)}}$.
This functional is nothing but the measurement of the dimensional length of the level set by unit of surface.
\end{remark}
\begin{proofarg}{Proof of Proposition \ref{decomposition base}}
We use a lemma for which a proof is given in Appendix.
\begin{lemma}
\label{convergence esperance}
\begin{equation*}
\frac{\mathbb{E}[J_{\overrightarrow{f^{\star}}}^{(1)}(u)]}{\mathbb{E}[J_{\bf 1}^{(1)}(u)]}= \frac{2}{\pi} \frac{A^{2}.v^{\star}}{\normp[2]{A.v^{\star}}} \frac{1}{\int_{S^{1}}\normp[2]{A .\alpha} \ud\alpha},
\end{equation*}
where $\ud \alpha$ denotes the normalized area measure on $S^1$.
\end{lemma}
Using last equality we show that $\ds \frac{\mathbb{E}[J_{\overrightarrow{f^{\star}}}^{(1)}(u)]}{\mathbb{E}[J_{\bf 1}^{(1)}(u)]}$ can be expressed in the orthonormal basis $(v^{\star}, v^{\star \star})$ as follows.
One has:
\begin{align*}
A^2.v^{\star}= \lambda_{1}^2\left[{\cos(\theta)\,v_1+\lambda^2 \sin(\theta)\,v_2
}\right],
\end{align*}
and 
\begin{align*}
\normp[2]{A.v^{\star}}= \lambda_{1}\left[{\cos^2(\theta)+\lambda^2 \sin^2(\theta)
}\right]^{\frac{1}{2}}.
\end{align*}
Also
\begin{align*}
\int_{S^1}  \normp[2]{A.\alpha}\, \ud\alpha =\lambda_{1} \int_{0}^{2\pi} \left[{\cos^2(x)+\lambda^2 \sin^2(x)
}\right]^{\frac{1}{2}}\frac{\ud x}{2\pi}= \frac{2\lambda_1}{\pi}\, I(\lambda),
\end{align*}
where $I(\lambda)$ is the elliptic integral defined by (\ref{integrale elliptique}).\\
Thus by using Lemma \ref{convergence esperance} we proved that
\begin{align*}
\frac{\mathbb{E}[J_{\overrightarrow{f^{\star}}}^{(1)}(u)]}{\mathbb{E}[J_{\bf 1}^{(1)}(u)]} =
\frac{1}{I(\lambda)} \frac{\cos(\theta)\,v_1+\lambda^2 \sin(\theta)\,v_2}{\left[{\cos^2(\theta)+\lambda^2 \sin^2(\theta)
}\right]^{\frac{1}{2}}}.
\end{align*}
At this stage of the proof we are going to change the basis, that is expressing last identity in the orthonormal basis $(v^{\star}, v^{\star \star})$, getting
\begin{align}
\label{formule esperance}
\frac{\mathbb{E}[J_{\overrightarrow{f^{\star}}}^{(1)}(u)]}{\mathbb{E}[J_{\bf 1}^{(1)}(u)]} =
\frac{1}{I(\lambda)} \left[{(\cos^{2}(\theta)+ \lambda^2 \sin^2(\theta))^{\frac{1}{2}}\, v^{\star} + \frac{\sin(\theta)\cos(\theta)(\lambda^2-1)}{(\cos^{2}(\theta)+ \lambda^2 \sin^2(\theta))^{\frac{1}{2}}} \,v^{\star \star}
}\right].
\end{align}
Proposition \ref{decomposition base} ensues from convergence given in (\ref{convergence quotient esperance}).
\end{proofarg}

We consider function $\overrightarrow{F}$  being defined by:
\begin{equation}
\label{function F}
\begin{aligned} \overrightarrow{F}:\introf{0}{1} \times \introf{-\tfrac{\pi}{2}}{\tfrac{\pi}{2}}& \longrightarrow \{(X, Y) \in \reels^2, X >0, X^2+Y^2<1\}\\
(\lambda, \theta) & \mapsto \overrightarrow{F}(\lambda, \theta)= (F_1(\lambda, \theta), F_2(\lambda, \theta)),
\end{aligned}
\end{equation}
where 
$$   
\begin{cases}
F_1(\lambda, \theta)&:=  \ds \frac{1}{I(\lambda)}\, (\cos^{2}(\theta)+ \lambda^2 \sin^2(\theta))^{\frac{1}{2}}\\
F_2(\lambda, \theta)& :=   \ds \frac{1}{I(\lambda)}\, \ds\frac{\sin(\theta)\cos(\theta)(\lambda^2-1)}{(\cos^{2}(\theta)+ \lambda^2 \sin^2(\theta))^{\frac{1}{2}}}.
\end{cases}
$$
\label{F1 F2}
Proposition \ref{decomposition base} says that $\ds \frac{J_{\overrightarrow{f^{\star}}}^{(n)}(u)}{J_{\bf 1}^{(n)}(u)}\xrightarrow[n\to+\infty]{a.s.} F_1(\lambda, \theta)\, v^{\star} +F_2(\lambda, \theta)\, v^{\star \star}$.\\
In order to define estimators for parameters $\lambda$ and $\theta$ we need to look more closely at the properties of function $\overrightarrow{F}$. To this end we are going to follow step by step Wschebor's method given in \cite[pages 82 to 85]{MR0871689}. That is the object of next section.
\subsection{The $F$-diffeomorphism}
\label{le diffeomorphisme}
We prove the following proposition.
\begin{proposition}
\label{diffeomorphisme}
Function $\overrightarrow{F}$ is a one to one function from $]0, 1[ \times\introf{-\frac{\pi}{2}}{\frac{\pi}{2}}$ onto $\{(X, Y) \in \reels^2, X >0, X^2+Y^2<1 \mbox{\, and } (X \neq \frac{2}{\pi} \mbox{\, or } Y\neq0) \}$.\\
Furthermore $\overrightarrow{F}$ is a $C^2$-diffeomorphism from the open set $O:=\,]0, 1[ \times ]-\frac{\pi}{2}, \frac{\pi}{2}[$ onto the open set $O^{\prime}:=\{(X, Y) \in \reels^2, X >0, X^2+Y^2<1 \mbox{\, and } (X>\frac{2}{\pi} \mbox{\, or } Y\neq0) \}$.
\end{proposition}
\begin{proofarg}{Proof of Proposition \ref{diffeomorphisme}}
We reproduce here all the computations given by Wschebor in \cite{MR0871689} and also attach his two forthcoming pictures.\\
Let $(X, Y) \in \reels^2$ be fixed such that $X>0, X^2+Y^2<1$ and $X \neq \frac{2}{\pi}$ or $Y \neq 0$.\\
We consider the following system of equations 
\begin{align}
\label{system}
\begin{cases}
X &= F_1(\lambda, \theta)\\
Y & =   F_2(\lambda, \theta)
\end{cases}
\end{align}
If the system admits a solution $\lambda$, this solution ought to verify the following equation in $\lambda$:
\begin{equation}
\label{72}
X^2\, I^4(\lambda) \,(X^2+Y^2)-X^2\,I^2(\lambda)\,(\lambda^2+1)+\lambda^2=0.
\end{equation}
At this step, two cases appear: the case where $(X^2+Y^2)(\frac{\pi}{2})^2 -1 <0$ and the case where $(X^2+Y^2)(\frac{\pi}{2})^2 -1 \ge 0$.
Let us consider the first one.
\begin{enumerate}
\item $(X^2+Y^2)(\frac{\pi}{2})^2 -1 <0$.\\
Since for all $x \in \intrff{0}{1}$, $I(x) \le \frac{\pi}{2}$, the following inequality holds, $(X^2+Y^2)\,I^2(\lambda) -1 <0$.\\
Now dividing (\ref{72}) by $I^2(\lambda) X^2$, one gets
\begin{align*}
I^2(\lambda) \,(X^2+Y^2)-(\lambda^2+1)+\frac{\lambda^2}{I^2(\lambda) X^2}=0,
\end{align*}
that is, since $(X^2+Y^2)\,I^2(\lambda) -1 <0$, 
\begin{align}
\label{f1 f2}
f_1(\lambda)=f_2(\lambda),
\end{align}
where
\begin{align*}
f_1(\lambda)&:= X^2 \left(\frac{I(\lambda)}{\lambda}\right)^2\\
f_2(\lambda)&:= \frac{X^2\, I^2(\lambda)-1}{(X^2+Y^2)\,I^2(\lambda) -1}.
\end{align*}
Equation (\ref{f1 f2}) admits a unique solution in the interval $0 < \lambda <1$.
Let us argue this last assertion.
\\
Since $X>0$, function $f_1$ is strictly decreasing, while $f_2$ is such that
$f_2'(\lambda)= \ds \frac{2\, I(\lambda)\, I'(\lambda)}{((X^2+Y^2)\,I^2(\lambda) -1)^2}\, Y^2 >0$, if $Y \neq 0$.\\
Let us consider the case where $Y \neq 0$.
In this case function $f_2$ is strictly increasing.\\
By summarizing the situation we know that $f_1-f_2$ is continuous on $\introf{0}{1}$, strictly decreasing and such that
$(f_1-f_2)(0^+)=+\infty$ and \\
$(f_1-f_2)(1)=\ds \frac{X^2\, Y^2\, (\frac{\pi}{2})^4 + (X^2\, (\frac{\pi}{2})^2-1)^2}{(X^2+Y^2)\,(\frac{\pi}{2})^2 -1} <0$.
Thus there exists an unique $0< \lambda <1$, such that $f_1(\lambda)=f_2(\lambda)$, see Figure \ref{fig:figure2}.
\begin{figure}[htbp]
\begin{center}
\includegraphics[width=.45\linewidth]{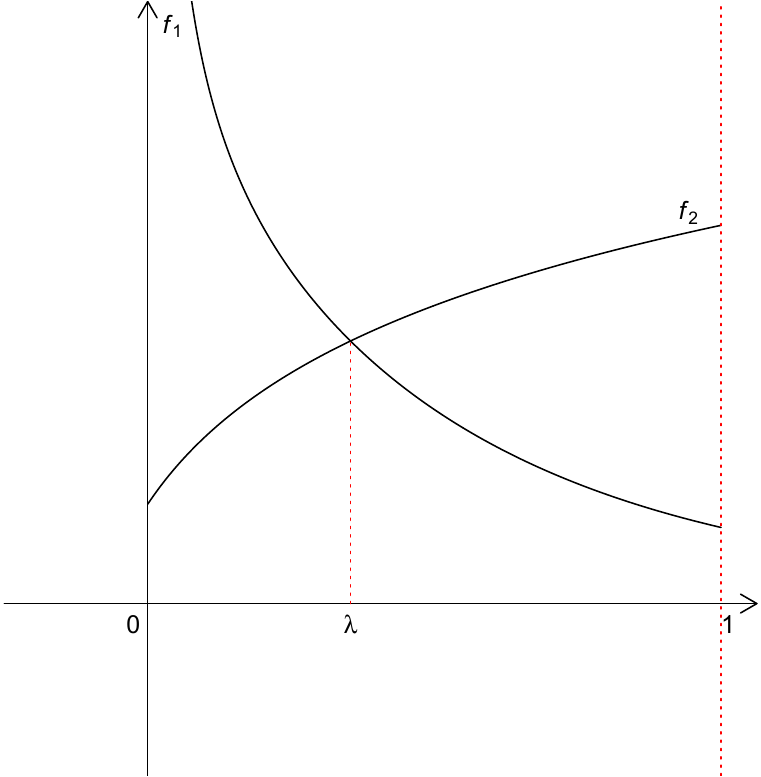}
\caption{$\lambda$-solution, first case}
\label{fig:figure2}
\end{center}
\end{figure}

Now, in the case where $Y:=0$, $f_2(\lambda)=1$.
We define $g$ by $g(\lambda):= \ds \frac{I(\lambda)}{\lambda}$.
Since $X < \frac{2}{\pi}$, $\lambda= g^{-1}(\frac{1}{X})$ is the unique solution in $\introo{0}{1}$ of (\ref{f1 f2}).

Let us consider now the second case.
\item $(X^2+Y^2)(\frac{\pi}{2})^2 -1 \ge 0$.

Since $X^2+Y^2 >0$, function $h: \lambda \to (X^2+Y^2)\, I^2(\lambda) -1$ is strictly increasing and continuous on $\introf{0}{1}$.\\
Moreover, $h(0)= X^2+Y^2-1 <0$ and $h(1)= (X^2+Y^2)(\frac{\pi}{2})^2 -1 \ge 0$.
Thus there exists a unique $0 < \lambda_0 \le 1$ such that
$h(\lambda_0)=(X^2+Y^2)\, I^2(\lambda_0) -1=0$ and then for $\lambda < \lambda_0$, one has $(X^2+Y^2)\, I^2(\lambda) -1 <0$ and for $\lambda > \lambda_0$, $(X^2+Y^2)\, I^2(\lambda) -1 >0$.\\
Arguing as in first part, we can deduce that if (\ref{system}) admits a solution $\lambda \neq \lambda_0$, this solution ought to verify $f_1(\lambda)=f_2(\lambda)$.\\
Function $f_1$ is strictly decreasing, while function $f_2$ is strictly increasing if $Y \neq 0$.\\
In this case $f_1-f_2$ is continuous on $\introf{0}{1}$, strictly decreasing and such that
$(f_1-f_2)(0^+)=+\infty$ and $(f_1-f_2)(\lambda_0^-)=-\infty$, since $Y \neq 0$.
Thus there exists an unique $0< \lambda < \lambda_0 \le1$, such that $f_1(\lambda)=f_2(\lambda)$.\\
On the other side, still if $Y \neq 0$ and if $\lambda_0 <1$, we have $(f_1-f_2)(\lambda_0^+)=+\infty$ and $(f_1-f_2)(1^-) >0$.
Then there is no more solution of (\ref{f1 f2}) into interval $]\lambda_0^+, 1]$.\\
Thus we have proved that in the case where $Y \neq 0$, if (\ref{system}) admits a solution $\lambda \neq \lambda_0$, this solution is the unique solution in $\introo{0}{1}$ of (\ref{f1 f2}) and this solution belongs to the interval $]0, \lambda_0[$ see Figure \ref{fig:figure3}.
\begin{figure}[htbp]
\begin{center}
\includegraphics[width=.45\linewidth]{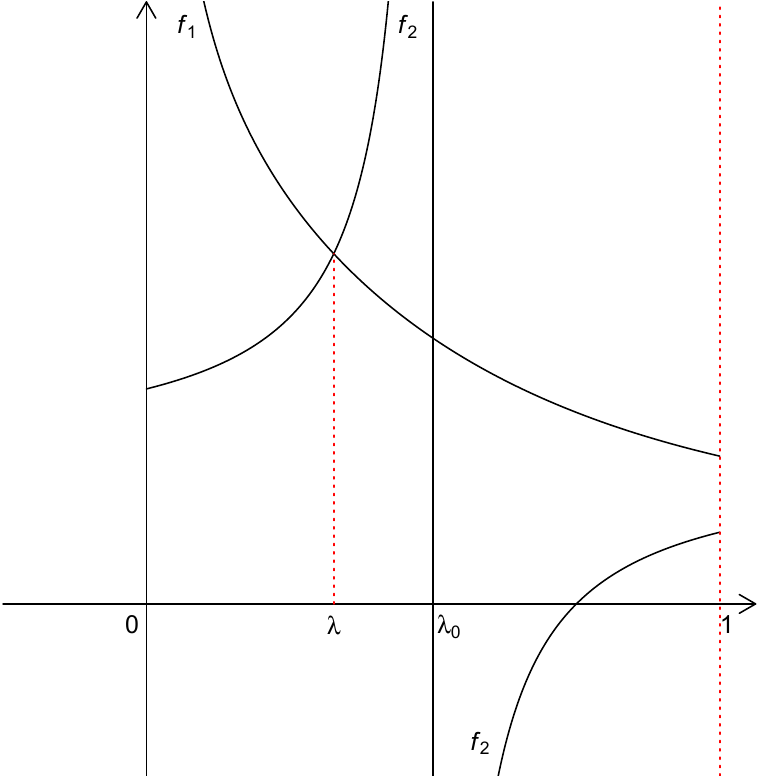}
\caption{$\lambda$-solution, second case}
\label{fig:figure3}
\end{center}
\end{figure}

Now we are going to prove that in the case where $Y=0$, there is no solution to (\ref{system}) different from $\lambda_0$.\\
In fact, if $\lambda$ is such a solution, that is different from $\lambda_0$, then $f_1(\lambda)=X^2 \left(\frac{I(\lambda)}{\lambda}\right)^2=f_2(\lambda)=1$, that is no possible.
Indeed, since $X^2 (\frac{\pi}{2})^2 -1 > 0$ then we would have $1 < X^2 (\frac{\pi}{2})^2 \le X^2\left(\frac{I(\lambda)}{\lambda}\right)^2=1$.

To finish the proof of this part of proposition, we look for conditions on $\lambda=\lambda_0$ to be solution of (\ref{system}).
If such a solution exists necessarily it has to verify (\ref{72}), and then
since $(X^2+Y^2)\, I^2(\lambda_0) -1=0$, we would have $\lambda_0^2 (1-X^2\, I^2(\lambda_0))=0$, so that $Y=0$ and $\lambda_0=I^{-1}(\frac{1}{X})$.\\
Thus, since $X^2 (\frac{\pi}{2})^2 >1$ and then $0< \lambda=\lambda_0 =I^{-1}(\frac{1}{X})<1$. \\
If we combine all results we get the following: if (\ref{system}) admits a solution $\lambda$ then, this solution is unique, satisfies $0 < \lambda <1$ and is such that:
\begin{itemize}
\item Or $(X^2+Y^2)(\frac{\pi}{2})^2 -1 <0$, and this solution verifies equation $f_1(\lambda)=f_2(\lambda)$.\\
Furthermore in the case where $Y=0$ then $\lambda=g^{-1}(\frac{1}{X})$, where we recall that $g(\lambda)=\ds \frac{I(\lambda)}{\lambda}$.
\item Or $(X^2+Y^2)(\frac{\pi}{2})^2 -1 \ge 0$ and $Y \neq 0$, and this solution verifies equation $f_1(\lambda)=f_2(\lambda)$.
\item  Or $X^2 (\frac{\pi}{2})^2 >1$ and $Y= 0$, and this solution is $\lambda =I^{-1}(\frac{1}{X})$.
\end{itemize}
\begin{remark}
\label{equation 72}
Note that in all cases $\lambda$ verifies (\ref{72}).
\end{remark}
\begin{remark}
\label{negatif}
Note that in all cases $(X^2+Y^2)\, I^2(\lambda) -1 \le 0$.
\end{remark}
\end{enumerate}
Now by continuing the reasoning by necessary condition, if we want that $\lambda<1$ verifies the first equation of (\ref{system}), we have to resolve the following equation:
\begin{align*}
X^2\, I^2(\lambda)= (\lambda^2-1)\, \sin^2(\theta)+1,
\end{align*}
that is
\begin{align*}
\sin^2(\theta)=\ds \frac{1-X^2\, I^2(\lambda)}{1-\lambda^2}=\ds \frac{\ds \frac{1}{\lambda^2}-f_1(\lambda)}{\ds \frac{1}{\lambda^2}-1}.
\end{align*}
Let us see that $0 \le \ds \frac{\ds \frac{1}{\lambda^2}-f_1(\lambda)}{\ds \frac{1}{\lambda^2}-1} \le 1$.
We have the following equivalence:
\begin{align*}
\left({\ds \frac{\frac{1}{\lambda^2}-f_1(\lambda)}{\frac{1}{\lambda^2}-1} \le 1}\right) \iff \left({f_1(\lambda) \ge 1}\right).
\end{align*}
In the case where $X^2 (\frac{\pi}{2})^2 >1$ and $Y = 0$, we have $\lambda =I^{-1}(\frac{1}{X})$. Thus $f_1(\lambda)= \ds \frac{1}{\lambda^2} \ge 1$,
since $0 < \lambda <1$.\\
In the other cases, $f_1(\lambda)=f_2(\lambda)$ and since $f_2$ is an increasing function, $f_1(\lambda)\ge f_2(0^+)= \ds \frac{X^2\, -1}{X^2+Y^2\,-1} \ge 1$, since $Y^2 \ge 0$ and $X^2+Y^2\,-1<0$.\\
Now $\left({\ds \frac{ \frac{1}{\lambda^2}-f_1(\lambda)}{ \frac{1}{\lambda^2}-1} \ge 0}\right) \iff \left({X^2\, I^2(\lambda) -1\le 0}\right)$.\\
Remark \ref{negatif} gives the result.\\
If we look at the sign of expression $F_2(\lambda, \theta)$, we have to set:
\begin{align*}
\theta:= \arcsin\left({\sqrt{\frac{1-X^2 I^2(\lambda)}{1-\lambda^2}}
 }\right)\left({ \1_{\{Y \le 0\}} - \1_{\{Y > 0\}}
 }\right).
\end{align*}
In this way, we always have $-\frac{\pi}{2} \le \theta \le \frac{\pi}{2}$ and in fact $\theta > -\frac{\pi}{2}$.
To verify last inequality, remark that
\begin{align*}
\left({\theta=-\frac{\pi}{2}}\right) \iff \left({Y>0 \mbox{ and \,}  X^2 I^2(\lambda)=\lambda^2}\right)\\
\Longrightarrow \left({Y>0  \mbox{ and \,} f_1(\lambda)=f_2(\lambda)=1}\right)
\Longrightarrow \left({Y>0 \mbox{ and \,} Y=0}\right),
\end{align*}
which is impossible.\\
Let us see now that we can go backward.
That is let us verify that these proposed $\lambda (<1)$ and $\theta$
are actually solutions of (\ref{system}).
\begin{align*}
\left({\sin^2(\theta)=\ds \frac{1-X^2\, I^2(\lambda)}{1-\lambda^2}}\right) \iff \left({X^2= \frac{\cos^2(\theta)+ \lambda^2 \sin^2(\theta)}{I^2(\lambda)} }\right)\\
\iff \left({X= \frac{(\cos^2(\theta)+ \lambda^2 \sin^2(\theta))^{\frac{1}{2}}}{I(\lambda)} =F_1(\lambda,\theta)}\right),
\end{align*}
last equivalence provides from the fact that $X>0$.

Let us consider the second equation of (\ref{system}).
In the case where $Y=0$ and $X> \ds \frac{2}{\pi}$, one has $\lambda=I^{-1}(\frac{1}{X})$ so that $\theta=\arcsin(0)=0$ and $Y=0=F_2(\lambda, 0)$.
In the other cases, one has $f_1(\lambda)= f_2(\lambda)$, and since $X^2\, I^2(\lambda)=\cos^2(\theta)+ \lambda^2 \sin^2(\theta)$ one has
\begin{align*}
\lefteqn{\left({f_1(\lambda)= f_2(\lambda)}\right)}\\
& \iff \left({\frac{X^2\, I^2(\lambda)}{\lambda^2}= \frac{X^2\, I^2(\lambda)-1}{(X^2+Y^2)\,I^2(\lambda) -1}}\right)\\
& \iff \left({Y^2=  \left[{\frac{1}{I(\lambda)}\, \ds\frac{\sin(\theta)\cos(\theta)(\lambda^2-1)}{(\cos^{2}(\theta)+ \lambda^2 \sin^2(\theta))^{\frac{1}{2}}} }\right]^2}\right)\\
& \iff \left({Y=  \frac{1}{I(\lambda)}\, \ds\frac{\sin(\theta)\cos(\theta)(\lambda^2-1)}{(\cos^{2}(\theta)+ \lambda^2 \sin^2(\theta))^{\frac{1}{2}}}=F_2(\lambda, \theta) }\right),
\end{align*}
last equivalence providing from the fact that if $Y \le 0$, $0 \le \theta \le \frac{\pi}{2}$ and if $Y >0$, $-\frac{\pi}{2} \le \theta \le 0$.\\
We thus proved that function $\overrightarrow{F}$ is a one to one function from $]0, 1[ \times\introf{-\frac{\pi}{2}}{\frac{\pi}{2}}$ onto $\{(X, Y) \in \reels^2, X >0, X^2+Y^2<1 \mbox{\, and } (X \neq \frac{2}{\pi} \mbox{\, or } Y\neq0) \}$. This yields the first part of proposition.
\begin{remark}
\label{equation isotropy}
Note that $\left({\overrightarrow{F}(\lambda, \theta)=(\ds \frac{2}{\pi}, 0)}\right) \iff \left({\lambda=1, -\frac{\pi}{2} <\theta \le \frac{\pi}{2}}\right)$.\\
\end{remark}
Now the following equivalence 
\begin{align*}
\left({\lambda <1, \theta=\frac{\pi}{2}}\right) \iff \left({Y=0 \mbox{ and \,}  X < \frac{2}{\pi}}\right),
\end{align*}
leads to the conclusion that $\overrightarrow{F}$ is a one to one function from the open set $O:= ]0, 1[ \times\introo{-\frac{\pi}{2}}{\frac{\pi}{2}}$ onto the open set $O^{\prime}:=\{(X, Y) \in \reels^2, X >0, X^2+Y^2<1 \mbox{\, and } (X > \frac{2}{\pi} \mbox{\, or } Y\neq0) \}$. The second part of proposition follows from Lemma \ref{jacobien} proved in Appendix which ensures that the Jacobian of the tranformation $\overrightarrow{F}$ is different from zero since $0< \lambda <1$. 
\end{proofarg}
We are now listen to give the definitions of the estimators of parameters $\lambda$ and $\theta$.
\subsection{Definition of the affinity parameters $\widehat{\lambda}_{n}$ and $\widehat{\theta}_{n}$}
\label{les parametres d'affinite}
As in \cite[chap. 3.6 F]{MR0871689} we write the {\bf observed ratio} of functionals $\ds \frac{J_{\overrightarrow{f^{\star}}}^{(n)}(u)}{J_{\bf 1}^{(n)}(u)}$ as
\begin{align}
\label{Xn Yn}
\ds \frac{J_{\overrightarrow{f^{\star}}}^{(n)}(u)}{J_{\bf 1}^{(n)}(u)}= X_n v^{\star} + Y_n v^{\star \star}.
\end{align}
As a first step we use a lemma proved in Appendix \ref{ann:proof}.
\begin{lemma}
\label{bornes}
The random variables
$X_n$ and $Y_n$ defined in (\ref{Xn Yn}) are such that for $n \in \naturels^{\star}$, a.s.
$X_n >0$ and $X_n^2 +Y_n^2 < 1$.
\end{lemma}
The first part of Proposition \ref{diffeomorphisme}, Lemma \ref{bornes} and Remark \ref{equation isotropy} provide the following estimators $\widehat{\lambda}_{n}$ and $\widehat{\theta}_{n}$  of parameters $\lambda$ and $\theta$ via the following proposition.
\begin{proposition}
\label{def estimateurs}
Let consider the following system of equations 
$$   
\begin{cases}
X_n &= F_1(\lambda, \theta)\\
Y_n & =   F_2(\lambda, \theta),
\end{cases}
$$
where $\overrightarrow{F}=(F_1, F_2)$ has been defined by (\ref{function F}).
\begin{enumerate}
\item In the case where $Y_n \neq 0$ or $X_n \neq \frac{2}{\pi}$, it admits a unique solution $(\widehat{\lambda}_{n}, \widehat{\theta}_{n})$ such that $0 < \widehat{\lambda}_{n} < 1$ and $-\frac{\pi}{2} < \widehat{\theta}_{n} \le \frac{\pi}{2}$.
\begin{itemize}
\item The estimator  $0 <\widehat{\lambda}_{n} < 1$ is solution of
\begin{align*}
X_n^2 \left(\frac{I(\lambda)}{\lambda}\right)^2= \frac{X_n^2 I^2(\lambda)-1}{(X_n^2+Y_n^2)I^2(\lambda)-1},
\end{align*}
\begin{itemize}
\item except for $Y_n=0$ and $X_n > \frac{2}{\pi}$, where $\widehat{\lambda}_{n}:= I^{-1}(\frac{1}{X_n})$.
\end{itemize}
\item The estimator  $\widehat{\theta}_{n}$ is defined as
\begin{align*}
\widehat{\theta}_{n}:= \arcsin\left({\sqrt{\frac{1-X_n^2 I^2(\widehat{\lambda}_n)}{1-\widehat{\lambda}_n^2}}
 }\right)\left({ \1_{\{Y_n \le 0\}} - \1_{\{Y_n > 0\}}
 }\right).
\end{align*}
\end{itemize}
\item In the case where $Y_n=0$ and $X_n=\frac{2}{\pi}$, it admits a unique solution $\widehat{\lambda}_{n}:=1$, $\widehat{\theta}_{n}$ being any number belonging to $]-\frac{\pi}{2}, \frac{\pi}{2}]$.
\end{enumerate}
\end{proposition}
Now we will apply the result of almost sure convergence highlighted in Proposition \ref{decomposition base} and the fact that function $\overrightarrow{F}$ is a diffeomorphism. We will deduce the almost sure convergence of the estimators of affinity parameters. Then we will use the law convergence result established for general functionals in Corollary \ref{convergence triple} in the aim to establish a law convergence result for $\frac{J_{\overrightarrow{f^{\star}}}^{(n)}(u)}{J_{\bf 1}^{(n)}(u)}$. Once this is done we will deduce some law convergence results for the estimators $\widehat{\lambda}_{n}$ and $\widehat{\theta}_{n}$ and also propose some confidence intervals for the parameters $\lambda$ and $\theta$.

\subsection{Convergence for the affinity estimators}
\label{convergence pour les parametres d'affinite}
\subsubsection{Consistency for the estimators}
We are ready to state the following results of consistency for the two proposed estimators $\widehat{\lambda}_{n}$ and $\widehat{\theta}_{n}$ of parameters $\lambda$ and $\theta$.
\begin{theorem}
\label{convergence presque-sure}
For $0 < \lambda <1$ and $-\frac{\pi}{2} < \theta < \frac{\pi}{2}$, one has
\begin{equation*}
 \widehat{\lambda}_{n} \xrightarrow[n\to+\infty]{a.s.} \lambda \mbox{\, and \, } \widehat{\theta}_{n} \xrightarrow[n\to+\infty]{a.s.} \theta.
\end{equation*}
\end{theorem}
\begin{proofarg}{Proof of Theorem \ref{convergence presque-sure}}
By the second part of Proposition \ref{diffeomorphisme} we know that $F$ is a $C^2$-diffeomorphism from the open set $O=\,\introo{0}{1} \times ]-\frac{\pi}{2}, \frac{\pi}{2}[$ onto the open set $O^{\prime}=\{(X, Y) \in \reels^2, X>0, X^2+Y^2 <1, 
 (Y \neq 0 \mbox{\, or \,} X > \frac{2}{\pi})\}$.\\
 Furthermore by using the result convergence given for $\ds \frac{J_{\overrightarrow{f^{\star}}}^{(n)}(u)}{J_{\bf 1}^{(n)}(u)}$ in Proposition \ref{decomposition base}, decomposition of this last \rv given in (\ref{Xn Yn}), we obtain the following convergence result: for all $0< \lambda \le 1$ and $-\frac{\pi}{2} <\theta \le \frac{\pi}{2}$,
 \begin{equation}
 \label{convergence presque sure}
\left({X_n, Y_n }\right) \xrightarrow[n\to+\infty]{a.s.} (X, Y):=\overrightarrow{F}(\lambda, \theta).
\end{equation}
If we reduce the definition domain of $(\lambda, \theta)$, that is if we suppose that $0< \lambda < 1$ and $-\frac{\pi}{2} < \theta< \frac{\pi}{2}$, then $(X, Y)$ belongs to the open set $O^{\prime}$.
If $\overrightarrow{F}^{-1}$ denotes the inverse function of $\overrightarrow{F}$, thus since $\overrightarrow{F}^{-1}$ is continuous from $O^{\prime}$ to $O$, one deduces that almost surely $\overrightarrow{F}^{-1}(X_n, Y_n)=(\widehat{\lambda}_{n}, \widehat{\theta}_{n})$ converges to $\overrightarrow{F}^{-1}(X, Y)=(\lambda, \theta)$.
This yields Theorem \ref{convergence presque-sure}.
\end{proofarg}

We have all the elements to prove a limit theorem about the distributions of estimators $\widehat {\lambda} _ {n} $ and $\widehat {\theta} _ {n} $.
Consequently, confidence intervals for the estimated parameters $\lambda $ and $\theta $ can be proposed.

\subsubsection{Convergence in law for the estimators}
First, we apply results of Section \ref{section xi}.\\
Recalling that $\overrightarrow{f^{\star}}=(f_1^{\star}, f_2^{\star}) $ has been defined in (\ref{fstar}). We define  $\overrightarrow{f_{\star}}:=(f_1^{\star}, f_2^{\star}, f_3^{\star})$, with $f^{\star}_3:={\bf 1}$. By Theorem \ref{convergence sure} and Corollary \ref{convergence triple} one can show the following proposition.
\begin{proposition}
\label{convergence en loi}
\begin{align*}
2n\left({
 \tfrac{\ds J_{\overrightarrow{f^{\star}}}^{(n)}(u)}{\ds J_{\bf 1}^{(n)}(u)} - \ds \tfrac{\mathbb{E}[\ds J_{\overrightarrow{f^{\star}}}^{(1)}(u)]}{\mathbb{E}[\ds J_{\bf 1}^{(1)}(u)]}
}\right) 
\xrightarrow[n \to +\infty]{Law} \mathcal{N}(0; \mathit{\Sigma^{\star}}(u)),
\end{align*}
where $\mathit{\Sigma^{\star}}(u)= \mathit{\Sigma^{\star}}(u, \lambda_1, \lambda_2, P):= B \mathit{\Sigma}_{\overrightarrow{f_{\star}}}(u)\, B^{t}$ and 
$$
B:=\ds \frac{\ds 1}{\ds a_{\bf 1}({\bf 0}, u)} \begin{pmatrix}1 & 0& - \frac{\ds a_{f_{1}^{\star}}({\bf 0}, u)}{\ds a_{\bf 1}({\bf 0}, u)} \\
0 & 1 & - \frac{\ds a_{f_{2}^{\star}}({\bf 0}, u)}{\ds a_{\bf 1}({\bf 0}, u)} \end{pmatrix},
$$
where 
$
	\mathit{\Sigma}_{\overrightarrow{f_{\star}}}(u)$ being defined by (\ref{Sigma vectoriel}).
\end{proposition}
\begin{remark}
\label{calculus}
The expression of coefficients $a_{f_{i}^{\star}}(\cdot, u)$ for $i=1, 2$ and of $a_{\bf 1}(\cdot, u)$ are respectively given in Lemma \ref{calculus of c} and Lemma \ref{calculus of a} of Appendix \ref{ann:proof}.
\end{remark}
\begin{remark}
\label{non degenerate}
The asymptotic variance matrix $\mathit{\Sigma^{\star}}(u)$ is a non-degenerate matrix.
\end{remark}
\begin{proofarg}{Proof of Proposition \ref{convergence en loi}}
Let $f: S^1 \to \reels$ be a continuous and bounded function.
Since $\mathbb{E}[J_{f}^{(n)}(u)] =a_f({\bf 0}, u)$ the following decomposition ensues
\begin{multline*}
2n\left({
\tfrac{\ds J_{f}^{(n)}(u)}{\ds J_{\bf 1}^{(n)}(u)} - \tfrac{\ds \mathbb{E}[J_{f}^{(n)}(u)]}{\ds \mathbb{E}[J_{\bf 1}^{(n)}(u)]}
}\right) =
\ds \tfrac{\ds 1}{\ds a_{\bf 1}({\bf 0}, u)} \left({\xi_{f}^{(n)}(u) - \tfrac{\ds a_{f}({\bf 0}, u)}{\ds a_{\bf 1}({\bf 0}, u)}\, \xi_{\bf 1}^{(n)}(u)}\right)\\
+ \left({\tfrac{\ds 1}{\ds J_{\bf 1}^{(n)}(u)} - \tfrac{\ds 1}{\ds a_{\bf 1}({\bf 0}, u)}}\right)  \left({\xi_{f}^{(n)}(u) - \tfrac{\ds a_{f}({\bf 0}, u)}{\ds a_{\bf 1}({\bf 0}, u)}\, \xi_{\bf 1}^{(n)}(u)}\right)\\
\equiv \tfrac{\ds 1}{\ds a_{\bf 1}({\bf 0}, u)} \left({\xi_{f}^{(n)}(u) - \tfrac{\ds a_{f}({\bf 0}, u)}{\ds a_{\bf 1}({\bf 0}, u)}\, \xi_{\bf 1}^{(n)}(u)}\right),
\end{multline*}
the last law equivalence providing from Theorem \ref{convergence sure} and from Theorem \ref{Peccati}.\\
Applying this reasoning again successively to $f:=f_1^{\star}$ and $f:=f_2^{\star}$, and using Corollary \ref{convergence triple} by taking $f_1:=f_1^{\star}$, $f_2:=f_2^{\star}$ and $f_3:={\bf 1}$, we get Proposition \ref{convergence en loi}.
\end{proofarg}
\begin{remark}
\label{two dimensional}
Note that this proof highlights the fact that the asymptotic behavior of \begin{align*}
2n\left({
 \tfrac{\ds J_{\overrightarrow{f^{\star}}}^{(n)}(u)}{\ds J_{\bf 1}^{(n)}(u)} - \ds \tfrac{\mathbb{E}[\ds J_{\overrightarrow{f^{\star}}}^{(1)}(u)]}{\mathbb{E}[\ds J_{\bf 1}^{(1)}(u)]}
}\right) 
\end{align*} is the same as that of $\xi_{\overrightarrow{f}}^{(n)}(u)$, where $\overrightarrow{f}:=(f_1, f_2)$, with $f_i:=\frac{\ds f_i^{\star}}{\ds a_{\bf 1}({\bf 0}, u)}-\frac{\ds a_{f_i^{\star}}({\bf 0}, u)}{\ds a^2_{\bf 1}({\bf 0}, u)}\,{\bf 1}$, $i=1, 2$.
\end{remark}
\begin{proofarg}{Proof of Remark \ref{non degenerate}}
It remains to prove that the matrix $\mathit{\Sigma^{\star}}(u)$ is positive definite.
In this aim, for $f_1, f_2:S^1 \to \reels$ continuous and bounded functions, let us note $\overrightarrow{f}:=(f_1, f_2)$ and for $q \in \naturels^{\star}$,
\begin{align*}
\left({\mathit{\Sigma}_{f_1, f_2}(u)}\right)_q:=\sum_{
\substack{
{\bk, \bm} \in \naturels^3 \\ \abs{{\bk}}= \abs{{\bm}}=q
}
} a_{f_1}({\bk}, u)\, a_{f_2}({\bm}, u)\, R({\bk}, {\bm}),
\end{align*}
and 
\begin{align*}
\left({ \mathit{\Sigma}_{\overrightarrow{f}}(u)}\right)_q:= \left({\left({\mathit{\Sigma}_{f_i, f_j}(u)}\right)_q}\right)_{1\le i, j \le 2}.
\end{align*}
The proof of Lemma \ref{determinant} is given in Appendix \ref{ann:proof}.
\begin{lemma}
\label{determinant}
For all $f_1, f_2:S^1 \to \reels$ continuous and bounded functions and for all $q \in \naturels^{\star}$, one has
\begin{align*}
\det\left({\mathit{\Sigma}_{\overrightarrow{f}}(u)}\right)\ge \sum_{q=1}^{+\infty} \det\left({ \mathit{\Sigma}_{\overrightarrow{f}}(u)}\right)_{q} \ge \det\left({ \mathit{\Sigma}_{\overrightarrow{f}}(u)}\right)_q.
\end{align*}
\end{lemma}
Thus applying this lemma to functions $f_i$ previously defined in Remark \ref{two dimensional} by
$f_i:=\frac{\ds f_i^{\star}}{\ds a_{\bf 1}({\bf 0}, u)}-\frac{\ds a_{f_i^{\star}}({\bf 0}, u)}{\ds a^2_{\bf 1}({\bf 0}, u)}\,{\bf 1}$, for $i=1, 2$, 
we obtain
\begin{align*}
\det\left({\mathit{\Sigma^{\star}}(u)}\right) \ge \det\left({\mathit{\Sigma}_{\overrightarrow{f}}(u)}\right)_2.
\end{align*}
Note that $\det\left({\mathit{\Sigma}_{\overrightarrow{f}}(u)}\right)_1=0$ and this is why we let $q=2$ in Lemma \ref{determinant}.
Using arguments similar to those given in the proof of Remark \ref{non degeneree}, simple calculations give
\begin{align*}
\det\left({\mathit{\Sigma^{\star}}(u)}\right) \ge 4 \times (2\pi)^2 \int_{\reels^4} f^2_x(s) f^2_x(t) \left\{{g^2(s, t) -g(s, t)g(t, s) }\right\} \ud s\, \ud t \ge 0,
\end{align*}
where 
\begin{align*}
g(s, t):= s_1^2t_2^2\, a + s_1^2t_1t_2\, b- s_1s_2 t_2^2\,c,
\end{align*}
with
\begin{align*}
a&:= \left[{(d_{11}d_{22}+d_{12}d_{21})\,A + d_{11}d_{12}\, B-d_{21}d_{22}\, C }\right] \det(D)\\
b&:= \left[{2d_{11}d_{21} A +d_{11}^2\, B -d_{21}^2\, C}\right] \det(D)\\
c&:= \left[{-2d_{12}d_{22}\, A-d_{12}^2\, B + d_{22}^2\, C}\right] \det(D),
\end{align*}
while
\begin{align*}
A&:= a_{f_1}((2, 0, 0), u)\, a_{f_2}((0, 2, 0), u)-  a_{f_1}((0, 2, 0), u)\, a_{f_2}((2, 0, 0), u)\\
B&:=  a_{f_1}((2, 0, 0), u)\, a_{f_2}(( 1, 1, 0), u)-  a_{f_1}((1, 1, 0), u)\, a_{f_2}((2, 0, 0), u)\\
C&:=  a_{f_1}((0, 2, 0), u)\, a_{f_2}((1, 1, 0), u)-  a_{f_1}((1, 1, 0), u)\, a_{f_2}((0, 2, 0), u),
\end{align*}
and $D$ is defined as in the end of the proof of Remark \ref{non degeneree} as $D=(d_{ij})_{1\le i, j \le 2}:=\ds \frac{1}{\sqrt{\mu}}\, \Lambda^{-1}P^t$.
The positivity of the last integral provides from an application of H\"{o}lder inequality.\\
Let us remark on the one hand that the nullity of the integral is equivalent to the equality in H\"{o}lder inequality and thus to $a=b=c=0$.
On the other hand since $\det(D) \neq 0$, $A, B, C$ is solution of a Cramer linear system.
We deduce that $\mathit{\Sigma^{\star}}(u)$ will be strictly positive if $(A, B, C) \neq (0, 0, 0)$.\\
Let us see that $C \neq 0$.\\
By using Lemmas \ref{calculus of c} and \ref{calculus of a}, straightforward computations show that,
\begin{align*}
C= -\frac{\det(P)}{a_{\bf 1}^2(0, 0)} \times \frac{\lambda_1 \lambda_2 \mu}{\pi (\lambda_1^2 (\omega^{\star}_1)^2+\lambda_2^2 (\omega^{\star}_2)^2)} \times \left[{\lambda_2^2 (\omega^{\star}_2)^2 + 2 (\lambda_1^2 (\omega^{\star}_1)^2-\lambda_2^2 (\omega^{\star}_2)^2) W}\right],
\end{align*}
where
\begin{align*}
W:= \frac{\sum_{n=0}^{+\infty}\frac{1}{4(n+1)}V_n}{\sum_{n=0}^{+\infty} V_n},
\mbox{ while for  }  n \in \naturels,
V_n:= \frac{(2n)!(2n+1)!}{(n!)^4 2^{4n+1}} (1-\lambda^2)^n.
\end{align*}
By using that $0 < W \le \frac{1}{4}$, one easily gets that $\lambda_2^2 (\omega^{\star}_2)^2 + 2 (\lambda_1^2 (\omega^{\star}_1)^2-\lambda_2^2 (\omega^{\star}_2)^2) W >0$, thus $C \neq 0$.\\
This yields proof of remark.
\end{proofarg}

We translate the convergence result expressed in the above Proposition \ref{convergence en loi} in the basis $(v^\star, v^{\star \star})$, recalling that function $\overrightarrow{F}$ has been defined in (\ref{function F}). \\
The decomposition given in (\ref{Xn Yn}), Proposition \ref{def estimateurs} and equality (\ref{formule esperance}) imply that for $0 < \lambda \le 1$ and $-\frac{\pi}{2} < \theta \le \frac{\pi}{2}$, 
\begin{multline*}
2n\left({
\frac{J_{\overrightarrow{f^{\star}}}^{(n)}(u)}{J_{\bf 1}^{(n)}(u)} - \frac{\mathbb{E}[J_{\overrightarrow{f^{\star}}}^{(n)}(u)]}{\mathbb{E}[J_{\bf 1}^{(n)}(u)]}
}\right) =
2n \left({
F_1(\widehat{\lambda}_{n} , \widehat{\theta}_{n}) - F_1(\lambda, \theta)
}\right) v^{\star} \\
+2n \left({
F_2(\widehat{\lambda}_{n} , \widehat{\theta}_{n}) - F_2(\lambda, \theta)
}\right) v^{\star \star}
\xrightarrow[n \to +\infty]{Law} \mathcal{N}(0; \mathit{\Sigma^{\star}}(u)).
\end{multline*} 
Thus we get for $0 < \lambda \le 1$ and $-\frac{\pi}{2} < \theta \le \frac{\pi}{2}$
\begin{align}
\label{convergence F}
2n \left({
\overrightarrow{F}(\widehat{\lambda}_{n} , \widehat{\theta}_{n}) - \overrightarrow{F}(\lambda, \theta)
}\right)
\xrightarrow[n \to +\infty]{Law} \mathcal{N}(0; \mathit{\Sigma_Q^{\star}}(u),
\end{align}
where $\mathit{\Sigma_Q^{\star}}(u)= \mathit{\Sigma_Q^{\star}}(u, \lambda_1, \lambda_2, P) :=Q \times \mathit{\Sigma^{\star}}(u) \times Q^{t}$ and $Q$ is the change of basis matrix from the canonical basis $(\vec{i}, \vec{j})$ to the basis $(v^\star, v^{\star \star})$.\\
By using last result convervence and the fact that $\overrightarrow{F}$ is a $C^2$-diffeomorphism, we get the following 
theorem.
\begin{theorem}
\label{convergence lambda et theta}
For $0 < \lambda < 1$ and $-\frac{\pi}{2} < \theta < \frac{\pi}{2}$, one has
\begin{align*}
2n \left({
\widehat{\lambda}_{n} - \lambda; \,\widehat{\theta}_{n} - \theta
}\right)
\xrightarrow[n \to +\infty]{Law} \mathcal{N}(0; \mathit{\Sigma_{\lambda, \theta}}(u)),
\end{align*}
where $\mathit{\Sigma_{\lambda, \theta}}(u):=C(\lambda, \theta) \times \mathit{\Sigma_Q^{\star}}(u) \times C^{t}(\lambda, \theta)$ and
\begin{align*}
C(\lambda, \theta):= \frac{1}{J_{\overrightarrow{F}}(\lambda, \theta)}
\begin{pmatrix}
\dfrac{\partial F_2}{\partial \theta}(\lambda, \theta) & -\dfrac{\partial F_1}{\partial \theta}(\lambda, \theta) \\
&\\
-\dfrac{\partial F_2}{\partial \lambda}(\lambda, \theta) &\dfrac{\partial F_1}{\partial \lambda}(\lambda, \theta)
\end{pmatrix}
\end{align*}
where the Jacobian $J_{\overrightarrow{F}}$ has been defined in Appendix (\ref{jacobien F}).
\end{theorem}
\begin{proofarg}{Proof of Theorem \ref{convergence lambda et theta}}
We supposed that $(\lambda, \theta) \in O= ]0, 1[ \times ]-\frac{\pi}{2}, \frac{\pi}{2}[$.
On the one hand, by Proposition \ref{diffeomorphisme}, we know that  $\overrightarrow{F}$ is a $C^2$-diffeomorphism from the open set $O$ onto the open set $O^{\prime}= \overrightarrow{F}(O)$.\\
On the other hand, by Theorem \ref{convergence presque-sure}, we know that 
\begin{align*}
(\widehat{\lambda}_{n}, \widehat{\theta}_{n}) \xrightarrow[n\to+\infty]{a.s.} (\lambda, \theta).
\end{align*}
Thus for $n$ large enough, almost surely $\overrightarrow{F}(\widehat{\lambda}_{n} , \widehat{\theta}_{n})$ and  $\overrightarrow{F}(\lambda, \theta)$ belong to $O^{\prime}$.
Using a second order Taylor-Young expansion of $\overrightarrow{F}^{-1}$ about $\overrightarrow{F}(\lambda, \theta)$, we get
\begin{multline*}
2n \left({
\widehat{\lambda}_{n} - \lambda; \,\widehat{\theta}_{n} - \theta
}\right)=\\
\sum_{j=1}^{2} \dfrac{\partial \overrightarrow{F}^{-1}}{\partial x_j}(\overrightarrow{F}(\lambda, \theta)) \times 2n(F_j(\widehat{\lambda}_{n} , \widehat{\theta}_{n}) -F_j(\lambda, \theta)) +\\
 \frac{1}{2} \sum_{j, k=1}^2 \dfrac{\partial^2 \overrightarrow{F}^{-1}}{\partial x_j \partial x_k}(\overrightarrow{F}(\lambda, \theta)) \times 2n(F_j(\widehat{\lambda}_{n} , \widehat{\theta}_{n}) -\\F_j(\lambda, \theta)) \times (F_k(\widehat{\lambda}_{n} , \widehat{\theta}_{n}) -F_k(\lambda, \theta)) \\
 + o(2n \normp[2]{\overrightarrow{F}(\widehat{\lambda}_{n} , \widehat{\theta}_{n}) -\overrightarrow{F}(\lambda, \theta)}^2).
\end{multline*}
Law convergence result expressed in (\ref{convergence F}) gives Theorem \ref{convergence lambda et theta}.
\end{proofarg}

Now we are in position to propose confidence intervals for $(\lambda, \theta)$ when the covariance $r_z$ of $Z$ is known.
\subsubsection{Confidence intervals for the affinity parameters}
In this section, we suppose that parameters $\lambda$ and $\theta$ are such that $0 < \lambda <1$, $-\frac{\pi}{2} < \theta < \frac{\pi}{2}$.
We will also assume that the covariance function $r_z$ is a {\bf known function}.

One can build confidence intervals for parameters $(\lambda, \,\theta)$.
We will show that
\begin{align}
\label{Sigma etoile}
\mathit{\Sigma^{\star}}(u)= \mathit{\Sigma^{\star}}(u, \lambda_1, \lambda_2, P)= \ds \frac{1}{\lambda_1^2}\, \mathit{\Sigma_{\star}}(u, \lambda, P),
\end{align} where $\mathit{\Sigma_{\star}}$ is a continuous matrix as function of $(\lambda, P)$ and is computable provided that $(\lambda, P)$ are given.

We then consider $\widehat{\lambda}_{1, n}$ and $\widehat{P}_n$ two estimators of respectively $\lambda_1$ and matrix $P$ obtained as follows.
We propose $\widehat{P}_n:= (\widehat{v}_{1,n}, \widehat{v}_{2,n})$ as estimator of $P=(v_1, v_2)$, the orthonormal basis of eigenvectors of matrix $A$, with:
$$   
\begin{cases}
\widehat{v}_{1,n}:=  \cos(\widehat{\theta}_{n}) v^{\star} - \sin(\widehat{\theta}_{n}) v^{\star \star}\\
\widehat{v}_{2,n}: =  \sin(\widehat{\theta}_{n}) v^{\star} + \cos(\widehat{\theta}_{n}) v^{\star \star} 
\end{cases}
$$
By Theorem \ref{convergence presque-sure}, $\widehat{P}_n$ is a consistent estimator of $P$.\\
Now for $\widehat{\lambda}_{1, n}$, first, we apply Theorem \ref{convergence sure} to the particular function $f \equiv {\bf 1}$ and then, we use the result of Proposition \ref{rice}.

We deduce that if $\widetilde{\Lambda}(\lambda):= \begin{pmatrix} 1 & 0 \\
0 & \lambda
\end{pmatrix}$, one has 
\begin{align}
\label{function Fi}
J_{1}^{(n)}(u) \xrightarrow[n\to+\infty]{a.s.} \lambda_1\, p_{Z(0)}(u)\, \mathbb{E}[\normp[2]{\widetilde{\Lambda}(\lambda)\, P^{t}.\nabla Z(0)}]= \lambda_1\, \Phi(u, \lambda, P),
\end{align}
where $\Phi$ is a continuous function of its arguments.
\\
Once again by applying Theorem \ref{convergence presque-sure} one gets a consistent estimator for $\lambda_1$ by taking:
$$
\widehat{\lambda}_{1, n}:= \frac{J_{1}^{(n)}(u)}{\Phi(u, \widehat{\lambda}_{n}, \widehat{P}_n)},
$$
and finally a consistent estimator of $\mathit{\Sigma^{\star}}(u)$ is given by
$$
 \mathit{\widehat{\Sigma_n^{\star}}}(u):= \ds \frac{1}{\widehat{\lambda}_{1, n}^2}\, \mathit{\Sigma_{\star}}(u, \widehat{\lambda}_{n}, \widehat{P}_n).
$$
The matrix $ \mathit{\widehat{\Sigma_n^{\star}}}(u)$ is computable and can be factorized as:
$ \mathit{\widehat{\Sigma_n^{\star}}}(u)= R_n \mathit{\Gamma}_{n}^{\star} R_n^t$, where $R_n$ is an unitary matrix and $\mathit{\Gamma}_{n}^{\star} $ is a diagonal matrix.\\
Remark that $\mathit{\Gamma}_{n}^{\star}$ is invertible since $\mathit{\Sigma}^{\star}(u)$ is non-degenerate (see Remark \ref{non degenerate}).
Thus let:
$$
D_n(u):= (C^{-1}(\widehat{\lambda}_{n}, \widehat{\theta}_{n}))^t Q R_n (\mathit{\Gamma}_{n}^{\star})^{-\frac{1}{2}} R_n^t.
$$
Theorem \ref{convergence lambda et theta} implies Corollary \ref{convergence lambda et theta Bis}.
\begin{corollary}
\label{convergence lambda et theta Bis}
For $0 < \lambda < 1$ and $-\frac{\pi}{2} < \theta < \frac{\pi}{2}$, one has
\begin{align*}
2n \left({
\widehat{\lambda}_{n} - \lambda; \,\widehat{\theta}_{n} - \theta
}\right)  . \, D_n(u)
\xrightarrow[n \to +\infty]{Law} \mathcal{N}(0; I_2).
\end{align*}
\end{corollary}
\begin{proofarg}{Proof of Corollary \ref{convergence lambda et theta Bis}}
First, let us prove that $\mathit{\Sigma}^{\star}(u)= \mathit{\Sigma}^{\star}(u, \lambda_1, \lambda_2, P)= \ds \frac{1}{\lambda_1^2}\, \mathit{\Sigma}_{\star}(u, \lambda, P)$, where $\mathit{\Sigma}_{\star}$ is a computable and a continuous matrix as function of $(\lambda, P)$.

On the one hand, this fact provides from the form of the ratio of coefficients 
${a_{f_i^{\star}}({\bk}, u)}/{a_{\bf 1}({\bf 0}, u)}$, $i=1,\,2,\,3$, 
these latter coefficients being defined in Lemmas \ref{calculus of c} and \ref{calculus of a} of Appendix \ref{ann:proof}.
Indeed these ratios only depend on the ratio $\lambda={\lambda_2}/{\lambda_1}$ and $P$ and they are computable as function of $\lambda$ and $P$.
They do not depend on $\mu$.
On the other hand by Lemma \ref{mehler} given in Section \ref{chaos section} one can see that the term
$R({\bk}, {\bm})$ only depends on the covariance function of $U$ through the following form.
By defining $W(v)$ as the $3$-dimensional vector defined as
$W(v) :=  \begin{pmatrix}
\frac{\ds P^t}{\ds \sqrt{\mu}}. \nabla Z(v),
\frac{\ds Z(v)}{\ds \sqrt{r_z(0)}}
\end{pmatrix}^t
$, 
one has
\begin{multline*}
R({\bk}, {\bm})= \int_{\reels^2} \mathbb{E}[\widetilde{H}_{{\bk}}(U(0)) \widetilde{H}_{{\bm}}(U(v))] \ud v\\
= \frac{1}{\lambda_1^2 \lambda} \int_{\reels^2} \mathbb{E}[\widetilde{H}_{{\bk}}(W(0)) \widetilde{H}_{{\bm}}(W(v))] \ud v 
= \frac{1}{\lambda_1^2} G(\lambda, P, {\bk}, {\bm}),
\end{multline*}
and since $r_z$ is supposed to be known, $G$ is computable as function of $\lambda$ and $P$.
Thus $\mathit{\Sigma}^{\star}(u)= \frac{1}{\lambda_1^2}\, \mathit{\Sigma}_{\star}(u, \lambda, P)$, and $\mathit{\Sigma}_{\star}$ is computable.

It remains to prove that $\mathit{\Sigma}_{\star}$ is a continuous matrix as function of $\lambda$ and $P$.
In this aim, let us compute $\mathit{\Sigma}_{f_1, f_1}(u)$, similar arguments would be raised for $\mathit{\Sigma}_{f_2, f_2}(u)$ and for $\mathit{\Sigma}_{f_1, f_2}(u)$, where we recall that
$
f_i= {\ds \frac{f_i^{\star}}{a_{\bf 1}({\bf 0}, u)}- \frac{a_{f_{i}^{\star}}({\bf 0}, u)}{a_{\bf 1}^2({\bf 0}, u)} {\bf 1}}
$, $i=1, 2$ (see Remark \ref{two dimensional}).

Applying Proposition \ref{variance asymptotique xi} to $f:=f_1$, the one order Rice formula for $\mathbb{E}[J_{f_1}^{(n)}(u)]$ (see Proposition \ref{rice}) and the second order Rice formula for $\mathbb{E}[J_{f_1}^{(n)}(u)]^2$ (see second part of Lemma \ref{continuite Y}), it is easy to see that another expression for $\mathit{\Sigma}_{f_1, f_1}(u)$ is
\begin{multline*}
\mathit{\Sigma}_{f_1, f_1}(u)= \int\limits_{\reels^2}  \left\{{ {\mathbb{E}[f_1(\nu_{X}(t))f_1(\nu_{X}(0)) \normp[2]{\nabla X(t)} \normp[2]{\nabla X(0)}\mid X(0)=X(t)=u] } }\right.\\
  \left.{ 
\times  p_{X(0), X(t)}(u, u) -(\mathbb{E}[f_1(\nu_X(0)) \normp[2]{\nabla X(0)}])^2 p^2_{X(0)}(u) 
 }\right\}   \ud t.
\end{multline*}
Denoting by $B(\lambda, P):= P\, \widetilde{\Lambda}(\lambda) P^t$ and letting $v:=A .t$, one gets
\begin{multline*}
\mathit{\Sigma}_{f_1, f_1}(u)= \frac{1}{a_{\bf 1}^2({\bf 0}, u)} \frac{1}{\lambda} \int_{\reels^2} \left\{{ 
 E\left[{g_1\Big(\frac{B(\lambda, P).\nabla Z(v)}{\normp[2]{B(\lambda, P).\nabla Z(v)}}\Big) 
g_1\Big(\frac{B(\lambda, P).\nabla Z(0)}{\normp[2]{B(\lambda, P).\nabla Z(0)}}\Big) 
  }\right.
 }\right.\\
{  \left.{  \normp[2]{B(\lambda, P).\nabla Z(v)} \normp[2]{B(\lambda, P).\nabla Z(0)}/ Z(0)=Z(v)=u }\right] p_{Z(0), Z(v)}(u, u) }\\
\left.{ -\left({\mathbb{E}[g_1\Big(\frac{B(\lambda, P).\nabla Z(0)}{\normp[2]{B(\lambda, P). \nabla Z(0)}}\Big) \normp[2]{B(\lambda, P).\nabla Z(0)}]}\right)^2 p^2_{Z(0)}(u) 
 }\right\}
 \ud v,
\end{multline*}
where $g_1:=f_1^{\star}-\frac{a_{f_1^{\star}}({\bf 0}, u)}{a_{\bf 1}({\bf 0}, u)} :=f_1^{\star} -b(u, \lambda, P)$, while $b$ is a continuous function of $\lambda$ and $P$.\\
We conclude the argument noticing that $a_{\bf 1}({\bf 0}, u)= \lambda_1\, c(u, \lambda)$ with $c$ a strictly positive continuous function of $\lambda$.

Now to prove Corollary \ref{convergence lambda et theta Bis}, we have Lemma \ref{convergence racine carree} proved in Appendix \ref{ann:proof}.
\begin{lemma}
\label{convergence racine carree}
Let $\Sigma:=R \mathit{\Gamma} R^t$ a definite positive matrix such that $R$ is an unitary matrix, while $\mathit{\Gamma}$ is a diagonal one.
Let also $(\Sigma_n)_n$ be an approximation of matrix $\Sigma$, i.e.
$\lim\limits_{n\to+\infty}\Sigma_n=\Sigma$, such that $\Sigma_n:=R_n \mathit{\Gamma}_n R_n^t$ with $R_n$ an unitary matrix and $\mathit{\Gamma}_n$ a diagonal one.
Consider $B_n$ a square root of $\Sigma_n$, that is $B_n:=R_n \mathit{\Gamma}_n^{\frac{1}{2}} R_n^t$.
Then $\lim\limits_{n\to+\infty}B_n=B:=R \mathit{\Gamma}^{\frac{1}{2}} R^t$.
\end{lemma}
Turning back to the proof of Corollary \ref{convergence lambda et theta Bis}, we apply Lemma \ref{convergence racine carree} to $\Sigma:= \mathit{\Sigma}^{\star}(u)=R \mathit{\Gamma}^{\star} R^t$, $\Sigma_n:= \mathit{\widehat{\Sigma_n^{\star}}}(u)= R_n \mathit{\Gamma}_n^{\star} R_n^t$.
Since $\widehat{\mathit{\Sigma}}_n^{\star}(u)$ is a consistent estimator of $\mathit{\Sigma}^{\star}(u)$, we deduce that $\lim\limits_{n\to+\infty}\Sigma_n=\Sigma$.
Using Lemma \ref{convergence racine carree} and Theorem \ref{convergence presque-sure}, we obtain that almost surely 
$$
\lim_{n\to+\infty} D_n(u)= D(u):= (C^{-1}(\lambda, \theta))^t Q R (\mathit{\Gamma}^{\star})^{-\frac{1}{2}} R^t.
$$
Using Theorem \ref{convergence lambda et theta} and the fact that $D^t(u) \mathit{\Sigma}_{\lambda, \theta}(u) D(u)=I_2$, one finally proved corollary.
\end{proofarg}

We deal in previous section \ref{convergence pour les parametres d'affinite} with the convergence for parameters $\lambda$ and $\theta$ when they are such that $0 < \lambda < 1$ and $-\frac{\ds \pi}{\ds 2} < \theta < \frac{\pi}{2}$. In next section we will complete these convergence results by focusing on the case where $\lambda=1$ that will lead naturally to an isotropy test.
\section{Towards a test of isotropy}
\label{Complementary results}
\subsection{Complementary results for estimating the parameter $\lambda$}
\subsubsection{Almost sure convergence for $\widehat{\lambda}_{n}$}
We emphasize that convergence result in Theorem \ref{convergence presque-sure} is valid under the assumption that 
$0 < \lambda < 1$ and $-\frac{\ds \pi}{\ds 2} < \theta < \frac{\pi}{2}$.
However, we can better elaborate what is happening to $\widehat{\lambda}_{n}$  in the isotropic case, when $\lambda=1$ (and $-\frac{\pi}{2} < \theta \le \frac{\pi}{2}$) and also when $0 < \lambda <1$ and $ \theta = \frac{\pi}{2}$, via the following theorem.
\begin{theorem}
\label{convergence presque-sure Bis}
\mbox{For} $\lambda=1$ \mbox{and} $-\frac{\pi}{2} < \theta \le \frac{\pi}{2}$ \mbox{\,or for \,} $0 <\lambda <1$ \mbox{and} $\theta = \frac{\pi}{2}$, \mbox{one has}
\begin{equation*}
 \widehat{\lambda}_{n} \xrightarrow[n\to+\infty]{a.s.} \lambda.
\end{equation*}
\end{theorem}
\begin{proofarg}{Proof of Theorem \ref{convergence presque-sure Bis}}
First let us consider the isotropic case, that is the case where $\lambda=1$ and $\theta$ being any parameter belonging to $]-\frac{\pi}{2}, \frac{\pi}{2}]$.
By Remark \ref{equation 72} stated at the end of Proposition \ref{diffeomorphisme} proof, the estimator $\widehat{\lambda}_{n}$ has to verify the following equation
\begin{equation*}
X_n^2\, I^4(\widehat{\lambda}_{n}) \,(X_n^2+Y_n^2)-X_n^2\,I^2(\widehat{\lambda}_{n})\,(\widehat{\lambda}_{n}^2+1)+\widehat{\lambda}_{n}^2=0.
\end{equation*}
Thus
\begin{multline}
\label{equation72bis}
\left({X_n^2 -(\tfrac{2}{\pi})^2}\right) \left[{ I^4(\widehat{\lambda}_{n}) \,(X_n^2+Y_n^2) +(\tfrac{2}{\pi})^2 \,I^4(\widehat{\lambda}_{n}) -I^2(\widehat{\lambda}_{n})\,(\widehat{\lambda}_{n}^2+1) }\right]\\ 
+ (\frac{2}{\pi})^2\, I^4(\widehat{\lambda}_{n})\, Y_n^2 - (\frac{2}{\pi})^4 \left({ (\frac{\pi}{2})^2- I^2(\widehat{\lambda}_{n}) }\right)\left({I^2(\widehat{\lambda}_{n}) - (\frac{\pi}{2})^2 \,\widehat{\lambda}_{n}^2 }\right)=0.
\end{multline}
We define the function $h$ for $0 < \lambda \le 1$, by
\begin{align}
\label{functionh}
h(\lambda):= \left({\frac{\pi}{2}- I(\lambda) }\right)\left({I(\lambda) - \frac{\pi}{2} \,\lambda }\right).
\end{align}
By using convergence given in (\ref{convergence presque sure}), we establish that
$$
\left({X_n, Y_n }\right) \xrightarrow[n\to+\infty]{a.s.} \left({\ds \tfrac{2}{\pi}\,, \,0}\right),
$$
and since $0 < \widehat{\lambda}_{n} \le 1$ and $1 < I(\widehat{\lambda}_{n}) \le \frac{\pi}{2}$, we finally showed that 
$$
h(\widehat{\lambda}_{n}) \xrightarrow[n\to+\infty]{a.s.} 0.
$$
Since $h$ is a strictly decreasing continuous function on $\introf{0}{1}$ and $h(1)=0$, we obtain that 
$$
\widehat{\lambda}_{n} \xrightarrow[n\to+\infty]{a.s.} 1,
$$
that is the required convergence.\\
Let us look now at the case where $0 <\lambda <1$ and $\theta = \frac{\pi}{2}$.
\\
Convergence established in (\ref{convergence presque sure}) now gives that
\begin{align*}
\left({X_n, Y_n }\right) \xrightarrow[n\to+\infty]{a.s.} (X, Y)
= \left({\frac{\lambda}{I(\lambda)}\,, \,0}\right).
\end{align*}
In the same way as before and using Remark \ref{equation 72}, one gets the following equality.
\begin{multline}
\left({X_n^2 -X^2}\right) \left[{ I^4(\widehat{\lambda}_{n}) \,(X_n^2+Y_n^2) +X^2 \,I^4(\widehat{\lambda}_{n}) -I^2(\widehat{\lambda}_{n})\,(\widehat{\lambda}_{n}^2+1) }\right] +X^2\, I^4(\widehat{\lambda}_{n})\, Y_n^2  \\ 
- (\tfrac{1}{I(\lambda)})^4 \left({ I^2(\lambda)- \lambda^2\,I^2(\widehat{\lambda}_{n}) }\right)\left({\lambda\, I(\widehat{\lambda}_{n}) + I(\lambda)\,\widehat{\lambda}_{n} }\right) \left({\lambda\, I(\widehat{\lambda}_{n}) - I(\lambda)\,\widehat{\lambda}_{n} }\right)=0.
\label{equation72ter}
\end{multline}
For fixed $0 < \lambda < 1$, let us define function $f$ by
\begin{align}
\label{function f}
f(x):=\lambda I(x) - I(\lambda)\, x, \mbox{\, for \,} 0 \le x \le 1.
\end{align}
Previous almost sure convergence and the fact that $0 < \widehat{\lambda}_{n} \le 1$, $1 < I(\widehat{\lambda}_{n}) \le \frac{\pi}{2}$ and for $\lambda <1$, $I(\lambda)- \lambda\,I(\widehat{\lambda}_{n}) > I(\lambda) -\frac{\pi}{2} >0$, imply that
\begin{align*}
f(\widehat{\lambda}_{n}) \xrightarrow[n\to+\infty]{a.s.} 0.
\end{align*}
A straightforward calculation shows that function $f$ is a strictly decreasing continuous function on $\intrff{0}{1}$ such that $f(\lambda)=0$.
We deduce that
\begin{align*}
\widehat{\lambda}_{n} \xrightarrow[n\to+\infty]{a.s.} \lambda,
\end{align*}
that yields Theorem \ref{convergence presque-sure Bis}.
\end{proofarg}
We thus set up a first approach to detect if the process $X$ is isotropic or not.
\subsubsection{Convergence in law for $\widehat{\lambda}_{n}$}
In case where $0 < \lambda < 1$ and $-\frac{\ds \pi}{\ds 2} < \theta < \frac{\pi}{2}$, Theorem \ref{convergence lambda et theta} gives as corollary a convergence law result for $2n (\widehat{\lambda}_{n} - \lambda)$. Here we complete the statement of this convergence under the assumption that 
 $\lambda=1$ (and $-\frac{\pi}{2} < \theta \le \frac{\pi}{2}$) or when $0 < \lambda <1$ and $ \theta = \frac{\pi}{2}$ via the following theorem.
\begin{theorem} 
\label{raccord}
\begin{enumerate}
\item For $0 < \lambda <1$ and $\theta = \frac{\pi}{2}$, one has
\begin{align*}
2n \left({
\widehat{\lambda}_{n} - \lambda
}\right)
\xrightarrow[n \to +\infty]{Law} \mathcal{N}(0; \left({\mathit{\Sigma}_{\lambda, \frac{\pi}{2}}(u)}\right)_{11}),
\end{align*}
where $\left({\mathit{\Sigma}_{\lambda, \frac{\pi}{2}}(u)}\right)_{11}$ stands for the first row and first column element locaded in matrix $\mathit{\Sigma}_{\lambda, \frac{\pi}{2}}(u)$ defined in Theorem \ref{convergence lambda et theta}.
\item For $\lambda =1$ and $-\frac{\pi}{2} < \theta \le \frac{\pi}{2}$, one has
\begin{align*}
2n\,  (1-\widehat{\lambda}_{n})
 \xrightarrow[n \to +\infty]{Law} \sqrt{V},
\end{align*}
where the density $f_V(t)$ of the positive \rv
$V$ is given by:
\begin{align*}
f_V(t):= \frac{1}{\sigma_{11} \widetilde{\sigma}_{22}}\left({\ds \frac{1}{4\pi} \int_{0}^{2\pi} e^{\ds -t \frac{\ds (\cos(\theta)-a \sin(\theta))^2}{\ds 2\, \sigma_{11}^2\,(1+a^2)}}\, e^{ \ds -t \frac{\ds \sin^2(\theta)\, (1+a^2)}{\ds 2\, \widetilde{\sigma}_{22}^2}}\, \ud\theta}\right) \1_{\{t >0 \}},
\end{align*}
where the coefficients $a$ and $\widetilde{\sigma}_{22}$ are defined by:
\begin{align*}
a:= \frac{\sigma_{12}}{\sigma_{11}^2},\,  \widetilde{\sigma}_{22}^2:= \frac{\sigma_{11}^2\sigma_{22}^2-\sigma_{12}^2}{\sigma_{11}^2},
\end{align*}
while coefficients $\sigma_{11}$, $\sigma_{22}$ and $\sigma_{12}$ are
\begin{align*}
\sigma_{11}^2:= 4 \left({ \frac{\pi}{2}}\right)^2 \mathit{\Sigma}_{11}^{(\star \star)}(u),\, \sigma_{22}^2:=  \left({ \frac{\pi}{2}}\right)^2 \mathit{\Sigma}_{22}^{(\star \star)}(u),\, \sigma_{12}:= 2 \left({ \frac{\pi}{2}}\right)^2 \mathit{\Sigma}_{12}^{(\star \star)}(u),\, 
\end{align*}
with 
\begin{align*}
\mathit{\Sigma}^{(\star \star)}(u)= \left({ \mathit{\Sigma}_{ij}^{(\star \star)}(u)}\right)_{1\le i, j \le2}:= \mathit{\Sigma_Q^{\star}}(u, \tau, \tau, I_2) (u),
\end{align*}
where $\tau$ is the common value of the eigenvalues of matrix $A$ under the isotropic hypothesis.
\end{enumerate}
\end{theorem}
\begin{remark}
\label{computation}
If $\lambda \neq 1$, by part one of Theorem \ref{raccord}, we readily get that $2n\,  (1-\widehat{\lambda}_{n})$ converges in law to Gaussian \rv with positive infinite mean when $n$ tends to infinity.
Thus one gets a one more way to detect if the process is isotropic or not.
\end{remark}
\begin{remark}
\label{non degenerate matrix}
By Remark \ref{non degenerate}, the matrix $ \mathit{\Sigma}^{(\star \star)}(u)$ is non-degenerate, ensuring that the coefficient $\widetilde{\sigma}_{22}$ appearing in the expression of the density $f_V$ does not vanish.
\end{remark}
\begin{remark}
\label{density}
We suppose that the covariance function $r_z$ of $Z$ is \textbf{known}.
We can estimate the density $f_V$ by $\widehat{f}_V(t):= f_{V^{\prime}}(\widehat{\tau}_{n}^2t) \widehat{\tau}_{n}^2$, where $\widehat{\tau}_{n}$ is a consistent estimator of $\tau$, the common value of the eigenvalues of matrix $A$ under the isotropic case and where $V^{\prime}$ has the same law as that of $V$, where we replaced $\mathit{\Sigma}^{(\star \star)}(u)$ by $Q \,\mathit{\Sigma}_{\star}(u, 1, I_2) \,Q^t$, where $\mathit{\Sigma}_{\star}$ being defined in equality (\ref{Sigma etoile})  is computable.
By using convergence in (\ref{function Fi}) we can take for $\widehat{\tau}_{n}$,
\begin{align}
\label{tau n}
\widehat{\tau}_{n}:= 2\, J_{1}^{(n)}(u)\, \sqrt{\frac{r_z(0)}{\mu}} \exp(\frac{1}{2} \frac{u^2}{r_z(0)}).
\end{align}
\end{remark}
\begin{proofarg}{Proof of Theorem \ref{raccord}}
Let us look first at the case where $0 < \lambda < 1$ and $\theta= \frac{\pi}{2}$, then we will be treat the case where $\lambda = 1$ (and $-\frac{\pi}{2} < \theta \le \frac{\pi}{2}$).\\
In the first case, by decomposition obtained in (\ref{equation72ter}) in the proof of Theorem \ref{convergence presque-sure Bis}, we have
\begin{multline*}
2n\left({X_n -X}\right) \left({X_n + X}\right) \left[{ I^4(\widehat{\lambda}_{n}) \,(X_n^2+Y_n^2) +X^2 \,I^4(\widehat{\lambda}_{n}) -I^2(\widehat{\lambda}_{n})\,(\widehat{\lambda}_{n}^2+1) }\right] \\
+2n Y_n^2  I^4(\widehat{\lambda}_{n})X^2
- 2n  \left({\lambda\, I(\widehat{\lambda}_{n}) - I(\lambda)\,\widehat{\lambda}_{n} }\right) \times \\
\frac{1}{I^4(\lambda)}\left({I^2(\lambda)-\lambda^2\,I^2(\widehat{\lambda}_{n})}\right) \left({\lambda\, I(\widehat{\lambda}_{n}) + I(\lambda)\,\widehat{\lambda}_{n} }\right)=0,
\end{multline*}
where
\begin{align*}
\left({X_n, Y_n }\right) \xrightarrow[n\to+\infty]{a.s.} (X, Y)= \left({{\lambda}/{I(\lambda)}\,, \,0}\right).
\end{align*}
Now by using this almost sure convergence result and those given in (\ref{convergence F}) plus Theorem \ref{convergence presque-sure Bis}, one obtains the following probability equivalence 
\begin{align*}
2n\left({X_n -X}\right) 2\lambda \,I(\lambda) (\lambda^2-1) \equiv 2n  \left({\lambda\, I(\widehat{\lambda}_{n}) - I(\lambda)\,\widehat{\lambda}_{n} }\right) \frac{2\lambda}{I(\lambda)}
(1-\lambda^2),
\end{align*}
and since $0<\lambda <1$, one finally gets
\begin{align*}
 2n  \left({\lambda\, I(\widehat{\lambda}_{n}) - I(\lambda)\,\widehat{\lambda}_{n} }\right) \equiv
 -2n\left({X_n -X}\right) I^2(\lambda).
\end{align*}
Now refering to function $f$ defined by (\ref{function f}) in the proof of Theorem \ref{convergence presque-sure Bis}, let
\begin{align*}
f(x)=\lambda I(x) - I(\lambda)\, x, \mbox{\, for \,} 0 \le x \le 1.
\end{align*}
One has $f(\lambda)=0$ and $f'(\lambda)=\lambda I'(\lambda)-I(\lambda) <0$.\\
So with the first order Taylor expansion of  $f$ about $\lambda$ evaluated at $\widehat{\lambda}_{n}$, one gets
\begin{align*}
2n  \left({\lambda\, I(\widehat{\lambda}_{n}) - I(\lambda)\,\widehat{\lambda}_{n} }\right)&= 2n f(\widehat{\lambda}_{n})\\
 &\equiv  2nf'(\lambda) (\widehat{\lambda}_{n}-\lambda) \\
 &\equiv -2n\left({X_n -X}\right) I^2(\lambda),
\end{align*}
and then one proved that
$$
 2n (\widehat{\lambda}_{n}-\lambda) \equiv \frac{2n\left({X_n -X}\right) I^2(\lambda)}{I(\lambda)-\lambda I'(\lambda)},
$$
and by convergence given in (\ref{convergence F}), \begin{align*}
2n \left({
\widehat{\lambda}_{n} - \lambda
}\right)
\xrightarrow[n \to +\infty]{Law} \mathcal{N}(0; \mathit{\Sigma}),
\end{align*}
where $\mathit{\Sigma}:= \frac{\ds I^4(\lambda)}{\ds (I(\lambda)-\lambda I'(\lambda))^2}\,(\mathit{\Sigma_Q^{\star}}(u))_{11}$.

To end the first part of this proof we just have to check that
$(\mathit{\Sigma}_{\lambda, \frac{\pi}{2}}(u))_{11}= \mathit{\Sigma}$.
In this aim, remember that 
$
\mathit{\Sigma}_{\lambda, \frac{\pi}{2}}(u)=C(\lambda, \frac{\pi}{2}) \times \mathit{\Sigma_Q^{\star}}(u) \times C^{t}(\lambda, \frac{\pi}{2})
$
with
$$
C(\lambda, \frac{\pi}{2})= \frac{1}{J_{\overrightarrow{F}}(\lambda, \frac{\pi}{2})}
\begin{pmatrix}
\frac{\partial F_2}{\partial \theta}(\lambda, \frac{\pi}{2}) & -\frac{\partial F_1}{\partial \theta}(\lambda, \frac{\pi}{2}) \\
&\\
- & -
\end{pmatrix}.
$$
Using Lemma \ref{jacobien}, one obtains that 
$J_{\overrightarrow{F}}(\lambda, \frac{\pi}{2})= \ds \frac{(\lambda^2-1)}{\lambda\, I(\lambda)}\times \left({\frac{\lambda I'(\lambda)-I(\lambda)}{I^2(\lambda)} }\right)$ and
\begin{align*}
C(\lambda, \frac{\pi}{2})= \frac{1}{J_{\overrightarrow{F}}(\lambda, \frac{\pi}{2})}
\begin{pmatrix}
\frac{1-\lambda^2}{\lambda I(\lambda)} & 0 \\
&\\
- &\ds  -
\end{pmatrix}.
\end{align*} 
In this way, one proved that $(\mathit{\Sigma}_{\lambda, \frac{\pi}{2}}(u))_{11}= \mathit{\Sigma}$, yielding the first part
of theorem. Now we consider the second part taking $\lambda=1$ (and $-\frac{\pi}{2} < \theta \le \frac{\pi}{2}$).\\

The decomposition obtained in (\ref{equation72bis}) gives
\begin{multline}
\label{Z}
2n\left({X_n -\tfrac{2}{\pi}}\right) \left({X_n +\tfrac{2}{\pi}}\right) \times \\2n \left[{ I^4(\widehat{\lambda}_{n}) \,(X_n^2+Y_n^2) +(\tfrac{2}{\pi})^2 \,I^4(\widehat{\lambda}_{n}) -I^2(\widehat{\lambda}_{n})\,(\widehat{\lambda}_{n}^2+1) }\right] \\ 
+ \,(\tfrac{2}{\pi})^2\, I^4(\widehat{\lambda}_{n})\, (2n Y_n)^2\\ 
= (2n)^2 \left({ \tfrac{\pi}{2}- I(\widehat{\lambda}_{n}) }\right)\left({I(\widehat{\lambda}_{n}) - \tfrac{\pi}{2}  \,\widehat{\lambda}_{n} }\right)(\tfrac{2}{\pi})^4 \left({ \tfrac{\pi}{2}+ I(\widehat{\lambda}_{n}) }\right)
\left({I(\widehat{\lambda}_{n}) + \tfrac{\pi}{2}  \,\widehat{\lambda}_{n} }\right)
\end{multline}
where
\begin{align}
\label{convergence presque sure Xn Yn}
\left({X_n, Y_n }\right) \xrightarrow[n\to+\infty]{a.s.} \left({X, Y}\right)=\left({\tfrac{2}{\pi}\,, \,0}\right).
\end{align}
Note that we can not hold the same reasoning as in the first part of the proof.
Indeed, by the first part of Theorem \ref{convergence presque-sure Bis} and the latter almost sure convergence result, we deduce that
\begin{align*}
Z_n:=I^4(\widehat{\lambda}_{n}) \,(X_n^2+Y_n^2) +(\frac{2}{\pi})^2 \,I^4(\widehat{\lambda}_{n}) -I^2(\widehat{\lambda}_{n})\,(\widehat{\lambda}_{n}^2+1)  \xrightarrow[n\to+\infty]{a.s.} 0,
\end{align*}
thus we have to normalize the studied expression by $(2n)^2$.

However, if in a first time we do not normalize this expression by $(2n)^2$ but rather by $(2n)$, we get
\begin{align*}
2n\left({X_n -\frac{2}{\pi}}\right) \left({X_n +\frac{2}{\pi}}\right) Z_n + (\frac{2}{\pi})^2\, I^4(\widehat{\lambda}_{n})\, 2n Y_n^2 \\
=2n\, h(\widehat{\lambda}_{n})\,(\frac{2}{\pi})^4 \left({ \frac{\pi}{2}+ I(\widehat{\lambda}_{n}) }\right)
\left({I(\widehat{\lambda}_{n}) + \frac{\pi}{2}  \,\widehat{\lambda}_{n} }\right),
\end{align*}
remembering that the function $h$ (see (\ref{functionh})) has been defined in the proof of Theorem \ref{convergence presque-sure Bis} by 
\begin{align*}
h(\lambda)= \left({\tfrac{\pi}{2}- I(\lambda) }\right)\left({I(\lambda) - \tfrac{\pi}{2} \,\lambda }\right),
\end{align*}
for $0 < \lambda \le 1$.\\
Now by using the almost sure convergence to 0 for $Z_n$ and those given in (\ref{convergence F}) and (\ref{convergence presque sure Xn Yn}) and in the first part  of Theorem \ref{convergence presque-sure Bis}, we get
\begin{align*}
2n\, h(\widehat{\lambda}_{n}) \xrightarrow[n\to+\infty]{Pr} 0.
\end{align*}
Since $h(1)=h'(1)=0$, one obtains by using a second order Taylor expansion of $h$ about 1 evaluated at point $\widehat{\lambda}_{n}$,
\begin{align}
\label{bonus bis}
h(\widehat{\lambda}_{n})= \frac{\pi^2}{16} (\widehat{\lambda}_{n}-1)^2 +o((\widehat{\lambda}_{n}-1)^2).
\end{align}
Thus as a bonus we proved that in probability 
\begin{align}
\label{bonus}
2n\, (\widehat{\lambda}_{n}-1)^2 =o(1).
\end{align}
This last equality will help us to study the normalized expression $2n\, Z_n$ and show that in probability 
\begin{align*}
2n\, Z_n \equiv 2n\, (X_n-\tfrac{2}{\pi}) \frac{\pi^3}{4} .
\end{align*}
Indeed
\begin{align*}
2n\, Z_n=I^2(\widehat{\lambda}_{n})\left[{2n\, \left({X_n^2-(\tfrac{2}{\pi})^2}\right)\,I^2(\widehat{\lambda}_{n}) +2n\, I^2(\widehat{\lambda}_{n}) Y_n^2 }\right.\\
\left.{+2n\, \left({2\, I^2(\widehat{\lambda}_{n})\, (\tfrac{2}{\pi})^2 - \left({\widehat{\lambda}_{n}^2+1 }\right)}\right)
}\right].
\end{align*}
Using an order two Taylor expansion of the elliptic integral $I$ (see (\ref{integrale elliptique}) for definition) in $\widehat{\lambda}_{n}$ about $\lambda=1$, one gets
\begin{align*}
I(\widehat{\lambda}_{n})=\frac{\pi}{2}+\frac{\pi}{4} (\widehat{\lambda}_{n}-1)+O((\widehat{\lambda}_{n}-1)^2),
\end{align*}
in such a way that 
\begin{align*}
2n\, \left({2\, I^2(\widehat{\lambda}_{n})\, (\tfrac{2}{\pi})^2 - \left({\widehat{\lambda}_{n}^2+1}\right)}\right)= -n(\widehat{\lambda}_{n}-1)^2+O(2n\,(\widehat{\lambda}_{n}-1)^2).
\end{align*}
By equality (\ref{bonus}) this expression converges in probability  to zero.

Part one of Theorem \ref{convergence presque-sure Bis}, convergences in (\ref{convergence F}) and (\ref{convergence presque sure Xn Yn}) give the required result
\begin{align*}
2n\, Z_n \equiv 2n\, (X_n-\tfrac{2}{\pi}) \tfrac{\pi^3}{4} .
\end{align*}
Back to expression (\ref{Z}), using once again the convergence obtained in part one of Theorem \ref{convergence presque-sure Bis} and equality (\ref{bonus bis}), one finally proved that in probability
\begin{align*}
\pi^2\,\left({2n\left({X_n -\tfrac{2}{\pi}}\right)}\right)^2 + \left({2n\, Y_n}\right)^2 \left({\tfrac{\pi}{2}}\right)^2
\equiv (2n)^2\, h(\widehat{\lambda}_{n})\, \tfrac{16}{\pi^2} \equiv (2n\, (\widehat{\lambda}_{n}-1))^2.
\end{align*}
By using convergence given in (\ref{convergence F}), we get the required convergence.
This ends the proof of this theorem.
\end{proofarg}
\begin{proofarg}{Proof of Remark \ref{density}}
The density $f_V$ of the positive \rv
$V$ can be expressed as $f_V(t)=f_{V^{\prime}}(t \tau^2) \tau^{2}$, where $\tau$ is the common eigenvalue of matrix $A$ under the isotropic case and  $V^{\prime}$ is defined as $V$, substituting $Q \mathit{\Sigma}_{\star}(u, 1, I_2) Q^t$ to $\mathit{\Sigma}^{\star}(u, \tau, \tau, I_2)$.
The matrix $\mathit{\Sigma}_{\star}$ is computable and is defined by (\ref{Sigma etoile}).
So, it is enough to estimate $\tau$ by a consistent estimator, say $\widehat{\tau}_n$.
Using convergence given in (\ref{function Fi}) we propose a consistent estimator of $\lambda_1= \tau$ in the isotropic case by taking 
\begin{align*}
\widehat{\tau}_{n}:= \frac{J_{\bf 1}^{(n)}(u)}{p_{Z(0)}(u) \mathbb{E}[\normp[2]{\nabla Z(0)}}]= 2\,J_{\bf 1}^{(n)}(u)\, \sqrt{\frac{r_z(0)}{\mu}} \exp\Big(\frac{1}{2} \frac{u^2}{r_z(0)}\Big),
\end{align*}
that is computable since $\mu$ and $r_z(0)$ are supposed to be known.
\end{proofarg}

In the next section, statistical tests are proposed for the null hypothesis ``$X$ is \textit{isotropic}'' against  ``$X$ is \textit{affine}''.
\subsection{Testing the isotropy}
We test
\begin{align*}
H_0: \lambda =1\qquad
\mbox{against} \qquad
H_1: \lambda <1.
\end{align*}
We still obtain a way to detect the possible isotropy of the process via the following corollaries.
\begin{corollary}
\label{convergence statistique}
Under the hypothesis $H_0$, the following convergence holds:
\begin{align*}
\frac{J_{\overrightarrow{f^{\star}}}^{(n)}(u)}{J_{\bf 1}^{(n)}(u)} \xrightarrow[n\to+\infty]{a.s.} \frac{2}{\pi} v^{\star}
\end{align*}
\end{corollary}
\begin{corollary}
\label{base}
Under the hypothesis $H_0$, we have the following convergence:
\begin{align*}
T_{\overrightarrow{f^{\star}}}^{(n)}(u):=2n\left({
\frac{J_{\overrightarrow{f^{\star}}}^{(n)}(u)}{J_{\bf 1}^{(n)}(u)} - \frac{2}{\pi} v^{\star}
}\right) \xrightarrow[n \to +\infty]{Law} \mathcal{N}(0; \mathit{\Sigma}^{\star}(u, \tau, \tau, I_2)),
\end{align*}
where $\tau$ is the common value of the eigenvalues of matrix $A$.
\end{corollary}
\begin{remark}
\label{non degenerescence}
By Remark \ref{non degenerate} we can observe that the matrix $ \mathit{\Sigma}^{\star}(u, \tau, \tau, I_2)$ is non-degenerate.
 \end{remark}
 \begin{remark}
 \label{esperance sous H1}
Under the alternative hypothesis $H_1$, the test statistic $T_{\overrightarrow{f^{\star}}}^{(n)}(u)$ converges in law toward a Gaussian \rv, the mean of which tends to infinity.
Indeed by using equality (\ref{formule esperance}) one can easily show that 
$$
\left({\ds \frac{\mathbb{E}[J_{\overrightarrow{f^{\star}}}^{(1)}(u)]}{\mathbb{E}[J_{\bf 1}^{(1)}(u)]} = \frac{2}{\pi} v^{\star}}\right) \iff \left({\lambda = 1 }\right).
$$
 \end{remark}
 \begin{proofarg}{Proof of Corollaries \ref{convergence statistique} and \ref{base}}
Using the almost sure convergence given in Proposition \ref{decomposition base}, convergence in law of Proposition \ref{convergence en loi}, equality (\ref{formule esperance}) and taking $\lambda=1$ one obviously gets Corollaries \ref{convergence statistique} and \ref{base}.
\end{proofarg}
\begin{proofarg}{Proof of Remark \ref{esperance sous H1}}
Appling Proposition \ref{convergence en loi}, under $H_1$, the test statistic $T_{\overrightarrow{f^{\star}}}^{(n)}(u)$ converges in law to a Gaussian \rv
with asymptotically mean equivalent to $2n\left({
\ds \frac{\mathbb{E}[J_{\overrightarrow{f^{\star}}}^{(1)}(u)]}{\mathbb{E}[J_{\bf 1}^{(1)}(u)]} - \frac{2}{\pi} v^{\star}
}\right)$.\\
Since $\left({
\ds\frac{\mathbb{E}[J_{\overrightarrow{f^{\star}}}^{(1)}(u)]}{\mathbb{E}[J_{\bf 1}^{(1)}(u)]} - \frac{2}{\pi} v^{\star}
}\right)$ does not depend on $n$ and is equal to zero if and if $\lambda=1$, this argument ends the proof of the remark.
\end{proofarg}

In the case where covariance {\bf function $r_z$ is known}, Corollary \ref{base} suggests another test statistic.
More precisely, let $R$ an unitary matrix obtained by diagonalizing the computable matrix $\mathit{\Sigma}_{\star}(u, 1, I_2)$ defined in (\ref{Sigma etoile}).
That is $\mathit{\Sigma}_{\star}(u, 1, I_2) = R\, \mathit{\Gamma}_{\star}\, R^t$.\\
 We consider the following test statistic $S_{\overrightarrow{f^{\star}}}^{(n)}(u)$:
 \begin{align*}
 S_{\overrightarrow{f^{\star}}}^{(n)}(u):=  \widehat{\tau}_{n}\, \mathit{\Gamma}_{\star}^{-\frac{1}{2}}\, R^t . T_{\overrightarrow{f^{\star}}}^{(n)}(u),
 \end{align*}
 where $\widehat{\tau}_{n}$ is given by (\ref{tau n}).\\
 We can state the following theorem.
 \begin{theorem}
 \label{xi}
Under the hypothesis $H_0$, we have the following convergence:
\begin{align*}
\Xi_{\overrightarrow{f^{\star}}}^{(n)}(u):= (S_{\overrightarrow{f^{\star}}}^{(n)}(u))^{t} \,S_{\overrightarrow{f^{\star}}}^{(n)}(u) \xrightarrow[n \to +\infty]{Law} \chi_2^2.
\end{align*}
 \end{theorem}
 \begin{remark}
 \label{rejection}
 The rejection region  is then $\Xi_{\overrightarrow{f^{\star}}}^{(n)}(u) > \gamma$.
This critical region provides a consistent test for any positive constant $\gamma$, because $\Xi_{\overrightarrow{f^{\star}}}^{(n)}(u)$ is stochastically unbounded, for $n \to +\infty$, except under the null hypothesis.
In fact when $\lambda < 1$, $\frac{1}{(2n)^2} \, \Xi_{\overrightarrow{f^{\star}}}^{(n)}(u)$ converges in probability to $b >0$, and this implies that $\Xi_{\overrightarrow{f^{\star}}}^{(n)}(u)$ converges in probability to $+\infty$.
 \end{remark}
\begin{proofarg}{Proof of Theorem \ref{xi} and Remark \ref{rejection}}
As we have already pointed out in Remark \ref{non degenerescence}, the matrix $\mathit{\Sigma}^{\star}(u, \tau, \tau, I_2)$ is invertible and we can factorize $\mathit{\Sigma}_{\star}(u, 1, I_2)$ into $\mathit{\Sigma}_{\star}(u, 1, I_2) = R\, \mathit{\Gamma}_{\star}\, R^t$, where $R$ is an unitary matrix.
Let 
\begin{align*}
S_{\overrightarrow{f^{\star}}}^{(n)}(u)= \widehat{\tau}_n\, \mathit{\Gamma}_{\star}^{-\frac{1}{2}} R^t .T_{\overrightarrow{f^{\star}}}^{(n)}(u).
 \end{align*}
 Since $\widehat{\tau}_n$ is a consistent estimator of $\tau$ under hypothesis $H_0$, by Corollary \ref{base}, asymptotically this random vector is a standard Gaussian one and Theorem \ref{xi} ensues.\\
To achieve the proof of Remark \ref{rejection}, let us check that $\frac{1}{(2n)^2}\, \Xi_{\overrightarrow{f^{\star}}}^{(n)}(u)$ converges in probability to $b >0$ under the alternative hypothesis $H_1$ that is when $\lambda <1$.\\
 In this aim, let us note 
 \begin{align*}
 \widetilde{S}_{\overrightarrow{f^{\star}}}^{(n)}(u):= \widehat{\tau}_n\, \mathit{\Gamma}_{\star}^{-\frac{1}{2}}  R^t . \widetilde{T}_{\overrightarrow{f^{\star}}}^{(n)}(u),
 \end{align*}
 where 
 \begin{align*}
\widetilde{T}_{\overrightarrow{f^{\star}}}^{(n)}(u):=2n\left({
\frac{J_{\overrightarrow{f^{\star}}}^{(n)}(u)}{J_{\bf 1}^{(n)}(u)} - \ds\frac{\mathbb{E}[J_{\overrightarrow{f^{\star}}}^{(1)}(u)]}{\mathbb{E}[J_{\bf 1}^{(1)}(u)]}
}\right).
\end{align*}
We have the following decomposition
\begin{align*}
\frac{1}{(2n)^2}\, \Xi_{\overrightarrow{f^{\star}}}^{(n)}(u)&= \frac{1}{(2n)^2} \,(S_{\overrightarrow{f^{\star}}}^{(n)}(u))^{t} S_{\overrightarrow{f^{\star}}}^{(n)}(u)\\
&= \frac{1}{(2n)^2} \left({S_{\overrightarrow{f^{\star}}}^{(n)}(u)-\widetilde{S}_{\overrightarrow{f^{\star}}}^{(n)}(u)
}\right)^{t}\, \left({S_{\overrightarrow{f^{\star}}}^{(n)}(u)-\widetilde{S}_{\overrightarrow{f^{\star}}}^{(n)}(u) }\right)\\
&\quad+
\frac{2}{(2n)^2} \,(\widetilde{S}_{\overrightarrow{f^{\star}}}^{(n)}(u))^{t} \,\left({S_{\overrightarrow{f^{\star}}}^{(n)}(u)-\widetilde{S}_{\overrightarrow{f^{\star}}}^{(n)}(u) 
}\right)\\
&\quad+\frac{1}{(2n)^2} \,(\widetilde{S}_{\overrightarrow{f^{\star}}}^{(n)}(u))^{t} \,(\widetilde{S}_{\overrightarrow{f^{\star}}}^{(n)}(u))\\
& = (1)+(2)+(3).
\end{align*}
Let us look at the third term $(3)$.
\begin{align*}
(3):= \left({\frac{1}{2n}\, (\widetilde{T}_{\overrightarrow{f^{\star}}}^{(n)}(u))^{t}}\right). \left({\widehat{\tau}_n^2\, R \mathit{\Gamma}_{\star}^{-1}R^t}\right) .\left({\frac{1}{2n}\, \widetilde{T}_{\overrightarrow{f^{\star}}}^{(n)}(u)}\right) \xrightarrow[n \to +\infty]{Pr} 0,
\end{align*}
by using Proposition \ref{convergence en loi} and the fact that $\widehat{\tau}_n$ almost surely converges to $\ds \sqrt{\frac{2}{\pi \mu}} \mathbb{E}[\normp[2]{\nabla X(0)}]$ (see Theorem \ref{convergence sure} and Proposition \ref{rice}).\\
For the second term $(2)$.
\begin{align*}
(2):=  \left({\frac{1}{n}\, (\widetilde{T}_{\overrightarrow{f^{\star}}}^{(n)}(u))^{t}}\right). \left({\widehat{\tau}_n^2\, R\, \mathit{\Gamma}_{\star}^{-1}R^t}\right) .\left({
\ds \frac{\mathbb{E}[J_{\overrightarrow{f^{\star}}}^{(1)}(u)]}{\mathbb{E}[J_{\bf 1}^{(1)}(u)]} - \frac{2}{\pi} v^{\star}
}\right),
\end{align*}
and this term tends in probability toward zero as for the third term $(3)$.\\
Finally, the first term $(1)$ gives
\begin{align*}
(1):= \left({\widehat{\tau}_n\, 
\mathit{\Gamma}_{\star}^{-\frac{1}{2}} R^t .\left({
\ds \frac{\mathbb{E}[J_{\overrightarrow{f^{\star}}}^{(1)}(u)]}{\mathbb{E}[J_{\bf 1}^{(1)}(u)]}  - \frac{2}{\pi} v^{\star}
}\right)
}\right)^{t}\,\times \\ \left({\widehat{\tau}_n\, 
\mathit{\Gamma}_{\star}^{-\frac{1}{2}} R^t .\left({
\ds \frac{\mathbb{E}[J_{\overrightarrow{f^{\star}}}^{(1)}(u)]}{\mathbb{E}[J_{\bf 1}^{(1)}(u)]} - \frac{2}{\pi} v^{\star}
}\right)
}\right)
=b_n,
\end{align*}
and 
\begin{align*}
\lim_{n\to+\infty} b_n= b:= \frac{2}{\pi \mu} (\mathbb{E}[\normp[2]{\nabla X(0)}])^2 \left({
\mathit{\Gamma}_{\star}^{-\frac{1}{2}} R^t .\left({
\ds \frac{\mathbb{E}[J_{\overrightarrow{f^{\star}}}^{(1)}(u)]}{\mathbb{E}[J_{\bf 1}^{(1)}(u)]} - \frac{2}{\pi} v^{\star}
}\right)
}\right)^{t}\,\times \\ \left({
\mathit{\Gamma}_{\star}^{-\frac{1}{2}} R^t .\left({
\ds \frac{\mathbb{E}[J_{\overrightarrow{f^{\star}}}^{(1)}(u)]}{\mathbb{E}[J_{\bf 1}^{(1)}(u)]}  - \frac{2}{\pi} v^{\star}
}\right)
}\right)
 \ge 0.
\end{align*}
Since $\mathit{\Gamma}_{\star}^{-\frac{1}{2}} R^{t}$ is invertible, one has
\begin{align*}
\left({ b >0}\right) \iff \left({
\left({
\ds \frac{\mathbb{E}[J_{\overrightarrow{f^{\star}}}^{(1)}(u)]}{\mathbb{E}[J_{\bf 1}^{(1)}(u)]}  - \frac{2}{\pi} v^{\star}
}\right) \neq 0
}\right) \iff \left({ \lambda < 1 }\right).
\end{align*}
That ends the proof of remark.
\end{proofarg}
%
%
\appendix
\section{Appendix}
\label{ann:proof}
\begin{proofarg}{Proof of Lemma \ref{continuite Y}}
Let  $F: T \times (\reels^{2})^{\star} \to \reels$ be a bounded continuous function of its arguments. We want to prove that
$
\textbf{y} \mapsto  \mathbb{E}[I(\textbf{y})]^2
$
is a continuous function, where
for $y \in \reels$, we have noted $I(y):=\int_{\mathcal{C}(y)} F(t, \nabla X(t)) \1_{\textrm{D}^{\textrm{r}}_{\textrm{X}}}(t) \ud \sigma_1(t)$. Also we would like to show a second order Rice formula for the functionnal $I(y)$.\\
Due to the presence of the indicator, function $\textbf{t} \mapsto F(\textbf{t}, \nabla X(\textbf{t})) \1_{\textrm{D}^{\textrm{r}}_{\textrm{X}}}(\textbf{t})$ is not continuous on $T$, so we can not apply the two-order Rice formula (\ref{formule de Rice d'ordre deux}) given in Corollary \ref{gaussien stationnaire}.\\
Thus we consider the following approximation.
For $m \in \naturels^{\star}$, let 
\begin{align*}
I_{m}(y):=\int_{\mathcal{C}(y)} G_m(t, \nabla X(t)) \ud \sigma_1(t),
\end{align*}
where function $G_m: T \times \reels^{2} \to \reels$, is defined by
\begin{align*}
G_m(t, z):= \varphi\left({\frac{1}{m \normp[2]{z}}}\right) F(t, z) \1_{\{z \neq 0\}},
\end{align*}
where $\varphi$ has been defined in the proof of Theorem \ref{continuite longueur de courbe}.\\
As in the proof of last cited theorem we can prove that function $G_m$ is a bounded continuous function on $T \times \reels^{2}$.\\
The two order Rice formula ensues by applying Corollary \ref{gaussien stationnaire} to the functional $I_m(y)$. That is for all $y \in \reels$, one has
\begin{align*}
\lefteqn{\mathbb{E}[I_m(y)]^2=}\\
& \int_{T \times T}
\mathbb{E} \left[G_m(t_1, \nabla X(t_1)) G_m(t_2, \nabla X(t_2)) \normp[2]{\nabla
X(t_1)} \normp[2]{\nabla
X(t_2)}|X(t_1)=X(t_2)=y\right]\\
&\qquad\qquad\qquad \, \times {p}_{X(t_1), X(t_2)}(y, y)\, \ud t_1\, \ud t_2 < +\infty.
\end{align*}
We observe that the following convergence holds. For all $t \in T$ and $z \in \reels$, \\$G_m(t, z) \xrightarrow[m \to +\infty]{} F(t, z) \1_{\{z \neq 0\}}$.
Then we point out that Remark \ref{finitude courbe de Niveau} following Corollary \ref{gaussien stationnaire} ensures that\\
$\mathbb{E}[ \sigma_1(\mathcal{C}(y))]^2 = \int_{T \times T}
	\mathbb{E} \left[\normp[2]{\nabla
X(t_1)} \normp[2]{\nabla
X(t_2)}|X(t_1)=X(t_2)=y\right] \,\times {p}_{X(t_1), X(t_2)}(y, y)\ud t_1\ud t_2 < +\infty$.\\
Thus by using the fact that function $G_m$ is bounded, an application of the Lebesgue convergence theorem induces that
\begin{align*}
\lefteqn{\mathbb{E}[I(y)]^2=
 \int_{T \times T} \mathbb{E} \left[{F(t_1, \nabla X(t_1)) \1_{\textrm{D}^{\textrm{r}}_{\textrm{X}}}(t_1) F(t_2, \nabla X(t_2)) \1_{\textrm{D}^{\textrm{r}}_{\textrm{X}}}(t_2) \times}\right.}\\
&\left.{ \normp[2]{\nabla
X(t_1)} \normp[2]{\nabla
X(t_2)}|X(t_1)=X(t_2)=y}\right]
\,  \times {p}_{X(t_1), X(t_2)}(y, y)\, \ud t_1\, \ud t_2 < +\infty.
\end{align*}
We have therefore established the second order Rice formula for $I(y)$.\\
We are going to prove that function $
\textbf{y} \mapsto  \mathbb{E}[I(\textbf{y})]^2
$
is continuous.
This would achieve the proof of Lemma \ref{continuite Y}.\\
In this aim, we use the decomposition given in Wschebor \cite[p{.} 60]{MR0871689}, that is, for $\tau:=t-s$, where $t, s \in T$
\begin{eqnarray*}
\nabla X(0)&=& \xi + (X(0) \alpha + X(\tau)\, \beta)\\
\nabla X(\tau)&=& \xi^{\star} - (X(\tau) \alpha + X(0)\, \beta),
\end{eqnarray*}
where $\xi$ and $\xi^{\star}$ are centred Gaussian vectors taking values in $\reels^2$, with joint Gaussian distribution, each of them independent of $(X(0), X(\tau))$, such that, if $\nabla r_x$ stands for the Jacobian of the $X$ covariance function,
\begin{align*}
\alpha&:= \frac{r_x(\tau) \nabla r_x(\tau)}{r^2_x(0)-r^2_x(\tau)}; \qquad \beta:= -\frac{r_x(0) \nabla r_x(\tau)}{r^2_x(0)-r^2_x(\tau)}\\
\Var(\xi)&=\Var(\xi^{\star})=-\nabla^{2}r_x(0) -\frac{r_x(0) }{r^2_x(0)-r^2_x(\tau)}\, \nabla r_x(\tau) \left({\nabla r_x(\tau)}\right)^{t}\\
\Cov(\xi, \xi^{\star})&=-\nabla^{2}r_x(\tau) -\frac{r_x(\tau) }{r^2_x(0)-r^2_x(\tau)}\, \nabla r_x(\tau) \left({\nabla r_x(\tau)}\right)^{t}.
\end{align*}
Thus with these notations one has
\begin{align*}
\mathbb{E}[I(y)]^2= \int_T \ud s\, \int_{T} L(y, s, t-s) \ud t,
\end{align*}
where for $\tau:=t-s$
\begin{align*}
\lefteqn{L(y, s, \tau):= \exp\left(-\ds \frac{y^2}{r_x(0)+r_x(\tau)}\right) \frac{1}{2\pi} \,\frac{1}{\left({r^2_x(0)-r^2_x(\tau)}\right)^{\frac{1}{2}}}}\\ 
& \quad\times \mathbb{E}\left[{
F(s, \xi - \frac{y}{r_x(0)+r_x(\tau)} \nabla r_x(\tau))} \1_{\{\xi - \frac{y}{r_x(0)+r_x(\tau)} \nabla r_x(\tau) \neq 0\}} \normp[2]{\xi - \frac{y}{r_x(0)+r_x(\tau)} \nabla r_x(\tau)}\right.\\
& \quad\times \left.{F(\tau +s, \xi^{\star} + \frac{y}{r_x(0)+r_x(\tau)} \nabla r_x(\tau))} \1_{\{\xi^{\star} + \frac{y}{r_x(0)+r_x(\tau)} \nabla r_x(\tau) \neq 0\}} \normp[2]{\xi^{\star}+ \frac{y}{r_x(0)+r_x(\tau)} \nabla r_x(\tau)} 
\right].
\end{align*}
Using that $F$ is a bounded continuous function on $T \times (\reels^{2})^{\star}$ one gets that
\begin{itemize}
\item $\textbf{z} \mapsto F(s, \textbf{z}) \1_{\{\textbf{z} \neq 0\}} \normp[2]{\textbf{z}}$ is continuous for all $s \in T$.
\end{itemize}
Furthermore since $r_x^2(0)-r^2_x(t) \neq 0$ for all $t \neq 0$, this implies in force that $r_x(0)+r_x(t) \neq 0$ for all $t \neq 0$, thus by using that function $F$ is bounded one obtains that
\begin{itemize}
\item $\textbf{y} \mapsto L(\textbf{y}, s, t-s)$ is continuous for almost $(s, t) \in T \times T$.
\end{itemize}
Now let us enunciate a lemma for which a proof can be founded for example in Berzin and Wschebor \cite[Lemma 1 (b)]{MR1222362}.
\begin{lemma}
\label{covariance r}
$\exists \,M >0$, $\forall\, t, s \in T, \left({\normp[2]{t-s} \le M \Longrightarrow r^2_x(0)-r^2_x(t-s) \ge M \normp[2]{t-s}^2 }\right)$
\end{lemma}
Now let us choose $B>0$ small enough such that for all $\tau:=s-t$, with  $s, t  \in T$ such that $\normp[2]{\tau} \le B$, one has the following inequalities,
\begin{itemize}
\item (1) $r_x(0)+r_x(\tau) \ge \frac{3}{2}\, r_x(0)$
\item (2) $r^2_x(0)-r^2_x(\tau) \ge B \normp[2]{\tau}^2$
\item (3) $\normp[2]{\nabla r_x(\tau)} \le {\bf C} \normp[2]{\tau}$
\end{itemize}
Note that these three inequalities are always possible to implement.
Indeed the first inequality just translates the fact that covariance function $r_x$ is continuous in zero, the second one comes from Lemma \ref{covariance r} and the third one from a first order Taylor expansion of $\nabla r_x(\cdot)$ about zero.

These inequalities allow showing the following bound on $\abs{L(y, s, \tau)}$.
\begin{itemize}
\item $\exists\, D >0$ such that for all $(s, t) \in T\times T$ and $\tau :=s-t$,\\
$\abs{L(y, s, \tau)} \le D (1+y^2)\left\{
\1_{\{\normp[2]{\tau} \ge B \}} + \frac{1}{\normp[2]{\tau}}\, \1_{\{\normp[2]{\tau} \le B \} }
\right\}$
\end{itemize}
Before explaining the last inequality, let us remark that since
$\int_{T\times T} \frac{\ud s\ud t}{\normp[2]{s-t}} < +\infty$, an application of Lebesgue's convergence theorem yields that function $
\textbf{y} \mapsto  \mathbb{E}[I(\textbf{y})]^2
$
is continuous.\\
Indeed, two cases occur: $\normp[2]{\tau} \ge B$ and $\normp[2]{\tau} \le B$.
If $\normp[2]{\tau} \ge B$, we consider the compact set $K:= \{s-t, (s, t) \in \overline{T} \times \overline{T} \mbox{ \, with \,} \normp[2]{s-t} \ge B \}$ where we have noted $\overline{T}$ for the closure of $T$. Since $r_x$ is continuous on $K$  and $r_x^2(0)-r^2_x(v) \neq 0$ for $v \in K$, thus there exists $E>0$ such for all $v \in K$, $r_x^2(0)-r^2_x(v) \ge E$ and $r_x(0)+r_x(v) \ge E$, so that $\abs{L(y, s, \tau)} \le D (1+y^2)$.
If $\normp[2]{\tau} \le B$, inequalities (2) and (3) give that 
$\mathbb{E}[\normp[2]{\xi}] \le {\bf C}$ and $\mathbb{E}[\normp[2]{\xi^{\star}}] \le {\bf C}$ and inequalities (1) and (2) give that $\abs{L(y, s, \tau)} \le D (1+y^2)
 \frac{1}{\normp[2]{\tau}}$.
 
 That completes the proof of Lemma \ref{continuite Y}.
\end{proofarg}

\begin{remark}
The approach chosen to prove Lemma \ref{continuite Y} does not work for showing that
\begin{align*}
({\bf{y_1}}, {\bf{y_2}}) \mapsto \mathbb{E}[J_{f}({\bf{y_1}}) J_{f}({\bf{y_2}})],
\end{align*}
is continuous.
In fact if we try to reproduce the previous proof, we realize that we have to check that
for all $y_1, y_2 \in \reels$,
\begin{multline*}
H(y_1, y_2):=\\
 \int_{T \times T}
\mathbb{E}[\normp[2]{\nabla
X(t_1)} \normp[2]{\nabla
X(t_2)}|X(t_1)=y_1, X(t_2)=y_2] 
\times {p}_{X(t_1), X(t_2)}(y_1, y_2) \ud t_1 \ud t_2 < +\infty.
\end{multline*}
As before, performing a regression, the expectation appearing in the expression of $H(y_1, y_2)$ is 
\begin{multline*}
\mathbb{E}\left[{
 \normp[2]{\xi + \frac{r_x(\tau) \nabla r_x(\tau)}{r^2_x(0)-r_x^2(\tau)} y_1 
-  \frac{r_x(0) \nabla r_x(\tau)}{r^2_x(0)-r_x^2(\tau)} y_2}\times 
\normp[2]{\xi^{\star} + \frac{r_x(0) \nabla r_x(\tau)}{r^2_x(0)-r_x^2(\tau)} y_1 
-  \frac{r_x(\tau) \nabla r_x(\tau)}{r^2_x(0)-r_x^2(\tau)} y_2}
}\right],
\end{multline*}
and the principal difficulty results in bounding $\ds \frac{\nabla r_x(\tau)}{r^2_x(0)-r^2_x(\tau)}$ when $\normp[2]{\tau} \le B$.
That is why we developped a little more sophisticated approach to prove Theorem \ref{continuite longueur de courbe}.\\
\end{remark}

\begin{proofarg}{Proof of Lemma \ref{mehler}}
Let $X=(X_i)_{i=1,\,2,\,3}$ and $Y=(Y_j)_{j=1,\,2,\,3}$ be two centred standard Gaussian vectors in $\reels^3$ such that for $1 \le i, j \le 3$, $\mathbb{E}[X_iY_j]=\rho_{ij}$.\\
Let ${\bk}=(k_1, k_2, k_3)$ and ${\bm}=(m_1, m_2, m_3)$ be two vectors of $\naturels^3$.
We want to give an explicit formula for $\mathbb{E}[\widetilde{H}_{{\bk}}(X) \widetilde{H}_{{\bm}}(Y)]$.\\
For ease of notations let us set $Y_j=X_{j+3}$ for $j=1,\,2,\,3$.

As in Aza{\"\i}s and Wschebor \cite[p{.} 269]{MR2478201}, a straightforward calculation on Gaussian characteristic  functions gives
\begin{multline}
\label{charasteristic}
\mathbb{E}[\prod_{i=1}^{6} \exp(t_i X_i - \tfrac{1}{2} t_i^2)]= \exp\Big(\sum_{i=1}^{3} \sum_{j=1}^{3} \rho_{ij} \,t_i t_{j+3}\Big)\\
= \sum_{r=0}^{+\infty} \frac{1}{r!}\, \Big({\sum_{i=1}^{3} \sum_{j=1}^{3} \rho_{ij} \,t_i t_{j+3}}\Big)^r
=  \sum_{r=0}^{+\infty} \sum_{\sum_{i, j}d_{ij}=r} \prod_{i,j=1}^{3} \frac{1}{d_{ij}!} \,\left({\rho_{ij} \,t_i t_{j+3} }\right)^{d_{ij}}\\
= \sum_{r=0}^{+\infty} \sum_{\sum_{i, j} d_{ij}=r} \left({\prod_{i,j=1}^{3} \frac{\rho_{ij}^{d_{ij}}}{d_{ij}!}}\right) \,
\left({\prod_{i=1}^{3} t_i^{\sum_{j=1}^3 d_{ij}}}\right) \left({\prod_{j=1}^{3} t_{j+3}^{\sum_{i=1}^3 d_{ij}}}\right)
\end{multline}
Now, by definition
$$
\exp(tx-\tfrac{1}{2}t^2)= \sum_{q=0}^{+\infty} t^q\, \frac{H_q(x)}{q!},
$$
thus 
\begin{multline*}
\mathbb{E}[\prod_{i=1}^{6} \exp(t_i X_i - \tfrac{1}{2} t_i^2)]\\
=\sum_{k_1=0}^{+\infty}\sum_{k_2=0}^{+\infty}\sum_{k_3=0}^{+\infty}\sum_{m_1=0}^{+\infty}\sum_{m_2=0}^{+\infty}\sum_{m_3=0}^{+\infty}
\frac{t_1^{k_1}}{k_1!}\,\frac{t_2^{k_2}}{k_2!}\,\frac{t_3^{k_3}}{k_3!}\,\frac{t_4^{m_1}}{m_1!}\,\frac{t_5^{m_2}}{m_2!}\,\frac{t_6^{m_6}}{m_6!}\,\mathbb{E}[\widetilde{H}_{{\bk}}(X) \widetilde{H}_{{\bm}}(Y)].
\end{multline*}
Identifying this last equality with (\ref{charasteristic}), it follows that $\mathbb{E}[\widetilde{H}_{{\bk}}(X) \widetilde{H}_{{\bm}}(Y)]=0$ if $\abs{{\bk}}\neq \abs{{\bm}}$.
In the case where $\abs{{\bk}} = \abs{{\bm}}$, one gets
\begin{align*}
\mathbb{E}[\widetilde{H}_{{\bk}}(X) \widetilde{H}_{{\bm}}(Y)]=\sum_{
\substack{
d_{ij} \ge 0 \\ \sum_{j}d_{ij}=k_i \\ \sum_{i}d_{ij}=m_j
}
} {\bk}!  {\bm}! \prod_{1 \le i, j \le 3} \frac{\rho_{ij}^{d_{ij}}}{d_{ij}!},
\end{align*}
yielding the lemma.
\end{proofarg}

\begin{proofarg}{Proof of Lemma \ref{lemme cramer}}
Let $f: S^{1} \to \reels$ be a positive continuous and bounded function. For $n \in \naturels^{\star}$, let $T_n:=\introo{0}{n}\times \introo{0}{n}$. One has
\begin{align*}
(1)&=\int_{0}^{n-1} \int_{0}^{n-1} \int\limits_{\mathcal{C}_{[t\,, \,t+1[\times [s\, ,\,s +1[}(u)} f(\nu_{X}(x)) \ud \sigma_{1}(x)\, \ud t \ud s \\
&\hfill- \int_{0}^{1}\int_{0}^{1} \int_{\mathcal{C}_{\introo{0}{t} \times \introo{0}{s}}(u)} f(\nu_{X}(x)) \ud \sigma_{1}(x)\, \ud t \ud s \\
& =
 \int_{0}^{n-1} \int_{0}^{n-1} \left({H(t+1, s+1)-H(t, s+1)-H(t+1, s)+H(t, s) }\right) \ud t\ud s\\
&\quad - \int_{0}^{1} \int_{0}^{1} H(t, s) \ud t\ud s,
 \end{align*}
 where we recall that we have noted
 \begin{equation*}
 H(t, s):=\int_{\mathcal{C}_{\introo{0}{t} \times \introo{0}{s}}(u)} f(\nu_{X}(x)) \ud \sigma_{1}(x).
 \end{equation*}
 Thus by making change of variable in the first integral, one obtains
 \begin{align*}
(1)=\int_{n-1}^{n} \int_{n-1}^{n}  H(t, s) \ud t \ud s - \int_{0}^{1} \int_{n-1}^{n}  H(t, s) \ud t \ud s -\int_{n-1}^{n} \int_{0}^{1}  H(t, s) \ud t \ud s.
 \end{align*}
 Since $f$ is a positive function then $H(t, s) \ge 0$ and we get the following upper bound
  \begin{multline*}
(1) \le \int_{n-1}^{n} \int_{n-1}^{n}  H(t, s) \ud t \ud s
\\ \le  \int_{n-1}^{n} \int_{n-1}^{n} \int_{\mathcal{C}_{\introo{0}{n} \times \introo{0}{n}}(u)} f(\nu_{X}(x)) \ud \sigma_{1}(x)\, \ud t \ud s = \int_{\mathcal{C}_{\introo{0}{n} \times \introo{0}{n}}(u)} f(\nu_{X}(x)) \ud \sigma_{1}(x).
 \end{multline*}
 The first inequality of lemma is then achieved.
To obtain the second inequality, arguing as previously we obtain the following lower bound
 \begin{multline*}
 (2)=\int_{0}^{n+1} \int_{0}^{n+1} \int_{\textrm{C}_{[t-1\,,\,t[\times[s-1\,,\,s[}(u)} f(\nu_{X}(x)) \ud \sigma_{1}(x)\, \ud t \ud s 
 =\int_{n}^{n+1} \int_{n}^{n+1}  H(t, s) \ud t \ud s\\ -\int_{-1}^{0} \int_{n}^{n+1}  H(t, s) \ud t \ud s - \int_{n}^{n+1}\int_{-1}^{0}   H(t, s) \ud t \ud s + \int_{-1}^{0}\int_{-1}^{0} H(t, s) \ud t \ud s\\
 \ge    \int_{n}^{n+1} \int_{n}^{n+1} \int_{\textrm{C}_{\introo{0}{t} \times \introo{0}{s}}(u)} f(\nu_{X}(x)) \ud \sigma_{1}(x)\, \ud t \ud s 
 \ge \int_{\mathcal{C}_{\introo{0}{n} \times \introo{0}{n}}(u)} f(\nu_{X}(x)) \ud \sigma_{1}(x),
\end{multline*}
this yields Lemma \ref{lemme cramer}.
\end{proofarg}
\begin{proofarg}{Proof of Lemma \ref{convergence esperance}}
Since $Z$ is an isotropic process, its density $p_{\nabla Z(0)}(v)$ depends only on $\normp[2]{v}$.
We denote it by $g(\normp[2]{v})$.\\
Thus if $\overrightarrow{f}: S^1 \to \reels^2$ is a continuous and bounded function, 
\begin{align*}
\mathbb{E}[\overrightarrow{f}(\nu_{X}(0)) \normp[2]{\nabla X(0)}]
	&= \mathbb{E}[\overrightarrow{f}\left(\tfrac{A. \nabla Z(0)}{\normp[2]{A . \nabla Z(0)}}\right) \normp[2]{A . \nabla Z(0)}]\\
	&= \int_{\reels^2} \overrightarrow{f}\left(\tfrac{A.v}{\normp[2]{A.v}}\right)\, \normp[2]{A.v}\, g(\normp[2]{v})\, \ud v.
\end{align*}
Letting $v:=r \alpha$, $\alpha \in S^1$, one gets
\begin{align*}
\lefteqn{\mathbb{E}[\overrightarrow{f}(\nu_{X}(0)) \normp[2]{\nabla X(0)}]}\\
&=\sigma_{1}(S^1)\, \left({\int_{0}^{\infty} g(r) r^2\, \ud r}\right)\, \left({\int_{S^1} \overrightarrow{f}\left(\tfrac{\ds A. \alpha}{\ds \normp[2]{A . \alpha}}\right)\, \normp[2]{A. \alpha}\, \ud\alpha}\right)\\
 &=\mathbb{E}[\normp[2]{\nabla Z(0)}]\, \left({\int_{S^1} \overrightarrow{f}\left(\tfrac{\ds A. \alpha}{\ds \normp[2]{A.\alpha}}\right)\, \normp[2]{A. \alpha}\, \ud\alpha}\right).
\end{align*}
In a similar way taking $f\equiv  {\bf 1}$, one gets
\begin{align*}
\mathbb{E}[\normp[2]{\nabla X(0)}]=\mathbb{E}[\normp[2]{\nabla Z(0)}]\, \left({\int_{S^1}  \normp[2]{A. \alpha}\, \ud\alpha}\right).
\end{align*}
Finally by using Proposition \ref{rice}, one obtains
\begin{align*}
\frac{\mathbb{E}[J_{\overrightarrow{f}}^{(n)}(u)]}{\mathbb{E}[J_{\bf 1}^{(n)}(u)]}&= \frac{\int_{S^1} \overrightarrow{f}\left(\tfrac{\ds A. \alpha}{\ds \normp[2]{A. \alpha}}\right)\, \normp[2]{A. \alpha}\, \ud\alpha}{\int_{S^1}  \normp[2]{A . \alpha}\, \ud\alpha}.
\end{align*}
Now choosing $\overrightarrow{f}:=\overrightarrow{f^{\star}}$ defined by (\ref{fstar}) and since $A$ is a self-adjoint matrix, one has
\begin{multline*}
\int_{S^1} \overrightarrow{f^{\star}}\left(\frac{A. \alpha}{\normp[2]{A.\alpha}}\right)\, \normp[2]{A.\alpha}\, \ud\alpha=2\, \int_{S^1} A.\alpha\, \1_{\{\ds \langle A.\alpha, v^{\star}\rangle \ge 0 \}}\, \ud\alpha\\
=2\, \int_{S^1} A. \alpha\, \1_{\{\ds \langle \alpha, A.v^{\star}\rangle \ge 0 \}}\, \ud\alpha
=2\, A. \int_{\ds S^1\cap \{\langle \alpha, w^{\star}\rangle \ge 0 \}} \alpha\, \ud\alpha,
\end{multline*}
where we have noted $w^{\star}:= \frac{\ds A.v^{\star}}{\ds \normp[2]{A.v^{\star}}}$.\\
If $R$ denotes the change of basis matrix from the canonical orthonormal basis ${\bf e}:=(\vec{i}, \vec{j})$ to an orthonormal basis choosen as ${\bf w}:=(w_{1}^{\star}, w^{\star})$ and making the change of variable $\beta:=R^{t} .\alpha$ in the last integral, one obtains
\begin{align*}
\lefteqn{\left({\int_{\ds S^1\cap \{\langle \alpha, w^{\star}\rangle_{{\bf e}} \ge 0 \}} \alpha\, \ud\alpha}\right)_{{\bf e}}}\\&=
\left({\int_{\ds S^1\cap \{\langle \alpha, (0,1)\rangle_{{\bf w}} \ge 0 \}} \alpha\, \ud\alpha}\right)_{{\bf w}}=
\left({0, \left({\int_{\ds S^{1} \cap \{\alpha_{2} \ge 0 \}} \alpha_{2} \ud\alpha
}\right)
}\right)_{{\bf w}}\\
&= \left({\int_{\ds S^{1} \cap \{\alpha_{2} \ge 0 \}} \alpha_{2} \ud\alpha }\right) w^{\star} = \left({\int_{0}^{\pi} \sin(x) \frac{\ud x}{2\pi} }\right)w^{\star}= \frac{1}{\pi} \,w^{\star}.
\end{align*}
Thus we have proved that
\begin{align*}
\int_{S^1} \overrightarrow{f^{\star}}\left(\tfrac{\ds A.\alpha}{\ds \normp[2]{A.\alpha}}\right)\, \normp[2]{A. \alpha}\, \ud\alpha = \frac{2}{\pi} A.w^{\star}=\frac{2}{\pi} \,\frac{A^2.v^{\star}}{\normp[2]{A.v^{\star}}},
\end{align*}
that yields Lemma \ref{convergence esperance}.
\end{proofarg}

\noindent 
\begin{lemma}
\label{jacobien}
The Jacobian $J_{\overrightarrow{F}}$ of transformation $\overrightarrow{F}$ (see (\ref{function F})) is given by 
\begin{multline}
\label{jacobien F}
J_{\overrightarrow{F}}(\lambda, \theta)= \frac{\lambda (1-\lambda^2)}{I^{3}(\lambda)(\cos^2(\theta)+\lambda^2 \sin^2(\theta))}\times \\
\left[{
\sin^{2}(\theta) \int_{0}^{1} \frac{\sqrt{1-u^2}}{\sqrt{1-(1-\lambda^2)u^2}}\, \ud u}\right.\\
 + \left.{
\cos^{2}(\theta) \int_{0}^{1} \frac{u^2}{\sqrt{1-(1-\lambda^2)u^2}}\, \frac{\ud u}{\sqrt{1-u^2}}
}\right] \neq 0,\, \,
\end{multline}
for $0< \lambda <1$.
\end{lemma}
\begin{proofarg}{Proof of Lemma \ref{jacobien}}
For $0 < \lambda \le 1$ and $-\frac{\pi}{2} < \theta \le \frac{\pi}{2}$, one has
$$
\begin{cases}
\dfrac{\partial F_1}{\partial \lambda}(\lambda, \theta)=  \ds \frac{1}{I^2(\lambda)}\, \frac{\left[{\lambda \sin^2(\theta)I(\lambda)-(\cos^2(\theta)+\lambda^2 \sin^2(\theta))I'(\lambda) }\right]}{(\cos^{2}(\theta)+ \lambda^2 \sin^2(\theta))^{\frac{1}{2}}}\\
\dfrac{\partial F_1}{\partial \theta}(\lambda, \theta)=  \ds \frac{1}{I(\lambda)}\, \frac{\left({\lambda^2-1}\right) \sin(\theta) \cos(\theta)}{(\cos^{2}(\theta)+ \lambda^2 \sin^2(\theta))^{\frac{1}{2}}},
\end{cases}
$$
and
$$
\begin{cases}
\dfrac{\partial F_2}{\partial \lambda}(\lambda, \theta)=  
	\ds \frac{\sin(\theta)\cos(\theta)}{I^2(\lambda)}
	\frac{\left({\small {\begin{array}{l}
				(\lambda I(\lambda)-(\lambda^2-1)I'(\lambda))\\
				\times (\cos^2(\theta)+\lambda^2 \sin^2(\theta))\\
				+\lambda I(\lambda)
			\end{array}
		}}\right)}
	{(\cos^{2}(\theta)+ \lambda^2 \sin^2(\theta))^{\frac{3}{2}}}\\
\dfrac{\partial F_2}{\partial \theta}(\lambda, \theta)=  \ds \frac{\left({\lambda^2-1}\right)}{I(\lambda)}\, \frac{\left[{\cos^4(\theta)-\lambda^2 \sin^4(\theta)}\right]}{(\cos^{2}(\theta) + \lambda^2 \sin^2(\theta))^{\frac{3}{2}}}.
\end{cases}
$$
From straightforward calculations we deduce that the Jacobian $J_{\overrightarrow{F}}$ of the transformation $\overrightarrow{F}$ can be written as
$$
J_{\overrightarrow{F}}(\lambda, \theta)= \ds \frac{\left({1-\lambda^2}\right)}{I^3(\lambda)}\,\frac{\left[{\lambda \sin^2(\theta)I(\lambda)+(\cos^2(\theta)-\lambda^2 \sin^2(\theta))I'(\lambda) }\right]}{(\cos^{2}(\theta)+ \lambda^2 \sin^2(\theta))} 
$$
Now
$I(\lambda)= \int_{0}^{\frac{\pi}{2}} (\cos^{2}(x)+ \lambda^2 \sin^2(x))^{\frac{1}{2}} \ud x$, and the change of variable:
$\sin(x)=u$ gives \\$I(\lambda)=\ds \int_{0}^{1} \frac{(1-(1-\lambda^2)\,u^2)^{\frac{1}{2}}}{(1-u^2)^{\frac{1}{2}}} \ud u$
and $I'(\lambda)=\ds \lambda \,\int_{0}^{1} \frac{u^2}{(1-(1-\lambda^2)\,u^2)^{\frac{1}{2}}} \frac{\ud u}{\sqrt{1-u^2}}$.\\
This yields 
\begin{multline*}
\lambda \sin^2(\theta)I(\lambda)+(\cos^2(\theta)-\lambda^2 \sin^2(\theta))I'(\lambda) =\\
\lambda\, \left[{
\sin^{2}(\theta) \int_{0}^{1} \frac{\sqrt{1-u^2}}{\sqrt{1-(1-\lambda^2)u^2}}\ud u}\right.\\
 + \left.{
\cos^{2}(\theta) \int_{0}^{1} \frac{u^2}{\sqrt{1-(1-\lambda^2)u^2}}\, \frac{\ud u}{\sqrt{1-u^2}}
}\right],
\end{multline*}
and lemma ensues from the fact that $\left({J_{\overrightarrow{F}}(\lambda, \theta)=0}\right) \iff \left({\lambda=1}\right)$.
\end{proofarg}

\begin{proofarg}{Proof of Lemma \ref{bornes}}
First let us see that for $n \in \naturels^{\star}$, a.s.
$X_n^2 +Y_n^2 <1$.
\begin{multline*}
X_n^2+Y_n^2= \frac{\normp[2]{\int_{\mathcal{C}_{n}(u)} \overrightarrow{f^{\star}}(\nu_{X}(t)) \ud\sigma_1(t)}^2}{\sigma_1^2(\mathcal{C}_{n}(u))} \le \frac{\int_{\mathcal{C}_{n}(u)} \normp[2]{\overrightarrow{f^{\star}}(\nu_{X}(t))}^2\ud \sigma_{1}(t)}{\sigma_1(\mathcal{C}_{n}(u))}=1,
\end{multline*}
last equality coming from the fact that $\overrightarrow{f^{\star}}$ takes its values into $S^1$.
But the strict inequality is true otherwise, we would have equality in H\"{o}lder inequality that remains impossible.

Now let us see that, a.s.
$X_n >0$.
\begin{align*}
X_n=\Big\langle \ds \frac{J_{\overrightarrow{f^{\star}}}^{(n)}(u)}{J_{\bf 1}^{(n)}(u)} , v^{\star} \Big\rangle = \frac{1}{\sigma_1(\mathcal{C}_{n}(u))}  \int_{\mathcal{C}_{n}(u)} \abs{ \langle \nu_{X}(t) , v^{\star} \rangle } \ud \sigma_1(t) \ge 0,
\end{align*}
and in a manner similar to that previously used it is proved that, a.s.
$X_n>0$.
This yields lemma.
\end{proofarg}

\begin{lemma}
\label{calculus of c}
The coefficients  $a_{f_{1}^{\star}}$ are given by: for $m, \ell \in \naturels$
\begin{align*}
\begin{pmatrix}
a_{f_1^{\star}}(2m, 2\ell)& a_{f_1^{\star}}(2m+1, 2\ell+1) & a_{f_1^{\star}}(2m, 2\ell+1)& a_{f_1^{\star}}(2m+1, 2\ell)\\
a_{f_2^{\star}}(2m, 2\ell)& a_{f_2^{\star}}(2m+1, 2\ell+1) & a_{f_2^{\star}}(2m, 2\ell+1)& a_{f_2^{\star}}(2m+1, 2\ell
\end{pmatrix}
= \\
P\,\times
\begin{pmatrix}
A(\lambda_1, \lambda_2, \omega_1^{\star}, \omega_2^{\star}, m, \ell, \mu) & B(\lambda_1, \lambda_2, \omega_1^{\star}, \omega_2^{\star}, m, \ell, \mu) & 0 & 0\\
A(\lambda_2, \lambda_1, \omega_2^{\star}, \omega_1^{\star}, \ell, m, \mu) & B(\lambda_2, \lambda_1, \omega_2^{\star}, \omega_1^{\star}, \ell, m, \mu) & 0 & 0
\end{pmatrix},
\end{align*}
where $\omega^{\star}=\begin{pmatrix} \omega_{1}^{\star} \,, \omega_{2}^{\star} \end{pmatrix} := v^{\star}.P$ and
\begin{multline*}
A(\lambda_1, \lambda_2, \omega_1^{\star}, \omega_2^{\star}, m, \ell, \mu):=\\
 \sqrt{\tfrac{2\mu}{\pi}} \lambda_{1}  \sum_{p=0}^{m} \frac{(-2)^{p-m}}{(2p)!(m-p)!}
\sum_{k=0}^{p} \frac{2^{k} p!}{(p-k)!}
\left({
\frac{\lambda_{2}\omega_{2}^{\star}}{\lambda_{1} \omega_{1}^{\star}}
}\right)^{2(p-k)} \times
\\
\sum_{n=0}^{\ell} \frac{(-2)^{n-\ell}}{(2n)!(\ell -n)!}
\left({
\frac{\lambda_{1} \omega_{1}^{\star}}{\left[{(\lambda_{1} \omega_{1}^{\star})^{2}+ (\lambda_{2} \omega^{\star}_{2})^{2}}\right]^{\frac{1}{2}}}
}\right)^{2(n+p-k)+1} \times\\ \frac{(2(n+p-k))!}{2^{n+p-k}(n+p-k)!},
\end{multline*}
and 
\begin{multline*}
B(\lambda_1, \lambda_2, \omega_1^{\star}, \omega_2^{\star}, m, \ell, \mu):=\\
 \sqrt{\frac{2\mu}{\pi}} \lambda_{1}  \sum_{p=0}^{\ell} \frac{(-2)^{p-\ell}}{(2p+1)!(\ell-p)!}
\sum_{k=0}^{p} \frac{2^{k} p!}{(p-k)!}
\left({
\frac{\lambda_{1}\omega_{1}^{\star}}{\lambda_{2} \omega_{2}^{\star}}
}\right)^{2(p-k)} \times
\\
\sum_{n=0}^{m} \frac{(-2)^{n-m}}{(2n+1)!(m -n)!}
\left({
\frac{\lambda_{2} \omega_{2}^{\star}}{\left[{(\lambda_{1} \omega_{1}^{\star})^{2}+ (\lambda_{2} \omega_{2}^{\star})^{2}}\right]^{\frac{1}{2}}}
}\right)^{2(p-k+n+1)+1} \times \\
\frac{(2(p-k+n+1))!}{2^{p-k+n+1}(p-k+n+1)!}.
\end{multline*}
\end{lemma}
\begin{proofarg}{Proof of Lemma \ref{calculus of c}}
Recall that for $f: S^{1} \to \reels$ a continuous and bounded functions and for all $u$ fixed level in $\reels$ the coefficients $a_f(\cdot, u)$ are defined by $a_f({\bk}, u)= a_f(k_1, k_2) \,a(k_3, u)
$, for ${\bk}=(k_1, k_2, k_3) \in \naturels^3$, 
with
$$
a(k_3, u)=\frac{1}{k_3!} H_{k_3}\bigg(\frac{u}{\sqrt{r_z(0)}}\bigg)\phi\bigg(\frac{u}{\sqrt{r_z(0)}}\bigg) \frac{1}{\sqrt{r_z(0)}}.
$$
So for giving an expression of the coefficients $a_{f_{i}^{\star}}(\cdot, u)$ for $i=1, 2$ and of $a_{\bf 1}(\cdot, u)$ just give one of $a_{f_{i}^{\star}}(k_1, k_2)$ and of that of $a_{\bf 1}(k_1, k_2)$.\\
First, let us compute $a_{f_1^{\star}}(k_1, k_2)$, for $k_1, k_2 \in \naturels$. Let $P:=(a_{i, j})_{1\le i, j \le 2}$.
By definition of coefficient $a_{f_1^{\star}}(k_1, k_2)$, we have
\begin{multline*}
a_{f_{1}^{\star}}(k_1, k_2) = \frac{\sqrt{\mu}}{k_1!k_2!} \int_{\reels^2} (a_{11}\,\lambda_1\,y_1+a_{12}\,\lambda_2\,y_2) \,\times \\
 \Big(\1_{\big\{\big\langle 
{\scriptsize \begin{pmatrix}
\lambda_1\,y_1\\
\lambda_2\,y_2
\end{pmatrix}},\, (\omega^{\star})^t\big\rangle \ge 0\big\}}
-\1_{\big\{\big\langle
{\scriptsize \begin{pmatrix}
\lambda_1\,y_1\\
\lambda_2\,y_2
\end{pmatrix}},\, (\omega^{\star})^t\big\rangle < 0\big\}}\Big)\,\times\\
H_{k_1}(y_1) \,\phi(y_1) \,H_{k_2}(y_2) \,\phi(y_2) \ud y_1 \ud y_2.
\end{multline*}
By using the fact that $H_n$ is even when $n$ is even and odd if not (see (\ref{hermite form}) and (\ref{H2l+1})), the coefficients $k_1$ and $k_2$ ought to be of the same parity, otherwise the coefficients  $a_{f_1^{\star}}(k_1, k_2)$ would be null.
Therefore let us suppose first that $k_1=2m$ and $k_2= 2\ell$, $m, \ell \in \naturels$.
In that way and using that polynomial $H_{2n}$ is even, one has
\begin{align*}
a_{f_1^{\star}}(2m, 2\ell) = \tfrac{2\sqrt{\mu}}{(2m)!(2\ell)!} \int_{\reels^2} (a_{11}\,\lambda_1\,y_1+a_{12}\,\lambda_2\,y_2) \,\times \\
 \1_{\big\{\big\langle 
{\scriptsize \begin{pmatrix}
\lambda_1\,y_1\\
\lambda_2\,y_2
\end{pmatrix}},\, (\omega^{\star})^t\big\rangle \ge 0\big\}} \times
H_{2m}(y_1) \,\phi(y_1) \,H_{2\ell}(y_2) \,\phi(y_2)\, \ud y_1\, \ud y_2.
\end{align*}
Let us compute 
\begin{align*}
E:=  \int_{\reels^2} y_1 \,\times 
  \1_{\big\{\big\langle 
{\scriptsize \begin{pmatrix}
\lambda_1\,y_1\\
\lambda_2\,y_2
\end{pmatrix}},\, (\omega^{\star})^t\big\rangle \ge 0\big\}}\times
H_{2m}(y_1) \,\phi(y_1) \,H_{2\ell}(y_2) \,\phi(y_2) \,\ud y_1\, \ud y_2.
\end{align*}
(similar computations would be done for the second integral on $y_2$, the arguments $\omega^{\star}_2, \omega^{\star}_1, \lambda_2, \lambda_1, \ell, m$, in this order, playing the role of $\omega^{\star}_1, \omega^{\star}_2, \lambda_1, \lambda_2, m, \ell$ in last integral).\\
At this step of the proof Lemma \ref{lemme Fp} is required.
\begin{lemma}
\label{lemme Fp}
Let $p \in \naturels$, $a \in \reels$ and $F_p(a):=\int_{a}^{+\infty} z^{2p+1} \phi(z) \ud z$.
Then
\begin{align*}
F_p(a)= \phi(a)\, \sum_{k=0}^{p} \frac{2^kp!}{(p-k)!} \,a^{2(p-k)}.
\end{align*}
\end{lemma}
\begin{proofarg}{Proof of Lemma \ref{lemme Fp}}
Integration by parts and straightforward calculations give the lemma.
\end{proofarg}
Three cases occur: $\omega^{\star}_1 >0$, $\omega^{\star}_1<0$ and $\omega^{\star}_1=0$.\\
Let us consider the first one.
If $\omega^{\star}_1 >0$, $E$ can be written as
\begin{align*}
E=  \int_{\reels} H_{2\ell}(y_2) \,\phi(y_2) \left[{\int_{-\frac{\lambda_2\, \omega^{\star}_2}{\lambda_1\, \omega^{\star}_1}y_2}^{+\infty}
y_1\,
H_{2m}(y_1) \,\phi(y_1)\ud y_1}\right]\ud y_2.
\end{align*}
For real $x$ and $m \in \naturels$, the polynomial form of $H_{2m}(x)$ is
\begin{align}
\label{hermite form}
H_{2m}(x)= (2m)! \sum_{p=0}^{m} \frac{(-2)^{p-m}}{(2p)!\, (m-p)!} \,x^{2p}.
\end{align}
Using Lemma \ref{lemme Fp}, one gets
\begin{align*}
E= (2m)! \sum_{p=0}^{m} \frac{(-2)^{p-m}}{(2p)!\, (m-p)!} \sum_{k=0}^p \frac{2^k p!}{(p-k)!}
\left({
\frac{\lambda_2 \omega^{\star}_2}
{\lambda_1 \omega^{\star}_1}
}\right)^{2(p-k)} \times\\
\sum_{n=0}^{\ell} (2\ell)! \frac{(-2)^{n-\ell}}{(2n)!\, (\ell-n)!} \,G_{n+p-k}\left({
\frac{\lambda_2 \omega^{\star}_2}
{\lambda_1 \omega^{\star}_1}
}\right),
\end{align*}
where for $q \in \naturels$ and $x \in \reels$, we defined 
\begin{align*}
G_q(x):=\int_{-\infty}^{+\infty} y^{2q}\, \phi(y)\, \phi(xy)\, dy=\frac{1}{\sqrt{2\pi}} \,\frac{1}{(1+x^2)^{q+\frac{1}{2}}} \frac{(2q)!}{2^q q!}.
\end{align*}
If $\omega^{\star}_1 <0$, using that for $m \in \naturels$, polynomial $H_{2m}$ is even, one obtains that  
\begin{align*}
E&=  \int_{\reels} H_{2\ell}(y_2) \phi(y_2) \left[{
\int_{-\infty}^{
-\frac{\lambda_2 \omega^{\star}_2}{\lambda_1 \omega^{\star}_1}
y_2}
y_1
H_{2m}(y_1) \phi(y_1) dy_1}\right] dy_2\\
&=-\int_{\reels} H_{2\ell}(y_2) \phi(y_2) \left[{\int_{-\frac{\lambda_2 \omega^{\star}_2}{\lambda_1 \omega^{\star}_1}y_2}^{+\infty}
y_1
H_{2m}(y_1) \phi(y_1) dy_1}\right] dy_2,
\end{align*}
and in a same way, when $\omega^{\star}_1 =0$,
\begin{align*}
E= \left({
\int_{\reels} H_{2\ell}(y_2)\, \phi(y_2)\, \1_{\{y_2 \omega_2^{\star} \ge 0\}} \,\ud y_2
}\right) \left({
\int_{\reels} y_1\, H_{2m}(y_1)\, \phi(y_1)\, \ud y_1
}\right) =0.
\end{align*}
Finally, knowing that $\abs{\omega_1^{\star}} \times\, \sign(\omega_1^{\star})=\omega_{1}^{\star}$, one gets the expression of coefficient $a_{f_1^{\star}}(2m, 2\ell)$.\\
For coefficient $a_{f_1^{\star}}(2m+1, 2\ell+1)$, similar arguments would be developed using the previous way and the polynomial form of $H_{2\ell+1}(x)$, that is
\begin{align}
\label{H2l+1}
H_{2\ell +1}(x)= (2\ell+1)! \sum_{p=0}^{\ell} \frac{(-2)^{p-\ell}}{(2p+1)!\, (\ell-p)!} \,x^{2p+1}.
\end{align}
To conclude the proof of Lemma \ref{calculus of c}, just remark that $a_{f_2^{\star}}(2m+1, 2\ell)=a_{f_2^{\star}}(2m, 2\ell+1)=0$ and that $a_{f_2^{\star}}(2m, 2\ell)$ (resp.
$a_{f_2^{\star}}(2m+1, 2\ell+1)$) would be computed replacing $a_{11}$ by $a_{21}$ and $a_{12}$ by $a_{22}$ in the expression of $a_{f_1^{\star}}(2m, 2\ell)$ (resp.
$a_{f_1^{\star}}(2m+1, 2\ell+1)$).
\end{proofarg}

In order to give the coefficients $a_{\bf 1}$ we introduce the following functions.\\
The $\beta$ function is defined by:
 \begin{align*}
 \beta(x; y):= \frac{\mathit{\Gamma}(x) \mathit{\Gamma}(y)}{\mathit{\Gamma}(x+y)},\quad x, y>0,
 \end{align*}
 while the $\mathit{\Gamma}$ function is defined by:
 \begin{align}
 \label{fonction Gamma}
 \mathit{\Gamma}(x):=\int_{0}^{\infty} e^{-t} t^{x-1} \ud t.
 \end{align}
\begin{lemma}
\label{calculus of a}
The coefficients  $a_{\bf 1}$ are given by: for $m, \ell \in \naturels$
\begin{align*}
\begin{pmatrix}
a_{\bf 1}(2m, 2\ell)& a_{\bf 1}(2m+1, 2\ell+1) & a_{\bf 1}(2m, 2\ell+1)& a_{\bf 1}(2m+1, 2\ell)
\end{pmatrix}
= \\
\begin{pmatrix}
c(\lambda_1, \lambda_2, m, \ell, \mu) & 0 & 0 & 0
\end{pmatrix},
\end{align*}
where
\begin{multline*}
c(\lambda_1, \lambda_2, m, \ell, \mu):=\\
\sqrt{2\pi \mu}\tfrac{(-2)^{-(m+\ell)}}{m! \ell!} \sum_{p=0}^{\ell} \sum_{q=0}^{m}
{{\ell}\choose{p}}
{{m}\choose{q}}
  (-1)^{p+q} \lambda_{2} \left({
 \tfrac{\lambda_2}{\lambda_1}
 }\right)^{2q+1} \times \\
 \sum_{n=0}^{+\infty} 
 \frac{
 {{q+n}\choose{q}}
  {{2q+2n}\choose{q+n}}
}
 {{{2q}\choose{q}}
 }
 \left({
\tfrac{1}{2}
}\right)^{2n}
\frac{1}{\beta(p+q+n+1; \frac{1}{2})}
\left({
1- \left({
\tfrac{\lambda_2}{\lambda_1}
}\right)^2
}\right)^{n}.
\end{multline*}
\end{lemma}
\begin{remark}
Note that in the case where $\lambda =1$, that is when the process is isotropic our result contains that of the authors expressed in Kratz and Le\'on \cite[Theorem 2]{MR1860517}.
\end{remark}
\begin{proofarg}{Proof of Lemma \ref{calculus of a}}
 Let us compute $a_{\bf 1}(k_1, k_2)$, for $k_1, k_2 \in \naturels$.
We have
 \begin{align*}
 a_{\bf 1}(k_1, k_2)= \tfrac{\ds \sqrt{\mu}}{\ds k_1! k_2!} \int_{\reels^2}
 \sqrt{\lambda_1^2\, y_1^2+\lambda_2^2 \,y_2^2}\,  \,H_{k_1}(y_1)\, H_{k_2}(y_2)\, \phi(y_1)\, \phi(y_2)\, \ud y_1\, \ud y_2.
 \end{align*}
 Similar arguments as those given in Lemma \ref{calculus of c} show that $a_{\bf 1}(k_1, k_2)=0$ except when $k_1$ and $k_2$ are even.\\
 So let us compute $a_{\bf 1}(2m, 2\ell)$ for $m, \ell \in \naturels$.\\
The change of variable $y_1:=\lambda r \sin(\theta)$ and $y_2:=r \cos(\theta)$ in last integral and expression of $H_{2\ell}$ given in (\ref{hermite form}) yield
\begin{multline*}
 a_{\bf 1}(2m, 2\ell)= \tfrac{\ds 2\sqrt{\ds \mu}}{\ds \pi} \tfrac{\ds \lambda_2^2}{\ds \lambda_1} \sum_{p=0}^{m} \sum_{q=0}^{\ell} \frac{(-2)^{p-m}(-2)^{q-\ell}}{(2p)!(2q)! (m-p)!(\ell-q)!} \,\lambda^{2p}\, \times\\
 \int_{0}^{\frac{\pi}{2}} \int_{0}^{+\infty} r^{2p+2q+2} \cos^{2q}(\theta) \sin^{2p}(\theta) 
\, e^{\ds -\frac{1}{2}r^2\left[{\cos^2(\theta)+\lambda^2 \sin^2(\theta)
 }\right]}\, \ud r \ud\theta.
  \end{multline*}
  Now making the change of variable $r^2:=\ds \frac{2v}{\cos^2(\theta)+\lambda^2 \sin^2(\theta)}$ and 
  $w:= \sin(\theta)$, one gets
  \begin{multline*}
 a_{\bf 1}(2m, 2\ell)= \tfrac{\ds \sqrt{\ds \mu}}{\ds \pi} \tfrac{\ds \lambda_2^2}{\ds \lambda_1} \sum_{p=0}^{m} \sum_{q=0}^{\ell} \frac{(-2)^{p-m}(-2)^{q-\ell}}{(2p)!(2q)! (m-p)!(\ell-q)!} \,\lambda^{2p}\, 2^{p+q+\tfrac{1}{2}}\,\times\\
 \mathit{\Gamma}(p+q+\tfrac{3}{2})\, F(p+q+\tfrac{3}{2}; p+\tfrac{1}{2}; p+q+1; 1-\lambda^2),
\end{multline*}
where the hypergeometric function $F$ is defined by
 \begin{align*}
 F(a; b; c; z):= \int_0^1 u^{b-1} (1-u)^{c-b-1} (1-uz)^{-a} \ud u.
 \end{align*}
 for $\abs{z} <1$, $0<b<c$ and $a>0$.\\
The proof of lemma ensues from the following one.
\begin{lemma}
\label{geom}
For $\abs{z} <1$, $0<b<c$ and $a>0$, one has
\begin{align*}
F(a; b; c; z)=\ds  \tfrac{\ds \mathit{\Gamma}(c-b)}{\ds \mathit{\Gamma}(a)} \sum_{n=0}^{+\infty} \tfrac{\ds \mathit{\Gamma}(a+n)\mathit{\Gamma}(b+n)}{\ds \mathit{\Gamma}(c+n)\mathit{\Gamma}(n+1)} \,z^n.
\end{align*}
\end{lemma}
\begin{proofarg}{Proof of Lemma \ref{geom}}
The proof consists in a serial development of function $f(z):=(1-uz)^{-a}$.
\end{proofarg}
\end{proofarg}
\begin{proofarg}{Proof of Lemma \ref{determinant}}
Let $f_1, f_2:S^1 \to \reels$ continuous and bounded functions and set $\overrightarrow{f}:=(f_1, f_2)$.
We have the following decomposition of $\det\left({\mathit{\Sigma}_{\overrightarrow{f}}(u)}\right)$.
\begin{align*}
\lefteqn{\det\left({\mathit{\Sigma}_{\overrightarrow{f}}(u)}\right)}\\
&=\mathit{\Sigma}_{f_1, f_1}(u) \mathit{\Sigma}_{f_2, f_2}(u) -\left({\mathit{\Sigma}_{f_1, f_2}(u)}\right)^2\\
&=\sum_{q=1}^{+\infty}  \sum_{q'=1}^{+\infty} \left({\mathit{\Sigma}_{f_1, f_1}(u)}\right)_q \left({\mathit{\Sigma}_{f_2, f_2}(u)}\right)_{q'}- \left({\sum_{q=1}^{+\infty} \left({\mathit{\Sigma}_{f_1, f_2}(u)}\right)_q}\right)^2\\
&=\sum_{q=1}^{+\infty} \det\left({\mathit{\Sigma}_{\overrightarrow{f}}(u)}\right)_{q}\\
 &+\sum_{q<q'} \left({
\sqrt{\left({\mathit{\Sigma}_{f_1, f_1}(u)}\right)_q}\,\sqrt{\left({\mathit{\Sigma}_{f_2, f_2}(u)}\right)_{q'}}-
\sqrt{\left({\mathit{\Sigma}_{f_1, f_1}(u)}\right)_{q'}}\,\sqrt{\left({\mathit{\Sigma}_{f_2, f_2}(u)}\right)_{q}}
 }\right)^2\\
 &+2  \sum_{q<q'} \left[{
 \sqrt{\left({\mathit{\Sigma}_{f_1, f_1}(u)}\right)_q}\,\sqrt{\left({\mathit{\Sigma}_{f_2, f_2}(u)}\right)_{q'}}
\sqrt{\left({\mathit{\Sigma}_{f_1, f_1}(u)}\right)_{q'}}\,\sqrt{\left({\mathit{\Sigma}_{f_2, f_2}(u)}\right)_{q}}  }\right.\\
&\qquad\qquad\left.{ -\left({\mathit{\Sigma}_{f_1, f_2}(u)}\right)_q\, \left({\mathit{\Sigma}_{f_1, f_2}(u)}\right)_{q'}
 }\right].
\end{align*}
To conclude the proof of lemma, we just have to verify that for all $q \in \naturels^{\star}$, 
\begin{align}
\label{inegalite sigma}
\abs{\left({\mathit{\Sigma}_{f_1, f_2}(u)}\right)_q} \le  \sqrt{\left({\mathit{\Sigma}_{f_1, f_1}(u)}\right)_q}\,\sqrt{\left({\mathit{\Sigma}_{f_2, f_2}(u)}\right)_{q}}.
\end{align}
So, let $q \in \naturels^{\star}$, and let define
\begin{align*}
\left({\xi_{f}^{(n)}(u)}\right)_q:=\tfrac{\ds 1}{\sqrt{\ds \sigma_2(T_n)}} \sum_{
\substack{
{\bk} \in \naturels^3 \\ \abs{{\bk}}=q
}
} a_{f}({\bk}, u)\, \int_{T_n} \widetilde{H}_{{\bk}}(U(t))\, \ud t.
\end{align*}
Corollary \ref{convergence triple} implies that,
\begin{align*}
\left({\left({\xi_{f_{1}}^{(n)}(u)}\right)_q\, ; \,\left({\xi_{f_{2}}^{(n)}(u)}\right)_q}\right)
\xrightarrow[n \to +\infty]{Law} \mathcal{N}(0; \left({\mathit{\Sigma}_{\overrightarrow{f}}(u)}\right)_q),
\end{align*}
and in force $\left({\mathit{\Sigma}_{\overrightarrow{f}}(u)}\right)_q$ is a semi-definite positive matrix.
This argument yields inequality (\ref{inegalite sigma}) and Lemma \ref{determinant}.
\end{proofarg}
\begin{proofarg}{Proof of Lemma \ref{convergence racine carree}}
Using that $B_n^2= \Sigma_n$, $\lim_n \Sigma_n=\Sigma=B^2$ and that $\lim_n \tr(B_n)= \tr(B) >0$, we obviously deduce that $\lim_n B_n=B$.\\
\end{proofarg}

\noindent {\bf Acknowledgements} \\

The author is grateful to the referees for many suggestions that led to substantial improvements in this work.

%
%
\providecommand{\bysame}{\leavevmode\hbox to3em{\hrulefill}\thinspace}
\providecommand{\MR}{\relax\ifhmode\unskip\space\fi MR }
\providecommand{\MRhref}[2]{%
  \href{http://www.ams.org/mathscinet-getitem?mr=#1}{#2}
}
\providecommand{\href}[2]{#2}




\end{document}